\documentclass[10pt]{amsart}

\usepackage{eucal,url,amssymb,stmaryrd,booktabs,
enumerate,amscd,paralist}
\usepackage[pagebackref,colorlinks=true]{hyperref}  

\usepackage{amsfonts}
\usepackage{amsmath,amsthm,amssymb,amscd,enumerate,eucal,url,stmaryrd}

\usepackage{eucal,url,amssymb,stmaryrd,
enumerate,amscd,
}

\numberwithin{equation}{section}

\newtheorem{thrm}{Theorem}[section]

\newtheorem{lemma}[thrm]{Lemma}

\newtheorem{prop}[thrm]{Proposition}

\newtheorem{cor}[thrm]{Corollary}

\newtheorem{dfn}[thrm]{Definition}

\newtheorem{rmrk}[thrm]{Remark}

\newtheorem{exm}[thrm]{Example}
\newtheorem{notation}[thrm]{Notation}

\newcommand{\Hn}{\mathbb{H}^n}
\newcommand{\Hnn}{\mathbb{H}^{n+1}}
\newcommand{\QH}{\boldsymbol {G\,(\mathbb{H})}}

\newcommand{\lap}{\mathcal{L}}

\newcommand{\abs}[1]{\lvert #1 \rvert}

\newcommand{\ds}{\frac{d}{ds}}

\newcommand{\diffta}{\,dt^\alpha}
\newcommand{\diffxa}{\,dx^\alpha}
\newcommand{\diffya}{\,dy^\alpha}
\newcommand{\diffza}{\,dz^\alpha}
\newcommand{\difftb}{\,dt^\beta}
\newcommand{\diffxb}{\,dx^\beta}
\newcommand{\diffyb}{\,dy^\beta}
\newcommand{\diffzb}{\,dz^\beta}

\newcommand{\dxa}[1]{{\partial_{x_{\alpha}}} {#1}}
\newcommand{\dta}[1]{{\partial_{t_{\alpha}}} {#1}}
\newcommand{\dya}[1]{{\partial_{y_{\alpha}}} {#1}}
\newcommand{\dza}[1]{{\partial_{z_{\alpha}}} {#1}}

\newcommand{\dxb}[1]{{\partial_{x_{\beta}}} {#1}}
\newcommand{\dtb}[1]{{\partial_{t_{\beta}}} {#1}}
\newcommand{\dyb}[1]{{\partial_{y_{\beta}}} {#1}}
\newcommand{\dzb}[1]{{\partial_{z_{\beta}}} {#1}}

\newcommand{\dt}[1]{{\partial_t} {#1}}
\newcommand{\dx}[1]{{\partial_x} {#1}}
\newcommand{\dy}[1]{{\partial_y} {#1}}
\newcommand{\dz}[1]{{\partial_z} {#1}}

\newcommand{\dtbta}[1]{{#1}_{ t_{\beta}{ t_{\alpha}}} }
\newcommand{\dtbxa}[1]{{#1}_{ t_{\beta}{ x_{\alpha}}} }
\newcommand{\dtbya}[1]{{#1}_{ t_{\beta}{ y_{\alpha}}} }
\newcommand{\dtbza}[1]{{#1}_{ t_{\beta}{ z_{\alpha}}} }

\newcommand{\dxbta}[1]{{#1}_{ x_{\beta}{ t_{\alpha}}} }
\newcommand{\dxbxa}[1]{{#1}_{ x_{\beta}{ x_{\alpha}}} }
\newcommand{\dxbya}[1]{{#1}_{ x_{\beta}{ y_{\alpha}}} }
\newcommand{\dxbza}[1]{{#1}_{ x_{\beta}{ z_{\alpha}}} }

\newcommand{\dybta}[1]{{#1}_{ y_{\beta}{ t_{\alpha}}} }
\newcommand{\dybxa}[1]{{#1}_{ y_{\beta}{ x_{\alpha}}} }
\newcommand{\dybya}[1]{{#1}_{ y_{\beta}{ y_{\alpha}}} }
\newcommand{\dybza}[1]{{#1}_{ y_{\beta}{ z_{\alpha}}} }

\newcommand{\dzbta}[1]{{#1}_{ z_{\beta}{ t_{\alpha}}} }
\newcommand{\dzbxa}[1]{{#1}_{ z_{\beta}{ x_{\alpha}}} }
\newcommand{\dzbya}[1]{{#1}_{ z_{\beta}{ y_{\alpha}}} }
\newcommand{\dzbza}[1]{{#1}_{ z_{\beta}{ z_{\alpha}}} }

\newcommand{\dbar}{\overline {\mathcal{D}}\ }
\newcommand{\dbara}{\overline {\mathcal{D}}_\alpha\ }
\newcommand{\dbarb}{\overline {\mathcal{D}}_\beta \ }

\newcommand{\dirb}{\mathcal{D}_\beta\ }
\newcommand{\dira}{\mathcal{D}_\alpha\ }

\newcommand{\diracb}[1]{( \dtb {#1}\ - \ i\dxb {#1}\ - \ j\dyb {#1}\ - \ k\dzb {#1} )}

\newcommand{\sameauthor}{{\leavevmode\vrule height 2pt depth -1.6pt width 23pt}}

\newcommand{\LieQ}{\mathcal{L}_Q\, }

\def\sideremark#1{\ifvmode\leavevmode\fi\vadjust{\vbox to0pt{\vss
 \hbox to 0pt{\hskip\hsize\hskip1em
 \vbox{\hsize2.5cm\tiny\raggedright\pretolerance10000
 \noindent #1\hfill}\hss}\vbox to8pt{\vfil}\vss}}}%

\renewcommand{\sideremark}[1]{}

\begin{document}
\begin{abstract}
A partial solution of the quaternionic contact Yamabe problem on the
quaternionic sphere is given. It is shown that the torsion of the
Biquard connection vanishes exactly when the trace-free part of the
horizontal Ricci tensor of the Biquard connection is zero and this
occurs precisely on 3-Sasakian manifolds. All conformal deformations
sending the standard flat torsion-free quaternionic contact
structure on the quaternionic Heisenberg group to a quaternionic
contact structure with vanishing torsion of the Biquard connection
are explicitly described. A '3-Hamiltonian form' of infinitesimal
conformal automorphisms of quaternionic contact structures is
presented.
\end{abstract}

\title[Quaternionic contact structures and the Yamabe problem]
{Quaternionic contact Einstein structures and the quaternionic contact Yamabe problem}
\date{\today}
\thanks{This project has been funded in part by the National Academy
of Sciences under the [Collaboration in Basic Science and Engineering Program 1 Twinning
Program] supported by Contract No. INT-0002341 from the National Science Foundation. The
contents of this publication do not necessarily reflect the views or policies of the
National Academy of Sciences or the National Science Foundation, nor does mention of
trade names, commercial products or organizations imply endorsement by the National
Academy of Sciences or the National Science Foundation.}

\author{Stefan Ivanov}
\address[Stefan Ivanov]{University of Sofia, Faculty of Mathematics and Informatics,
blvd. James Bourchier 5, 1164,
Sofia, Bulgaria} \email{ivanovsp@fmi.uni-sofia.bg}
\address{and Max-Planck-Institut f\"ur Mathematik\\Vivatsgasse
  7\\D-53111 Bonn\\Germany}

\author{Ivan Minchev}
\address[Ivan Minchev]{University of Sofia\\
Sofia, Bulgaria\\
and Institut f\"ur Mathematik,
Humboldt Universit\"at zu Berlin\\
Unter den Linden~6, Berlin~D-10099, Germany}
\email{minchevim@yahoo.com} \email{minchev@fmi.uni-sofia.bg}

\author{Dimiter Vassilev}
\address[Dimiter Vassilev]{
Department of Mathematics and Statistics\\ University of New Mexico\\
Albuquerque, New Mexico, 87131-0001\\
and\\
University of California, Riverside\\
Riverside, CA 92521} \email{vassilev@math.unm.edu} \maketitle


\setcounter{tocdepth}{2}

\tableofcontents

\section{Introduction}

The Riemannian \cite{LP} and CR Yamabe problems \cite{JL1,JL2,JL3,JL4} have been a
fruitful subject in geometry and analysis. Major steps in the  solutions is the
understanding of the conformally flat cases. A model for this setting is given by the
corresponding spheres, or equivalently, the Heisenberg groups with, respectively,
0-dimensional and 1-dimensional centers. The equivalence is established through the
Cayley transform \cite{Ko1}, \cite{CDKR} and \cite{CDKR2}, which in the Riemannian case
is the usual stereographic projection.

In the present paper we consider the Yamabe problem on the quaternionic Heisenberg group
(three dimensional center). This problem turns out to be equivalent to the quaternionic
contact Yamabe problem on the unit (4n+3)-dimensional sphere in the quaternionic space
due to the quaternionic Cayley transform, which is a conformal quaternionic contact
transformation (see the proof of Theorem~\ref{t:Yamabe}).

The central notion is the quaternionic contact structure (QC
structure for short), \cite{Biq1,Biq2}, which  appears naturally as
the conformal boundary at infinity of the quaternionic hyperbolic
space, see also \cite{P,GL,FG}. Namely, a QC structure
$(\eta,\mathbb Q)$ on a (4n+3)-dimensional smooth manifold $M$ is a
codimension 3 distribution $H$, such that, at each point $p\in M$
the nilpotent Lie algebra $H_p\oplus (T_pM/H_p)$ is isomorphic to
the quaternionic Heisenberg algebra $\mathbb{H}^m\oplus Im\
\mathbb{H}$. This is equivalent to the existence of a 1-form
$\eta=(\eta_1,\eta_2,\eta_3)$ with values in $\mathbb R^3$, such
that, $H=Ker~\eta$ and the three 2-forms ${d\eta_i}_{|H}$ are the
fundamental 2-forms of a quaternionic structure $\mathbb Q$ on $H$.
A special phenomena here, noted by Biquard \cite{Biq1}, is that the
3-contact form $\eta$ determines the quaternionic structure as well
as the metric on the horizontal bundle in a unique way. Of crucial
importance is the existence of a distinguished linear connection,
see \cite{Biq1}, preserving the QC structure
 and its Ricci tensor and scalar curvature $Scal$, defined in
\eqref{e:horizontal ricci}, and called correspondingly qc-Ricci tensor  and qc-scalar
curvature. The Biquard connection will play a role similar to the Tanaka-Webster
connection \cite{W} and \cite{T} in the CR case.

The quaternionic contact Yamabe problem, in the considered setting, is about the
possibility of finding in the conformal class of a given QC structure one with constant
qc-scalar curvature.

The question reduces to the solvability of the Yamabe equation \eqref{e:conf change
scalar curv}. As usual if we take the conformal factor in a suitable form the gradient
terms in \eqref{e:conf change scalar curv} can be removed and one obtains the more
familiar form of the Yamabe equation.  In fact, taking the conformal factor of the form
$\bar\eta=u^{1/(n+1)}\eta$ reduces \eqref{e:conf change scalar curv} to the equation
\[
\lap u\ \equiv\  4\frac {n+2}{n+1}\ \triangle u -\ u\, Scal \ =\ -
u^{2^*-1}\,\overline{Scal},
\]
where $\triangle $ is the horizontal sub-Laplacian, cf. \eqref{e:horizontal laplacian and
gradient}, and $Scal$ and $\overline{Scal}$ are the qc-scalar curvatures correspondingly
of $(M,\, \eta)$ and $(M, \, \bar\eta)$, and
$2^* = \frac {2Q}{Q-2},$ with $Q=4n+6$. In the case of the quaternionic Heisenberg
group, cf. Section 4.1,  the equation is
\begin{equation*}\label{Yamabe}
\lap u\ \equiv\ \sum_{\alpha = 1}^n \big ( T^2_\alpha u \ +\
X^2_\alpha u \ +\ Y^2_\alpha u \ +\ Z^2_\alpha u \big )\ = \ -\frac
{n+1}{4(n+2)}\, u^{2^*-1}\,\overline{Scal}.
\end{equation*}
\noindent  This is also, up to a scaling, the Euler-Lagrange equation of the non-negative
extremals in the $L^2$ Folland-Stein embedding theorem \cite{F2} and \cite{FS}, see
\cite{GV} and \cite{Va2}. On the other hand, on a compact quaternionic contact manifold
$M$ with a fixed conformal class $[\eta]$ the Yamabe equation characterizes the
non-negative extremals of the Yamabe functional defined by
$$ \Upsilon (u)\ =\
\int_M 4\frac {n+2}{n+1}\ \abs{\nabla u}^2\ +\ \text{Scal}\, u^2\
dv_g,\qquad \int_M u^{2^*}\, dv_g \ =\ 1, \ u>0. $$ When the Yamabe
constant $\lambda(M)\overset{def}=\lambda(M,[\eta])=\inf \{ \Upsilon
(u)\,:\ \int_M u^{2^*}\, dv_g \ =\ 1, \ u>0\} $ is less than that of
the sphere the existence of solutions can be constructed with the
use of suitable coordinates see \cite{Wei} and \cite{JL2}.

The present paper can be considered as a contribution towards the
soluiton of the Yamabe problem in the  case when the Yamabe constant
of the considered quaternionic contact manifold is equal to the
Yamabe constant of the unit sphere with its standard quaternionic
contact structure, which is induced from the embedding in the
quaternion $(n+1)$-dimensional space. It is also natural to
conjecture that if the quaternionic contact structure is not locally
equivalent to the standard sphere then the Yamabe constant is less
than that of the sphere, see \cite{JL4} for a proof in the CR case.
The results of the present paper will be instrumental for the
analysis of these and some other questions concerning the geometric
analysis on quaternionic contact structures.

In this article we provide a partial solution of the Yamabe problem on the quaternionic
sphere with its standard contact quaternionic structure or, equivalently, the
quaternionic Heisenberg group. Note that according to \cite{GV2} the extremals of the
above variational problem are $\mathcal{C}^\infty$ functions, so we will not consider
regularity questions in this paper. Furthermore, according to \cite{Va} or \cite{Va2} the
infimum is achieved, and the extremals are solutions of the Yamabe equation.  Let us
observe that \cite{GV2} solves the same problem in a more general setting, but under the
assumption that the solution is invariant under a certain group of rotation.  If one is
on the flat models, i.e., the groups of Iwasawa  type \cite{CDKR} the assumption in
\cite{GV2} is equivalent to the a-priori assumption that, up to a translation, the
solution is radial with respect to the variables in the first layer. The proof goes on by
using the moving plane method and showing that the solution is radial also in the
variables from the center, after which a very non-trivial identity is used to determine
all cylindrical solutions. In this paper the a-priori assumption is of a different
nature, see further below, and the method has the potential of solving the general
problem.

The strategy, following the steps of \cite{LP} and \cite{JL3}, is to solve the Yamabe
problem on the quaternionic sphere by replacing the non-linear Yamabe equation by an
appropriate geometrical system of equations which could be solved.

Our first observation is that if if $n>1$ and the qc-Ricci tensor is trace-free (qc-Einstein
condition) then the qc-scalar curvature is constant 
(Theorem~\ref{maincon}). Studying conformal deformations of QC structures preserving the
qc-Einstein condition, we describe explicitly all global functions on the quaternionic
Heisenberg group sending conformally the standard flat QC structure to another
qc-Einstein structure. Our result here is the following Theorem.

\begin{thrm}\label{t:einstein preserving}
Let $\Theta=\frac{1}{2h}\tilde\Theta$ be a conformal deformation of the standard
qc-structure $\tilde\Theta$ on the quaternionic Heisenberg group $\QH$. If $\Theta$ is
also qc-Einstein, then up to a left translation the function $h$ is given by
$$h \ =\ c\ \Big [  \big ( 1\ +\
\nu\, |q|^2 \big )^2\  +\ \nu^2\, (x^2\ +\ y^2\ +\ z^2)\Big ],$$ where $c$ and $\nu$ are
positive constants. All functions $h$ of this form have this property.
\end{thrm}

The crucial observation reducing the Yamabe  equation to the system preserving the
qc-Einstein condition is Proposition~\ref{qcYamab} which asserts that, under some "extra"
conditions, QC structure with constant qc-scalar curvature obtained by a conformal
transformation of a qc-Einstein structure on compact manifold must be again qc-Einstein.
The prove of this relies on detailed analysis of the Bianchi identities for the Biquard
connection. Using the quaternionic Cayley transform combined with Theorem~\ref{t:einstein
preserving} lead to our second main  result.

\begin{thrm}\label{t:Yamabe}
Let $\eta= f\,\tilde\eta$ be a conformal deformation of the standard
qc-structure $\tilde\eta$ on the quaternionic  sphere $S^{4n+3}$. Suppose
$\eta$ has constant qc-scalar curvature.
\begin{enumerate}
\item[a)] { If $n>1$ then any one of the following two conditions}
\begin{enumerate}
\item[i)] the
vertical space of $\eta$ is integrable, \item[ii)] the function  $f$ is
the real part of an anti-CRF function,
\end{enumerate} implies that  up to a multiplicative constant $\eta$ is obtained from
$\tilde\eta$ by a conformal  quaternionic contact automorphism.
\item[b)] { If $n=1$ and the
vertical space of $\eta$ is integrable then  up to a multiplicative
constant $\eta$ is obtained from $\tilde\eta$ by a  conformal
quaternionic contact automorphism.}
\end{enumerate}
\end{thrm}

The definition of  conformal quaternionic contact automorpism can be
found in Definition \ref{d:3-ctct auto}.  The solutions (conformal
factors) we find agree with those conjectured in \cite{GV}. It might
be possible to dispense of the "extra" assumption in
Theorem~\ref{t:Yamabe}. In a subsequent paper \cite{IMV1} we give
such a proof for the seven dimensional sphere.

Studying the geometry of the Biquard connection, our main geometrical tool towards
understanding the geometry of the Yamabe equation, we show that the qc-Einstein condition
is equivalent to the vanishing of the torsion of Biquard connection. Furthermore, in our
third main result we give a local characterization of such spaces.

\begin{thrm}\label{Ein2MO}
{Let $(M^{4n+3},g,\mathbb{Q})$ be a QC manifold}  with positive qc
scalar curvature $Scal>0$. The next conditions are equivalent:
\begin{enumerate}[a)]
\item {$(M^{4n+3},g,\mathbb{Q})$ is qc-Einstein manifold with
constant qc scalar curvature if $n=1$}; \item $M$ is locally
3-Sasakian in the sense that locally there exists a $SO(3)$-matrix
$\Psi$ with smooth entries, such that, the local QC structure
$(\frac{16n(n+2)}{Scal}\Psi\cdot\eta, \mathbb{Q})$ is 3-Sasakian; \item The
torsion of the Biquard connection is identically zero and the qc
scalar curvature is constant if $n=1$.
\end{enumerate}
In particular, a qc-Einstein manifold  with positive qc scalar curvature, assumed also constant if
$n=1$, is Einstein manifold with positive
Riemannian scalar curvature and if complete it is compact with
finite fundamental group.
\end{thrm}
In addition to the above results, in Theorem \ref{flatqc} we show that the above conditions are equivalent to
the property that every Reeb vector field, defined in \eqref{bi1}, is an infinitesimal
generator of a conformal quaternionic contact automorphism, cf. Definition \ref{d:3-ctct
v field}. It might be helpful to point that  the statement of the above Theorem reflects the outstanding open question whether a seven dimensional qc-Einstein structure must have a constant qc scalar curvature, cf. Theorem \ref{maincon} and Remark \ref{r:open question}.

In the paper we also develop useful tools necessary for the geometry and analysis on QC
manifolds. We define and study some special functions, which will be relevant in the
geometric analysis on quaternionic contact and hypercomplex manifolds as well as
properties of infinitesimal automorphisms of QC structures. In particular, the considered
anti-regular functions will be relevant in the study of qc-pseudo-Einstein structures,
cf. Definition \ref{d:qc-pseudo_einstein}.

\textbf{Organization of the paper}: In the following two chapters we describe in details
the notion of a quaternionic contact manifold, abbreviate sometimes to QC-manifold, and
the Biquard connection, which is central to the paper.

In Chapter 4 we write explicitly the Bianchi identities and derive a system of equations
satisfied by the divergences of some important tensors. As a result we are able to show
that qc-Einstein manifolds, i.e., manifolds for which the restriction to the horizontal
space of the qc-Ricci tensor is proportional to the metric, have constant scalar
curvature, see Theorem \ref{maincon}. The proof uses Theorem \ref{t:horizontal system} in
which we derive a relation between the horizontal divergences of certain
$Sp(n)Sp(1)$-invariant tensors. By introducing an integrability condition on the
horizontal bundle we define hyperhermitian contact structures, see Definition \ref{d:hyper
cplx ctct}, and with the help of Theorem \ref{t:horizontal system} we prove
Theorem~\ref{Ein2MO}.

Chapter 5 describes the conformal transformations preserving the qc-Einstein condition.
Note that here a conformal quaternionic contact transformation between two quaternionic
contact manifold is a diffeomorphism $\Phi$ which satisfies
$\Phi^*\eta=\mu\ \Psi\cdot\eta,$
for some positive smooth function $\mu$ and some matrix $\Psi\in SO(3)$ with smooth
functions as entries and $\eta=(\eta_1,\eta_2,\eta_3)^t$ is considered as an element of
$\mathbb R^3$. One defines in an obvious manner a point-wise conformal transformation.
Let us note that the Biquard connection does not change under rotations as above, i.e.,
the Biquard connection of $\Psi\cdot\eta$ and $\eta$ coincides. In particular, when
studying conformal transformations we can consider only transformations with $\Phi^*\eta\
=\ \mu\ \eta$. We find all conformal transformations  preserving the qc-Einstein
condition on the quaternionis Heisenberg group or, equivalently, on the quaternionic
sphere with their standard contact quaternionic structures proving  Theorem
\ref{t:einstein preserving}.

Chapter 6 concerns a special class of functions, which we call anti-regular, defined
respectively on the quaternionic space, real hyper-surface in it, or on a quaternionic
contact manifold, cf. Definitions \ref{d:anti-regular functions} and \ref{d:anti-CRF
functions} as functions preserving the quaternionic structure. The anti-regular functions
play a role somewhat similar to those played by the CR functions, but the analogy is not
complete. The real parts of such functions will be also of interest in connection with
conformal transformation preserving the qc-Einstein tensor and should be thought of as
generalization of pluriharmonic functions. Let us stress explicitly that regular
quaternionic functions have been studied extensively, see \cite{S} and many subsequent
papers, but they are not as relevant for the considered geometrical structures.
Anti-regular functions on hyperk\"ahler and quaternionic K\"ahler manifolds are studied
in \cite{CL1,CL2,LZ} in a different context, namely in connection with minimal surfaces
and quaternionic maps between quaternionic K\"ahler manifolds. The notion of hypercomplex
contact structures will appear in this section again since on such manifolds the real
part of anti-CRF functions, see \eqref{crfhyp} for the definition, have some interesting
properties, cf. Theorem \ref{crfth1}

In Chapter 7 we study infinitesimal  conformal automorphisms of QC
structures (QC-vector fields) and show that they depend on three
functions satisfying some differential conditions thus establishing
a '3-hamiltonian' form of the QC-vector fields
(Proposition~\ref{qaut1}). The formula becomes very simple
expression on a 3-Sasakian manifolds (Corollary~\ref{connforms}). We
characterize the vanishing of the torsion of Biquard connection in
terms of the existence of three vertical vector fields whose flow
preserves the metric and the quaternionic structure. Among them,
3-Sasakian manifolds are  exactly those admitting three transversal
QC-vector fields (Theorem~\ref{flatqc}, Corollary~\ref{3con3sas}).

In the last section we complete the proof of our main result Theorem~\ref{t:Yamabe}.
\begin{notation}\label{n:notation}\begin{enumerate}[a)]
\item Let us note explicitly, that  in this paper for a one form $\theta$ we use
\[
d\theta(X,Y)\ =\ X\theta(Y)\ -\ Y\theta(X)\ -\ \theta([X,Y]).
\]
\item We shall denote with $\nabla h$ the horizontal gradient of the function $h$, see \eqref{New19},
    while $dh$ means as usual the differential of the function $h$.
\item  The triple $\{i,j,k\}$ will denote a cyclic permutation of  $\{1,2,3\}$, unless it is explicitly  stated otherwise.
\end{enumerate}
\end{notation}

\textbf{Acknowledgements} S.Ivanov is a Senior Associate to the Abdus Salam ICTP.
The paper was completed during the visit of S.I. in Max-Plank-Institut f\"ur Mathematik, Bonn.
S.I. thanks ICTP and MPIM, Bonn for providing the support and an excellent research environment.
I.Minchev is a member of the Junior Research Group "Special
Geometries in Mathematical Physics" founded by the Volkswagen Foundation" The authors
would like to thank The National Academies for the support.  It is  a pleasure to
acknowledge the role of Centre de Recherches Math\'ematiques and CIRGET, Montr\'eal where
the project initiated. The authors also would  like to thank University of California,
Riverside and University of Sofia for hosting the respective visits of the authors.

\section{Quaternionic contact structures and the Biquard connection}
The notion of \emph{Quaternionic Contact Structure} has been introduced by O.Biquard in
~{\cite{Biq1} and \cite{Biq2}}. Namely, a quaternionic contact structure (QC structure
for short) on a (4n+3)-dimensional smooth manifold $M$ is a codimension 3 distribution
$H$, such that, at each point $p\in M$ the nilpotent step two Lie algebra $H_p\oplus
(T_pM/H_p)$ is isomorphic to the quaternionic Heisenberg algebra $\mathbb{H}^{n}\oplus\,
Im\ \mathbb{H}$.
The quaternionic Heisenberg algebra structure on
$\mathbb{H}^{n}\oplus\, Im\  \mathbb{H}$ is obtained by the identification of
$\mathbb{H}^{n}\oplus\, Im\  \mathbb{H}$ with the algebra of the left invariant vector
fields on the quaternionic Heisenberg group, see Section 5.2. In particular, the Lie
bracket is given by the formula
$\left[
(q_o, \omega_o), (q, \omega)\right]\ =\ 2\ Im \ q_o\cdot \bar q,$
where $q=(q^1, q^2,\dots,q^n),\ q_o=(q_o^1, q_o^2,\dots,q_o^n)\in \mathbb{H}^{n}$ and
$\omega, \ \omega_o\in Im\ \mathbb{H}$ with
$ q_o\cdot \bar q\ =\ \sum_{\alpha =1}^n q^\alpha_o\cdot \overline{ q^\alpha},$
see Section \ref{ss:Pluriharmonic functions in Hn} for notations concerning $\mathbb{H}$.
It is important to observe that if $M$ has a quaternionic contact structure as above then
the definition implies that the distribution $H$ and its commutators generate the tangent
space at every point.

The following is another, more explicit, definition of a quaternionic contact structure.
\begin{dfn}\cite{Biq1}
{A quaternionic contact structure (\emph{QC-structure}) on a $4n+3$ dimensional manifold
$M$, $n>1$, is the data of a codimension three distribution $H$ on $M$ equipped with a
CSp(n)Sp(1) structure, i.e., we have
\begin{enumerate}[i)]
\item a fixed conformal class $[g]$ of metrics on $H$;
\item a 2-sphere bundle $\mathbb{Q}$ over $M$ of almost complex structures, such that,
locally we have $\mathbb{Q}= \{aI_1+bI_2+cI_3:\ a^2+b^2+c^2=1  \}$, where the almost
complex structures $I_s\,:H \rightarrow H,\quad I_s^2\ =\ -1, \quad s\ =\ 1,\ 2,\ 3,$
satisfy the commutation relations of the imaginary quaternions  $I_1I_2=-I_2I_1=I_3$;
\item  $H$ is locally the kernel of a 1-form $\eta=(\eta_1,\eta_2,\eta_3)$ with values
in $\mathbb R^3$ and  the following compatibility condition holds
   \begin{equation}\label{con1}
    2g(I_sX,Y)\ =\ \ d\eta_s(X,Y), \quad s=1,2,3, \quad X,Y\in H.
   \end{equation}
\end{enumerate}}
\end{dfn}

A manifold $M$ with a structure as above will be called also quaternionic contact
manifold (QC manifold) and denoted by $(M, [g], \mathbb{Q})$. With a slight abuse of notation we shall use the letter $\mathbb{Q}$ to also  denote the rank-three bundle consisting of
endomorphisms of $H$ locally generated
by three almost complex structures $I_1,I_2,I_3$ on $H$.  The meaning should be clear from the context.  We note that if in some
local chart $\bar\eta$ is another form, with corresponding $\bar g\in [g]$ and almost
complex structures $\bar I_s$, $s=1,2,3$, then $\bar\eta\ =\ \mu\, \Psi\,\eta$ for some
$\Psi\in SO(3)$ and a positive function $\mu$
{Typical
examples of manifolds with QC-structures are totally umbilical hypersurfaces in
quaternionic K\"ahler or hyperk\"ahler manifold, see Proposition \ref{p:QRhypersurface}
for the latter.}

It is instructive to consider the case when there is a globally defined one-form $\eta$.
The obstruction to the global existence of $\eta$ is encoded in the first Pontrjagin
class \cite{AK}. Besides clarifying the notion of a QC-manifold, most of the time, for
example when considering the Yamabe equation, we shall work with a QC-structure for which
we have a fixed globally defined contact form. In this case, if we rotate the
$\mathbb{R}^3$-valued contact form and the almost complex structures by the same rotation
we obtain again a contact form, almost complex structures and a metric (the latter is
unchanged) satisfying the above conditions. On the other hand, it is important to observe
that given a contact form the almost complex structures and the horizontal metric are
unique if they exist. Finally, if we are given the horizontal bundle and a metric on it,
there exists at most one sphere of associated contact forms with a corresponding sphere
$\mathbb{Q}$ of almost complex structures \cite{Biq1}.

Besides the non-uniqueness due to the action of $SO(3)$, the 1-form $\eta$  can be
changed by a conformal factor, in the sense that if $\eta$ is a form for which we can
find associated almost complex structures and metric $g$ as above, then for any $\Psi\in
SO(3)$ and a positive function $\mu$, the form $\mu\, \Psi\,\eta$ also has an associated
complex structures and metric. In particular, when $\mu=1$ we obtain a whole unit sphere
of contact forms, and we shall denote, as already mentioned, by $\mathbb{Q}$ the
corresponding sphere bundle of associated triples of almost complex structures.  With the
above consideration in mind we introduce the following notation.

\begin{notation} {We shall denote with $(M, \eta)$ a
QC-manifold with a fixed globally defined contact form. $(M, g, \mathbb{Q})$ will denote
a QC-manifold with a fixed metric $g$ and a sphere bundle of almost complex structures
$\mathbb{Q}$. In this  case we have in fact a $Sp(n)Sp(1)$ structure, i.e., we are
working with a fixed metric on the horizontal space. Correspondingly, we shall denote
with $\eta$ any (locally defined) associated contact form.} \end{notation}

We recall the definition of the Lie groups $Sp(n)$, $Sp(1)$ and
$Sp(n)Sp(1)$.
Let us identify $\mathbb H^n=\mathbb R^{4n}$ and let
$\mathbb H$ acts on $\mathbb H^n$ by right multiplications, $\lambda(q)(W) = W\cdot
q^{-1}$. This defines a homomorphism $
 \lambda:\{ \text {unit  quaternions} \} \longrightarrow SO(4n)$
 with the convention that
$SO(4n)$ acts on $\mathbb R^{4n}$ on the left. The image is the Lie group $Sp(1)$. Let
$\lambda(i)=I_0, \lambda(j)=J_0,\lambda(k)=K_0$. The Lie algebra of $Sp(1)$ is
$sp(1)=span\{I_0,J_0,K_0\}.$ The group $Sp(n)$ is $Sp(n)=\{O\in SO(4n): OB=BO\quad {\text for
\quad \text{all}}\quad B\in Sp(1)\}$ or $Sp(n)=\{O\in M_n(\mathbb{H}): O\,\bar O ^t =
I\}$, and $O\in Sp(n)$ acts by $(q^1, q^2,\dots,q^n)^t\mapsto \, O \, (q^1,
q^2,\dots,q^n)^t $. Denote by $Sp(n)Sp(1)$ the product of the two groups in $SO(4n)$.
Abstractly, $Sp(n)Sp(1)=(Sp(n)\times Sp(1))/\mathbb Z_2$. The Lie algebra of the group
$Sp(n)Sp(1)$ is $sp(n)\oplus sp(1)$.

Any endomorphism $\Psi$ of $H$ can be  decomposed with respect to the quaternionic
structure $(\mathbb{Q},g)$ uniquely into four $Sp(n)$-invariant parts
$\Psi=\Psi^{+++}+\Psi^{+--}+\Psi^{-+-}+\Psi^{--+},$
where $\Psi^{+++}$ commutes with all three $I_i$, $\Psi^{+--}$ commutes with $I_1$ and
anti-commutes with the others two and etc. Explicitly,
\begin{gather*}4\Psi^{+++}=\Psi-I_1\Psi I_1-I_2\Psi I_2-I_3\Psi
I_3,\quad
4\Psi^{+--}=\Psi-I_1\Psi I_1+I_2\Psi I_2+I_3\Psi I_3,\\
4\Psi^{-+-}=\Psi+I_1\Psi I_1-I_2\Psi I_2+I_3\Psi I_3,\quad 4\Psi^{--+}=\Psi+I_1\Psi
I_1+I_2\Psi I_2-I_3\Psi I_3.
\end{gather*}
\noindent The two $Sp(n)Sp(1)$-invariant components are given by
\begin{equation}{\label{New21}}
\Psi_{[3]}=\Psi^{+++}, \qquad \Psi_{[-1]}=\Psi^{+--}+\Psi^{-+-}+\Psi^{--+}.
\end{equation}
\noindent Denoting the corresponding (0,2) tensor via $g$ by the same letter one sees
that the $Sp(n)Sp(1)$-invariant components are the projections on the eigenspaces of the
Casimir operator
\begin{equation}\label{e:cross}
\dag\ =\ I_1\otimes I_1\ +\ I_2\otimes I_2\ +\ I_3\otimes I_3
\end{equation}
corresponding, respectively, to the eigenvalues $3$ and $-1$, see \cite{CSal}. If $n=1$
then the space of symmetric endomorphisms commuting with all $I_i, i=1,2,3$ is
1-dimensional, i.e. the [3]-component of any symmetric endomorphism $\Psi$ on $H$ is
proportional to the identity, $\Psi_{3}=\frac{|\Psi|^2}{4}Id_{|H}$.

There exists a canonical connection compatible with a given quaternionic contact
structure. This connection was discovered by O. Biquard \cite{Biq1} when the dimension
$(4n+3)>7$ and by D. Duchemin \cite{D} in the 7-dimensional case.  The next result due to
O. Biquard is crucial in the quaternionic contact geometry.

\begin{thrm}\cite{Biq1}\label{biqcon}
{Let $(M, g,\mathbb Q)$ be a quaternionic contact manifold} of
dimension  $4n+3>7$ and a fixed metric $g$ on $H$ in the conformal class $[g]$. Then
there exists a unique connection $\nabla$ with torsion $T$ on $M^{4n+3}$ and a unique
supplementary subspace $V$ to $H$ in $TM$, such that:
\begin{enumerate}[i)]
\item $\nabla$ preserves the decomposition $H\oplus V$ and
the metric $g$;
\item for $X,Y\in H$, one has $T(X,Y)=-[X,Y]_{|V}$;
\item
$\nabla$ preserves the $Sp(n)Sp(1)$-structure on $H$, i.e., $\nabla g\ = \ 0$ and $\nabla
\mathbb{Q}\subset\mathbb{Q}$;
\item for $\xi\in V$, the endomorphism $T(\xi,.)_{|H}$
of $H$ lies in $(sp(n)\oplus sp(1))^{\bot}\subset gl(4n)$;
\item the connection on $V$ is induced by the natural identification $\varphi$ of
$V$ with the subspace $sp(1)$ of the endomorphisms of $H$, i.e. $\nabla\varphi=0$.
\end{enumerate}
\end{thrm}
\noindent In  (iv) the inner product on $End(H)$ is given by
\begin{equation}\label{e:inner product on End}
g(A,B)=tr(B^*A)=\sum_{a=1}^{4n}g(A(e_a),B(e_a)),
\end{equation}
where $A,B\in End(H)$, $\{e_1,...,e_{4n}\}$ is some $g$-orthonormal
basis of $H$.

We shall call the above connection \emph{the Biquard connection}.
Biquard \cite{Biq1} also described the supplementary
subspace $V$ explicitly, namely, {locally }$V$ is generated by vector fields
$\{\xi_1,\xi_2,\xi_3\}$, such that
\begin{equation}\label{bi1}
\begin{aligned}
\eta_s(\xi_k)=\delta_{sk}, \qquad (\xi_s\lrcorner d\eta_s)_{|H}=0,\\
(\xi_s\lrcorner d\eta_k)_{|H}=-(\xi_k\lrcorner d\eta_s)_{|H}, \qquad
s,\, k\in\{1,2,3\}.
\end{aligned}
\end{equation}
The vector fields $\xi_1,\xi_2,\xi_3$ are called Reeb vector fields or fundamental vector fields.

{ If the dimension of $M$ is seven, the conditions \eqref{bi1} do not always hold.
Duchemin shows in \cite{D} that if we assume, in addition, the existence of Reeb vector
fields as in \eqref{bi1}, then Theorem~\ref{biqcon} holds.  Henceforth, by a QC structure
in dimension $7$ we shall always mean a QC structure satisfying \eqref{bi1}.}

Notice that equations \eqref{bi1} are invariant under the natural $SO(3)$ action. Using
the  Reeb vector fields we extend  $g$ to a metric on $M$ by 
$span\{\xi_1,\xi_2,\xi_3\}=V\perp H \text{ and } g(\xi_s,\xi_k)=\delta_{sk}.$

\noindent The extended  metric does not depend on the action of $SO(3)$ on $V$, but it
changes in an obvious manner if $\eta$ is multiplied by a conformal factor. Clearly, the
Biquard connection preserves the extended metric on $TM, \nabla g=0$. We shall also
extend the quternionic structure by setting $I_{s|V}=0$.

Suppose the Reeb vector fields $\{\xi_1,\xi_2,\xi_3\}$ have been
fixed. The restriction of the torsion of the Biquard connection to
the vertical space $V$ satisfies \cite{Biq1}
\begin{equation}\label{torv}
\sideremark{torv} T(\xi_i,\xi_j)=\lambda \xi_k-[\xi_i,\xi_j]_{|_H},
\end{equation}
where $\lambda$ is a smooth function on $M$. Let us recall that here
and further in the paper the indices follow  the convention
\ref{n:notation}. Further properties of the Biquard connection are
encoded in the properties of the torsion endomorphism
$T_{\xi}=T(\xi,.) : H\rightarrow H, \quad \xi\in V$. Decomposing the
endomorphism $T_{\xi}\in(sp(n)+sp(1))^{\perp}$ into symmetric part
$T^0_{\xi}$ and skew-symmetric part $b_{\xi}$,
$T_{\xi}=T^0_{\xi} + b_{\xi},$
we summarize the description of the torsion due to O. Biquard in the
following Proposition.

\begin{prop}\cite{Biq1}\label{torb}
\sideremark{torb} The torsion $T_{\xi}$ is completely trace-free,
\begin{equation}\label{torb0}
\sideremark{torb0} tr T_{\xi}=\sum_{a=1}^{4n}g(T_{\xi}(e_a),e_a)=0, \quad tr T_{\xi}\circ
I=\sum_{a=1}^{4n}g(T_{\xi}(e_a),Ie_a)=0, \quad I\in Q,
\end{equation}
where $e_1\dots e_{4n}$ is an orthonormal basis of $H$.

The symmetric part of the torsion has the properties: \sideremark{tors1\\tors2}
\begin{gather}\label{tors1}
T^0_{\xi_i}I_i=-I_iT^0_{\xi_i}, \qquad i=1,2,3;\\\begin{aligned}\label{tors2}
I_2(T^0_{\xi_2})^{+--}=I_1(T^0_{\xi_1})^{-+-},\quad
I_3(T^0_{\xi_3})^{-+-}=I_2(T^0_{\xi_2})^{--+},\quad
I_1(T^0_{\xi_1})^{--+}=I_3(T^0_{\xi_3})^{+--}.
\end{aligned}
\end{gather}
The skew-symmetric part can be represented in the following way
\begin{equation}\label{toras}
\sideremark{toras} b_{\xi_i}=I_iu, \quad i=1,2,3,
\end{equation}
where $u$ is a traceless symmetric (1,1)-tensor on $H$ which commutes with $I_1,I_2,I_3$.

If $n=1$ then the tensor $u$ vanishes identically, $u=0$ and the torsion is a symmetric
tensor, $T_{\xi}=T^0_{\xi}$.
\end{prop}
\section{The torsion and curvature of the Biquard connection}

Let $(M^{4n+3},g,\mathbb Q)$ be a quaternionic contact structure on
a $4n+3$-dimensional smooth manifold.  Working in a local chart, we
 fix the vertical space $V=span\{\xi_1,\xi_2,\xi_3\}$ by requiring the conditions
\eqref{bi1}. The fundamental 2-forms $\omega_i, i=1,2,3$ of the
quaternionic structure $Q$ are defined by
\begin{equation}\label{thirteen}
2\omega_{i|H}\ =\ d\eta_{i|H},\qquad
\xi\lrcorner\omega_i=0,\quad \xi\in V.
\end{equation}
Define three 2-forms $\theta_i,
i=1,2,3$ by the formulas
\begin{multline}\label{fiftyone}
\theta_i=\frac12\{d((\xi_j\lrcorner d\eta_k)_{|_H})+(\xi_i\lrcorner
d\eta_j)\wedge (\xi_i\lrcorner d\eta_k)\}_{|_H}\\ =\
\frac12\{d(\xi_j\lrcorner d\eta_k)+(\xi_i\lrcorner d\eta_j)\wedge
(\xi_i\lrcorner
d\eta_k)\}_{|_H}-d\eta_k(\xi_j,\xi_k)\omega_k+d\eta_k(\xi_i,\xi_j)\omega_i.
\end{multline}
Define, further, the corresponding
$(1,1)$ tensors  $A_i$ by
$g(A_i(X),Y)=\theta_i(X,Y), X,Y\in H.$

\subsection{The torsion}
Due to \eqref{thirteen}, the torsion restricted to $H$ has the form
\begin{equation}\label{torha}
\sideremark{torha} T(X,Y)=-[X,Y]_{|V}=2\sum_{s=1}^3\omega_s(X,Y)\xi_s, \qquad X,Y\in H.
\end{equation}
The next two Lemmas provide some useful technical facts.
\begin{lemma}\label{twenty}\sideremark{twenty}
Let $D$ be any differentiation of the tensor algebra of $H$. Then we have
\begin{gather*}D(I_i)\,I_i\ =\ -I_i\,D(I_i),\ i=1,2,3,\\ \begin{aligned}
I_1\,D(I_1)^{-+-}\ =\ I_2\,D(I_2)^{+--},\quad I_1\,D(I_1)^{-+-}\ =\
I_2\,D(I_2)^{+--},\quad I_1\,D(I_1)^{-+-}\ =\ I_2\,D(I_2)^{+--}.
\end{aligned}
\end{gather*}
\end{lemma}
\begin{proof}The proof is a straightforward consequence of the next identities
\begin{gather*}
0=I_2(D(I_1)-I_2D(I_1)I_2)+I_1(D(I_2)-I_1D(I_2)I_1)
=I_2D(I_1)^{-+-}+I_1D(I_2)^{+--},\\
0=D(-\text{Id}_V)=D(I_iI_i)=D(I_i)I_i+I_iD(I_i).
\end{gather*}
\end{proof}

With $\mathcal{L}$ denoting the Lie derivative, we set
$\mathcal{L}'=\mathcal{L}_{|_H}$.

\begin{lemma}\label{twentytwo}
The following identities hold true.
\begin{gather}\label{thirtyone}
 \mathcal{L}^{'}_{\xi_1}I_1\ =\
-2T^0_{\xi_1}I_1+d\eta_1(\xi_1,\xi_2)I_2+d\eta_1(\xi_1,\xi_3)I_3,\\\label{twentytree}
 \mathcal{L}^{'}_{\xi_1}I_2\ =\
-2{T^0_{\xi_1}}^{--+}I_2-2I_3\tilde{u}+d\eta_1(\xi_2,\xi_1)I_1\\\nonumber \hspace{3cm} +\
{1\over
2}(d\eta_1(\xi_2,\xi_3)-d\eta_2(\xi_3,\xi_1)-d\eta_3(\xi_1,\xi_2))I_3\\\label{twentyfour}
 \mathcal{L}^{'}_{\xi_2}I_1\ =\
 -2{T^0_{\xi_2}}^{--+}I_1+2I_3\tilde{u}+d\eta_2(\xi_1,\xi_2)I_2\\\nonumber \hspace{3cm}
 -\ {1\over
2}(-d\eta_1(\xi_2,\xi_3)+d\eta_2(\xi_3,\xi_1)-d\eta_3(\xi_1,\xi_2))I_3,
\end{gather}
where the symmetric endomorphism $\tilde u$ on $H$, commuting with
$I_1,I_2,I_3$, is defined by
\begin{gather}\label{thirq}
2\tilde{u}=I_3(({\mathcal{L}^{'}_{\xi_1}I_2})^{--+})+
{\frac{1}{2}}(d\eta_1(\xi_2,\xi_3)-d\eta_2(\xi_3,\xi_1)-d\eta_3(\xi_1,\xi_2))Id_H
\end{gather}
\noindent  In addition, we have six more identities, which can be
obtained with a cyclic permutation of (1,2,3).
\end{lemma}
\begin{proof}
\noindent For all $k,l=1,2,3$ we have
\begin{equation}\label{eighteen}
\mathcal{L}_{\xi_k}\omega_l(X,Y)=\mathcal{L}_{\xi_k}g(I_lX,Y)+g((\mathcal{L}_{\xi_k}I_l)X,Y)
\end{equation}
\noindent Cartan's formula yields
\begin{equation}\label{New1}
\mathcal{L}_{\xi_k}\omega_l={\xi_k}\lrcorner (d\omega_l)+d({\xi_k}\lrcorner \omega_l).
\end{equation}
\noindent A direct calculation using \eqref{thirteen} gives
\begin{equation}\label{fifteen}\sideremark{fifteen}
2\omega_l=(d\eta_l)_{|_H}=d\eta_l - \sum_{s=1}^3{\eta_s\wedge({\xi_s}\lrcorner d\eta_{l})} +
\sum_{1\leq s<t\leq 3}{d\eta_l(\xi_s,\xi_t) \eta_s \wedge \eta_t}.
\end{equation}
\noindent Combining \eqref{fifteen} and \eqref{New1} we obtain,
after a short calculation, the following identities
\begin{gather}
\label{New2}({\mathcal{L}_{\xi_1}\omega_1})_{|H}=(d\eta_1(\xi_1,\xi_2)\omega_2+d\eta_1(\xi_1,\xi_3)\omega_3)_{|H} \\
 \label{fifty}
2{(\mathcal{L}_{\xi_i}\omega_j)}_{|H}={(d({\xi_i}\lrcorner
d\eta_j)-({\xi_i}\lrcorner d\eta_k)\wedge ({\xi_k}\lrcorner d\eta_j)
)}_{|H},
\end{gather} where $i\not=j\not=k\not=i, \quad i,j,k\in\{1,2,3\}$.
Clearly, \eqref{New2} and \eqref{eighteen} imply \eqref{thirtyone}.

Furthermore, using  \eqref{bi1} and  \eqref{fifty} twice for
$i=1,j=2$ and $i=2,j=1$, we find
\begin{multline}
\label{twentyfive} {(\mathcal{L}_{\xi_1}\omega_2\ + \ \mathcal{L}_{\xi_2}\omega_1)}_{|H}\
=\ {\frac{1}{2}}{(d({\xi_1}\lrcorner d\eta_2)\ +\ d({\xi_2}\lrcorner
d\eta_1))}_{|H} \\
=\ d\eta_1(\xi_2,\xi_1)\omega_1\ +\ d\eta_2(\xi_1,\xi_2)\omega_2\ +\
(d\eta_1(\xi_2,\xi_3)\ +\ d\eta_2(\xi_1,\xi_3))\omega_3.
\end{multline}
\noindent On the other hand,  \eqref{eighteen} implies
\begin{multline}\label{twentysix}
2T^0_{\xi_1}I_2\ +\ \mathcal{L}^{'}_{\xi_1}I_2\ +\ 2T^0_{\xi_2}I_1\
+\ \mathcal{L}^{'}_{\xi_2}I_1\\
=\ d\eta_1(\xi_2,\xi_1)I_1\ +\ d\eta_2(\xi_1,\xi_2)I_2\ +\ (d\eta_1(\xi_2,\xi_3)\ +\
d\eta_2(\xi_1,\xi_3))I_3.
\end{multline}
\noindent Let us decompose \eqref{twentysix} into $Sp(n)$-invariant components:
\begin{gather}\label{twentyeight}\begin{aligned} {(\mathcal{L}^{'}_{\xi_1}I_2})^{+--} =
-2{T^0_{\xi_1}}^{--+}I_2+d\eta_1(\xi_2,\xi_1)I_1, \quad
{(\mathcal{L}^{'}_{\xi_2}I_1})^{-+-} =
-2{T^0_{\xi_2}}^{--+}I_1+d\eta_2(\xi_1,\xi_2)I_2,
\end{aligned}
\\\label{thirty}
{(\mathcal{L}^{'}_{\xi_1}I_2+\mathcal{L}^{'}_{\xi_2}I_1)}^{--+}\ =\
(d\eta_1(\xi_2,\xi_3)+d\eta_2(\xi_1,\xi_3))I_3.
\end{gather}
\noindent Using (\ref{thirty}) and \eqref{thirq}, we obtain
\begin{gather*}
2\tilde
u=-I_3(({\mathcal{L}^{'}_{\xi_2}I_1})^{--+})+{\frac{1}{2}}(-d\eta_1(\xi_2,\xi_3)+d\eta_2(\xi_3,\xi_1)
-d\eta_3(\xi_1,\xi_2))Id_H.
\end{gather*}
The latter, together with \eqref{thirq}, tells us that $\tilde u$
commutes with all $I\in Q$.  Now, Lemma \ref{twenty} with
$D=\mathcal{L}^{'}$ implies \eqref{twentytree} and
\eqref{twentyfour}.

The vanishing of the symmetric part of the left hand side in (\ref{eighteen}) for
$k=1,\quad l=2$, combined with \eqref{tornew1} and \eqref{twentytree} yields
$0=-2g(I_3\tilde u X,Y)-2g(I_3\tilde u Y,X).$
As $\tilde u$ commutes with all $I\in \mathbb Q$ we conclude that $\tilde u$ is symmetric.

The rest of the identities can be obtained through a cyclic permutation of (1,2,3).
\end{proof}

We describe the properties of the quaternionic contact torsion more
precisely in the next Proposition.

\begin{prop}\label{torbase}
The torsion  of the Biquard connection satisfies the identities:
\begin{gather}\label{to} T_{\xi_i}=T^0_{\xi_i} + I_iu, \quad i=1,2,3, \\\label{tornew1} T^0_{\xi_i}=\frac12\mathbb
\mathcal{L}_{\xi_i}g, \quad i=1,2.3,
\\\label{tornew2} u=\tilde u - \frac{tr(\tilde u)}{4n}Id_H,
\end{gather}
where the symmetric endomorphism $\tilde u$ on $H$ commuting with $I_1,I_2,I_3$ satisfies
\begin{equation}\label{toru}
\begin{aligned}
\tilde{u}\ =\ \frac12I_1A_1^{+--}\ +\ \frac14\big (-d\eta_1(\xi_2,\xi_3)\ + \
d\eta_2(\xi_3,\xi_1)\ +\ d\eta_3(\xi_1,\xi_2)\big )\, Id_{H}\\
 =\ \frac12I_2A_2^{-+-}\ +
\ \frac14 \big (d\eta_1(\xi_2,\xi_3)\ -\ d\eta_2(\xi_3,\xi_1)\ +\
d\eta_3(\xi_1,\xi_2)\big ) \, Id_{H}\\
=\ \frac12I_3A_3^{--+}\ +\ \frac14 \big
(d\eta_1(\xi_2,\xi_3)\ +\ d\eta_2(\xi_3,\xi_1)\ -\ d\eta_3(\xi_1,\xi_2)\big )\,Id_{H}.
\end{aligned}
\end{equation}
For $n=1$ the tensor $u=0$ and $\tilde u=\frac{tr(\tilde
u)}{4}\text{Id}_H$.
\end{prop}
\begin{proof}
Expressing the Lie derivative in terms of the Biquard connection,
using that $\nabla$ preserves the splitting $H\oplus V$, we see that
for $X,Y\in H$ we have
$$\mathcal L_{\xi_i}g(X,Y)=g(\nabla_X\xi_i,Y)+g(\nabla_Y\xi_i,X) +g(T_{\xi_i}X,Y)+g(T_{\xi_i}Y,X)=
2g(T^0_{\xi_i}X,Y).
$$
To show  that $\tilde u$ satisfies  \eqref{toru}, insert \eqref{fifty} into
\eqref{fiftyone} to get
\begin{equation}\label{tet1}\sideremark{tet1}
\theta_3=({\mathcal{L}_{\xi_1}\omega_2})_{|H}-d\eta_2(\xi_1,\xi_2)
\omega_3+d\eta_2(\xi_3,\xi_1)\omega_3.
\end{equation}
A substitution of \eqref{eighteen} and \eqref{twentytree} in
\eqref{tet1} gives
\begin{multline}\label{New5}
A_3= 2{T^0_{\xi_1}}^{-+-}I_2-2I_3\tilde{u}
\\
 +\ d\eta_1(\xi_2,\xi_1)I_1-d\eta_2(\xi_1,\xi_2)I_2+{1\over
2}(d\eta_1(\xi_2,\xi_3)+d\eta_2(\xi_3,\xi_1)-d\eta_3(\xi_1,\xi_2))I_3.
\end{multline}
Now, by comparing the (+++) component on both sides of \eqref{New5}
we see the last equality of \eqref{toru}. The rest of the identities
can be obtained with a cyclic permutation of (1,2,3).

Turning to the rest of the identities, let $\Sigma^2$ and
$\Lambda^2$ denote, respectively, the subspaces of symmetric and
skew-symmetric endomorphisms of $H$. Let $skew:
End(H)\rightarrow\Lambda^2$ be the natural projection with kernel
$\Sigma^2$. We have
\begin{multline}\nonumber
4[T_{\xi_i}]_{{(\Sigma^2\oplus sp(n))}^{\perp}}\ =\ 3skew(T_{\xi_i})\ +\
I_1skew(T_{\xi_i})I_1\ +\ I_2skew(T_{\xi_i})I_2\ +\ I_3skew(T_{\xi_i})I_3
\\
=\ \sum_{s=1}^3{(skew(T_{\xi_i})\ +\ I_sskew(T_{\xi_i})I_s)}.\hskip2truein
\end{multline}
According to Theorem~\ref{biqcon}, $T_{\xi}X\in H$ for $X\in H, \xi \in V$. Hence,
\begin{equation}\label{fo1}
T(\xi,X)=\nabla_{\xi}X-[\xi,X]_H=\nabla_{\xi}X-\mathcal{L}^{'}_\xi(X).
\end{equation}
An application of \eqref{fo1} gives
\begin{multline}\label{thirtyseven}
g(4[T_{\xi_i}]_{{(\Sigma^2\oplus
sp(n))}^{\perp}}X,Y)=-\sum_{s=1}^3{g\big ((\nabla_{\xi_i}I_s)X,I_sY\big )}\\
+{\frac{1}{2}}\sum_{s=1}^3{\{g\big ((\mathcal L_{\xi_i}I_s)X,I_sY\big )-g((\mathcal
L_{\xi_i}I_s)Y,I_sX)\}}.
\end{multline}
 Let $B(H)$ be the orthogonal
complement of $\Sigma^2\oplus sp(n)\oplus sp(1)$ in $\text{End}(H)$.
Obviously, $B(H)\subset \Lambda^2$ and we have the following
splitting of $\text{End}(H)$ into  mutually orthogonal components
\begin{equation}\label{SpltEmdH}
\text{End}(H)=\Sigma^2\oplus sp(n)\oplus sp(1)\oplus B(H).
\end{equation}
If $\Psi$ is an arbitrary section of the bundle $\Lambda^2$ of $M$, the orthogonal
projection of $\Psi$ into $B(H)$ is given by
$[\Psi]_{B(H)}=\Psi^{+--}+\Psi^{-+-}+\Psi^{--+}-[\Psi]_{sp(1)},$
where $[\Psi]_{sp(1)}$ is the orthogonal projection of $\Psi$ onto $sp(1)$. We have also
$[\Psi]_{sp(1)}=\frac{1}{4n}\sum_{s=1}^{3}\sum_{a=1}^{4n}g(\Psi e_a,I_se_a)I_s.$

\noindent Theorem~\ref{biqcon} - (iv) and  the splitting \eqref{SpltEmdH} yield
\begin{equation}\label{eqTsplit}T_{\xi_i}=[T_{\xi_i}]_{(sp(n)\oplus
sp(1))^{\perp}}=[T_{\xi_i}]_{\Sigma^2} + [T_{\xi_i}]_{B(H)} =
T^0_{\xi_i}+[T_{\xi_i}]_{{(\Sigma^2\oplus sp(n))}^{\perp}}-[T_{\xi_i}]_{sp(1)}.
\end{equation}
\noindent Using \eqref{thirtyseven}, Lemma~\ref{twentytwo} and the fact that
$I_s(\nabla_{\xi_i}I_s)\in sp(1)$, we compute
\begin{equation}
\begin{aligned}\label{TspPerp}
4[T_{\xi_i}]_{{(\Sigma^2\oplus sp(n))}^{\perp}} \ -\
[T_{\xi_i}]_{sp(1)}\ =\ -\sum_{s=1}^3\bigl\{\, skew(I_s(\mathcal
L^{'}_{\xi_i}I_s))\ -\ [I_s(\mathcal
L^{'}_{\xi_i}I_s)]_{sp(1)}\bigr\}\\ =\
\sum_{s=1}^3skew(2I_sT^0_{\xi_i}I_s)+4u\ =\ 4u.
\end{aligned}
\end{equation}
\noindent Inserting \eqref{TspPerp} in  \eqref{eqTsplit} completes
the proof.
\end{proof}

\noindent The $Sp(n)$-invariant splitting of \eqref{New5} leads to
the following Corollary.

\begin{cor}\label{fiftytwo}
 The (1,1)-tensors $A_i$ satisfy the equalities
\begin{gather*}
A_3^{+++}\  =\ 2{T^0_{\xi_1}}^{-+-}I_2,\quad
A_3^{+--}\  =\ d\eta_1(\xi_2,\xi_1)I_1,\quad
A_3^{-+-}\  =\ -d\eta_2(\xi_1,\xi_2)I_2,\\
A_3^{--+}\  =\ -2I_3\tilde{u}+{1\over 2
}(d\eta_1(\xi_2,\xi_3)+d\eta_2(\xi_3,\xi_1)-d\eta_3(\xi_1,\xi_2))I_3.
\end{gather*}
Analogous formulas for $A_1$ and $A_2$ can be obtained by a cyclic permutation of
$(1,2,3).$
\end{cor}

\begin{prop}\label{forty}
 The covariant derivative of the quaternionic contact structure with
respect to the Biquard connection is given by
\begin{equation}\label{der}
 \nabla I_i=-\alpha_j\otimes I_k+\alpha_k\otimes I_j,
\end{equation}
where the $sp(1)$-connection 1-forms $\alpha_s$ are determined by
\begin{gather}\label{coneforms}
\alpha_i(X)=d\eta_k(\xi_j,X)=-d\eta_j(\xi_k,X), \quad X\in H, \quad \xi_i\in
V;\\\label{coneform1}
\alpha_i(\xi_s)=d\eta_s(\xi_j,\xi_k)-\delta_{is}\left({tr(\tilde{u})\over
2n}+\frac12\left(d\eta_1(\xi_2,\xi_3)+d\eta_2(\xi_3,\xi_1)+d\eta_3(\xi_1,\xi_2)\right)\right), \quad s=1,2,3.
\end{gather}
\end{prop}
\begin{proof}
The equality \eqref{coneforms} is proved by Biquard in \cite{Biq1}. Using \eqref{fo1}, we
obtain
$$\nabla_{\xi_s}I_i=[T_{\xi_s},I_i]+\mathcal{L}^{'}_{\xi_s}I_i=
[T^0_{\xi_s},I_i]+u[I_s,I_i]+\mathcal{L}^{'}_{\xi_s}I_i.$$ An application of Lemma
\ref{twentytwo} completes the proof.
\end{proof}

\begin{cor}\label{fortyone}\sideremark{fortyone}
The covariant derivative of the distribution $V$ is given by
$$\nabla\xi_i=-\alpha_j\otimes\xi_k+\alpha_k\otimes\xi_j.$$
\end{cor}
We finish this section by expressing the Biquard connection in terms of the Levi-Civita
connection $D^g$ of the metric $g$, namely, we have
\begin{equation}\label{Biq-LC}
\nabla_B Y= D^g_B Y + \sum_{s=1}^3\{((D^g_B\eta_s)Y)\xi_s + \eta_s(B)(I_s u - I_s)Y\},\quad B\in TM,\quad Y\in H.
\end{equation}
\noindent Indeed, for $B=X\in H$ formula \eqref{Biq-LC} follows from the equation
$\nabla_X Y=[D^g_X Y]_H.$ If $B\in V$ we may assume $B=\xi_1$ and for $Z\in H$ we compute
\begin{align*}
2g(D^g_{\xi_1}Y,Z)&=\xi_1g(Y,Z)+g([\xi_1,Y],Z)-g([\xi_1,Z],Y)-g([Y,Z],\xi_1)=\\
&=(\mathcal{L}_{\xi_1}g)(Y,Z)+2g([\xi_1,Y],Z)+d\eta_1(Y,Z)=2g(T_{\xi_1}Y+[\xi_1,Y],Z)\\
&-2g(I_1uY,Z)+2g(I_1Y,Z)=2g(\nabla_{\xi_1}Y,Z)-2g((I_1u-I_1)Y,Z).
\end{align*}
In the above calculation we used \eqref{fo1} and Proposition~\ref{torbase}.

Note that the covariant derivatives $\nabla_B\xi_s$ are also determined by \eqref{Biq-LC} in view of
the relation $g(\nabla_B\xi_s,\xi_k)=\frac{1}{4n}g(\nabla_BI_s,I_k),\ s,k=1,2,3$

\subsection{The Curvature Tensor}

Let $R=[\nabla,\nabla]-\nabla_{[\ ,\ ]}$ be the curvature tensor of $\nabla$. For any
$B,C \in \Gamma(TM)$ the curvature operator $R_{BC}$ preserves the QC structure
on $M$ since $\nabla$ preserves it. In particular $R_{BC}$ preserves the
distributions $H$ and $V$, the  quaternionic structure $\mathbb Q$ on $H$ and the $(2,1)$ tensor
$\varphi.$ Moreover, the action of $R_{BC}$ on $V$ is completely determined by its action
on $H$,
$R_{BC}\xi_i=\varphi^{-1}([R_{BC},I_i]),\ \ i=1,2,3.$
Thus, we may regard $R_{BC}$ as an endomorphism of $H$ and we have $R_{BC} \in
sp(n)\oplus sp(1).$

\begin{dfn}\label{d:Ricci forms}
The Ricci 2-forms $\rho_i$ are defined by
\begin{equation*}
\rho_i(B,C)=\frac{1}{4n}\sum_{a=1}^{4n}g(R(B,C)e_a,I_ie_a),\ B,C\in
\Gamma(TM).
\end{equation*}
\end{dfn}
\noindent Hereafter $\{e_1,\dots,e_{4n}\}$ will denote an orthonormal  basis of $H$. We
decompose the curvature into $sp(n)\oplus sp(1)$-parts. Let  $R^0_{BC} \in sp(n)$ denote
the $sp(n)$-component.

\begin{lemma} The curvature of the Biquard connection
decomposes as follows
\begin{gather}
R_{BC}=R^0_{BC}+\rho_1(B,C)I_1+\rho_2(B,C)I_2+\rho_3(B,C)I_3.
\nonumber\\\label{fiftynine}
[R_{BC},I_i]=2(-\rho_j(B,C)I_k+\rho_k(B,C)I_j),\ B,C\in \Gamma(TM),\\
\label{sixtyone} \rho_i={1\over 2}(d\alpha_i+\alpha_j\wedge \alpha_k),
\end{gather}
where the connection 1-forms $\alpha_s$ are determined in \eqref{coneforms},
\eqref{coneform1}.
\end{lemma}
\begin{proof}
The first two identities follow directly from the definitions. Using \eqref{der}, we
calculate \begin{gather*}[R_{BC},I_1]= \nabla_B(\alpha_3(C)I_2-\alpha_2(C)I_3)
-\nabla_C(\alpha_3(B)I_2-\alpha_2(B)I_3)-(\alpha_3([B,C])I_2-\alpha_2([B,C])I_3)\\
=-(d\alpha_2+\alpha_3\wedge \alpha_1)(B,C)I_3+(d\alpha_3+\alpha_1\wedge
\alpha_2)(B,C)I_2.
\end{gather*} Now \eqref{fiftynine} completes the proof.
\end{proof}

\begin{dfn}\label{d:horizontal Ricci}
The quaternionic contact Ricci tensor (\emph{qc-Ricci tensor} for short) and the
qc-scalar curvature $Scal$  of the Biquard connection are defined by
\begin{equation}\label{e:horizontal ricci}
Ric(B,C)=\sum_{a=1}^{4n}{g(R(e_a,B)C,e_a)}, \qquad Scal=\sum_{a=1}^{4n}Ric(e_a,e_a).
\end{equation}
\end{dfn}
\noindent It is known, cf. \cite{Biq1}, that the qc-Ricci tensor restricted to $H$ is
symmetric. In addition, we define six Ricci-type tensors $\zeta_i,\tau_i,\quad i=1,2,3$
as follows
\begin{equation}\label{e:zeta and tau}
\zeta_i(B,C)={1\over 4n}\sum_{a=1}^{4n}g(R(e_a,B)C,I_ie_a),\quad \tau_i(B,C)={1\over
4n}\sum_{a=1}^{4n}g(R(e_a,I_ie_a)B,C).
\end{equation}
\noindent We shall show that all Ricci-type contractions evaluated on the horizontal
space $H$ are determined by the components of the torsion.  First, define the following
2-tensors on $H$
\begin{gather}\label{tor}
T^0(X,Y)\overset{def}=g((T_{\xi_1}^{0}I_1+T_{\xi_2}^{0}I_2+T_{\xi_3}^{0}I_3)X,Y),\quad
U(X,Y)\overset{def}=g(uX,Y), \quad X,Y\in H.
\end{gather}

\begin{lemma}
The tensors $T^0$ and $U$ are $Sp(n)Sp(1)$-invariant trace-free symmetric tensors with
the properties: \sideremark{propt\\propu}
\begin{gather}\label{propt}
T^0(X,Y)+T^0(I_1X,I_1Y)+T^0(I_2X,I_2Y)+T^0(I_3X,I_3Y)=0,\\\label{propu}
3U(X,Y)-U(I_1X,I_1Y)-U(I_2X,I_2Y)-U(I_3X,I_3Y)=0.
\end{gather}
\end{lemma}
\begin{proof}
The lemma follows directly from \eqref{tors1}, \eqref{toras} of Proposition~\ref{torb}.
\end{proof}

\noindent
We turn to a Lemma, which shall be used later.
\begin{lemma}\label{tehnic}
\sideremark{tehnic\\rizet\\sixtythree}  
For any $X,Y\in H, \quad B\in H\oplus V$, we have
\begin{align}\label{rizet}
& Ric(B,I_iY)+4n\zeta_i(B,Y)\ =\ 2\rho_j(B,I_kY)-2\rho_k(B,I_jY),\\\label{sixtythree}
 & \zeta_i(X,Y)\ =\ -{1\over
2}\rho_i(X,Y)+{1\over
2n}g(I_iuX,Y)+{2n-1\over 2n}g(T_{\xi_i}^{0}X,Y)\\
& \nonumber  \hskip1.5truein  +\ {1\over 2n}g(I_jT_{\xi_k}^{0}X,Y)-{1\over
2n}g(I_kT_{\xi_j}^{0}X,Y).
\end{align}
The Ricci 2-forms evaluated on $H$ satisfy \sideremark{rhoh}
\begin{align}
& \rho_1(X,Y)\  =\ 2g(T_{\xi_2}^{0--+}I_3X,Y)-2g(I_1uX,Y)-{tr(\tilde
u)\over n}\omega_1(X,Y),\nonumber\\
& \rho_2(X,Y)\  =\ 2g(T_{\xi_3}^{0+--}I_1X,Y)-2g(I_2uX,Y)-{tr(\tilde
u)\over n}\omega_2(X,Y),\label{rhoh}\\
& \rho_3(X,Y)\ =\ 2g(T_{\xi_1}^{0-+-}I_2X,Y)-2g(I_3uX,Y)-{tr(\tilde u)\over
n}\omega_3(X,Y).\nonumber
\end{align}
The 2-forms $\tau_s$ evaluated on $H$ satisfy \sideremark{sixtytwo}
\begin{align}
\tau_1(X,Y)\ & =\ \rho_1(X,Y)+2g(I_1uX,Y)+{4\over
n}g(T_{\xi_2}^{0--+}I_3X,Y),\nonumber\\
\tau_2(X,Y)\ & =\ \rho_2(X,Y)+2g(I_2X,Y)+\frac4ng(T_{\xi_3}^{0+--}I_1X,Y),\label{sixtytwo}\\
\tau_3(X,Y)\ &  =\ \rho_3(X,Y)+2g(I_3X,Y)+\frac4ng(T_{\xi_1}^{0-+-}I_2X,Y).\nonumber
\end{align}
For $n=1$ the above formulas hold with  $U=0$.
\end{lemma}
\begin{proof}
>From (\ref{fiftynine}) we have
\begin{multline*}
Ric(B,I_1Y)\ +\ 4n\zeta_1(B,Y)\ =\ \sum_{a=1}^{4n}\{R(e_a,B,I_1Y,e_a)\ +\ R(e_a,B,Y,I_1e_a)\}\\
=\ \sum_{a=1}^{4n}\{-2\rho_2(e_a,B)\omega_3(Y,e_a)\ +\ 2\rho_3(e_a,B)\omega_2(Y,e_a))\} \
=\ 2\rho_2(B,I_3Y)\ -\ 2\rho_3(B,I_2Y),
\end{multline*}
Use (\ref{fiftyone}) and (\ref{sixtyone}) to get
$\rho_1(X,Y)=A_1(X,Y)-{1\over 2}\alpha_1([X,Y]_V) = A_1(X,Y)+\sum_{s=1}^3\omega_s(X,Y)\alpha_1(\xi_s).$
Now, Corollary~\ref{fiftytwo} and Corollary~\ref{forty} imply the first equality in
\eqref{rhoh}. The other two equalities in \eqref{rhoh} can be obtained in the same
manner.

Letting $b(X,Y,Z,W)\ =\ 2\sigma_{X,Y,Z}\{\sum_{l=1}^3{\omega_l(X,Y)g(T_{\xi_l}Z,W)}\},$
where $\sigma_{X,Y,Z}$ is the cyclic sum over $X,Y,Z$, we have
\begin{align}
& \label{sixtyfive} \sum_{a=1}^{4n}b(X,Y,e_a,I_1e_a)\ =\
4g(I_1uX,Y)+8g(I_2T_{\xi_3}^{0-+-}X,Y),\\
& \label{sixtysix} \sum_{a=1}^{4n}b(e_a,I_1e_a,X,Y)\ =\
(8n-4)g(T_{\xi_1}^{0}X,Y)+(8n+4)g(I_1uX,Y)\\
& \nonumber \hskip1truein + \ 4g(T_{\xi_2}^{0}I_3X,Y)-4g(T_{\xi_3}^{0}I_2X,Y).
\end{align}
\noindent The first Bianchi identity gives \sideremark{b1\\b2}
\begin{multline}\label{b1}
 4n(\tau_1(X,Y)+2\zeta_1(X,Y))=\\ \
 \sum_{a=1}^{4n}\{R(e_a,I_1e_a,X,Y)+R(X,e_a,I_1e_a,Y)+R(I_1e_a,X,e_a,Y)\}
 =\ \sum_{a=1}^{4n}b(e_a,I_1e_a,X,Y),
\end{multline}
\begin{multline}
\label{b2}
4n(\tau_1(X,Y)-\rho_1(X,Y))=\sum_{a=1}^{4n}\{R(e_a,I_1e_a,X,Y)-R(X,Y,e_a,I_1e_a)\}\\
=\ {1\over
2}\sum_{a=1}^{4n}\{b(e_a,I_1e_a,X,Y)-b(e_a,I_1e_a,Y,X)-b(e_a,X,Y,I_1e_a)+b(I_1e_a,X,Y,e_a)\}.
\end{multline}
\noindent Consequently, (\ref{sixtyfive}), (\ref{sixtysix}), \eqref{b1} and \eqref{b2}
yield the first set of equalities in (\ref{sixtytwo}) and (\ref{sixtythree}). The  other
equalities in (\ref{sixtytwo}) and (\ref{sixtythree}) can be shown similarly.
\end{proof}

\begin{thrm}\label{sixtyseven}
{Let $(M^{4n+3},g,\mathbb Q)$ be a quaternionic contact} $(4n+3)$-dimensional manifold,
$n>1.$  For any $X,Y \in H$ the qc-Ricci tensor and the qc-scalar curvature satisfy
\begin{equation}\label{sixtyfour}
\begin{aligned}
& Ric(X,Y) \ =\ (2n+2)T^0(X,Y) +(4n+10)U(X,Y)+(2n+4){tr(\tilde{u})\over
n}g(X,Y)\\
& Scal \ =\ (8n+16)tr(\tilde{u}).
\end{aligned}
\end{equation}
For $ \quad n=1,\qquad Ric(X,Y)\ =\ 4T^0(X,Y)+6{tr(\tilde{u})\over n}g(X,Y).$
\end{thrm}
\begin{proof}
The proof follows from Lemma \ref{tehnic}, (\ref{sixtythree}), \eqref{rhoh} and
\eqref{rizet}.  If $n=1$, recall that $U=0$ to obtain the last equality.
\end{proof}

\begin{cor}\label{scalcur}
The qc-scalar curvature satisfies the equalities
$$\frac{Scal}{2(n+2)}=\sum_{a=1}^{4n}\rho_i(I_ie_a,e_a)=\sum_{a=1}^{4n}\tau_i(I_ie_a,e_a)
=-2\sum_{a=1}^{4n}\zeta_i(I_ie_a,e_a),\quad i=1,2,3.$$
\end{cor}

We determine the unknown function $\lambda$ in \eqref{torv} in the
next Corollary.

\begin{cor}\label{fortytwo}
The torsion of the Biquard connection restricted to $V$ satisfies
the equality
\begin{equation}\label{cornv}
T(\xi_i,\xi_j)=-\frac{Scal}{8n(n+2)}\xi_k-[\xi_i,\xi_j]_H.
\end{equation}
\end{cor}

\begin{proof} A small calculation using Corollary \ref{fortyone} and
  Proposition~\ref{forty}, gives
$$T(\xi_i,\xi_j)=\nabla_{\xi_i}\xi_j-\nabla_{\xi_j}\xi_i-[\xi_i,\xi_j]=
-{tr(\tilde{u})\over n}\xi_k-[\xi_i,\xi_j]_H.$$ Now, the assertion
follows from the second equality in \eqref{sixtyfour}.
\end{proof}

\begin{cor}
The tensors $T^0,U, \tilde{u}$  do not depend on the choice of the local basis.
\end{cor}
\section{QC-Einstein quaternionic contact structures}

The aim of this section is to show that the vanishing of the torsion
of the quaternionic contact structure implies that the qc-scalar
curvature is constant  and to prove our classification
Theorem~\ref{Ein2MO}. The Bianchi identities will have an important
role in the analysis.
\begin{dfn}\label{d:qc Einstein}
A quaternionic contact structure is \emph{qc-Einstein} if the qc-Ricci tensor is
trace-free,
$$Ric(X,Y)=\frac{Scal}{4n}g(X,Y), \quad X,Y\in H.$$
\end{dfn}

\begin{prop}\label{tor-ein}
{A quaternionic contact manifold $(M,g,\mathbb Q)$} is a qc-Einstein if and only if the
quaternionic contact torsion vanishes identically, $T_{\xi}=0, \xi\in V$.
\end{prop}
\begin{proof}
If $(\eta,\mathbb Q)$ is qc-Einstein structure then $T^0=U=0$
because of \eqref{sixtyfour}.  We will use the same symbol $T^0$ for
the corresponding endomorphism of the 2-tensor $T^0$ on $H$.
According to \eqref{tor}, we have
$T^0=T^0_{\xi_1}I_1+T^0_{\xi_2}I_2+T^0_{\xi_3}I_3.$
Using first \eqref{tors1} and then \eqref{tors2}, we compute
\begin{equation}\label{New8}
{(T^0)}^{+--}\ =\ (T^0_{\xi_2})^{--+}I_2\ +\ (T^0_{\xi_3})^{-+-}I_3\
=\ 2(T^0_{\xi_2})^{--+}I_2.
\end{equation}
Hence, $T_{\xi_2}=T^0_{\xi_2}+I_2u$ vanishes. Similarly $T_{\xi_1}=T_{\xi_3}=0$.   The
converse  follows from \eqref{sixtyfour}.
\end{proof}

\begin{prop}\label{rhoxi}
For $X\in V$ and any cyclic permutation $(i,j,k)$ of $(1,2,3)$ we have
\begin{align}\label{sixtyeight}
& \rho_i(X,\xi_i)=-\frac{X(Scal)}{32n(n+2)}+{1\over
2}(\omega_i([\xi_j,\xi_k],X)-\omega_j([\xi_k,\xi_i],X)-\omega_k([\xi_i,\xi_j],X)),\\
& \label{seventy} \rho_i(X,\xi_j)=\omega_j([\xi_j,\xi_k],X),
\quad \rho_i(X,\xi_k)=\omega_k([\xi_j,\xi_k],X),\\
&\label{rhotor}
\rho_i(I_kX,\xi_j)=-\rho_i(I_jX,\xi_k)=g(T(\xi_j,\xi_k),I_iX)=\omega_i([\xi_j,\xi_k],X).
\end{align}
\end{prop}
\begin{proof}
Since $\nabla$ preserves the splitting $H\oplus V$, the first Bianchi identity,
\eqref{cornv} and \eqref{fiftynine} imply
\begin{multline}\label{as1}
\sideremark{as1}
2\rho_i(X,\xi_i)+2\rho_j(X,\xi_j)=g(R(X,\xi_i)\xi_j,\xi_k)+g(R(\xi_j,X)\xi_i,\xi_k)\\
=\ \sigma_{\xi_i,\xi_j,X}\{
g((\nabla_{\xi_i}T)(\xi_j,X),\xi_k)+g(T(T(\xi_i,\xi_j),X),\xi_k)\}\\
=\
g((\nabla_{X}T)(\xi_i,\xi_j),\xi_k)+g(T(T(\xi_i,\xi_j),X),\xi_k)
=\ -\frac{X(Scal)}{8n(n+2)}-2\omega_k([\xi_i,\xi_j],X).
\end{multline}
\noindent Summing up the first two equalities in \eqref{as1} and subtracting the third
one, we obtain \eqref{sixtyeight}. Similarly,
\begin{gather*}
2\rho_k(\xi_j,X)=g(R(\xi_j,X)\xi_i,\xi_j)=
\sigma_{\xi_i,\xi_j,X}\{g((\nabla_{\xi_i}T)(\xi_j,X),\xi_j)+g(T(T(\xi_i,\xi_j),X),\xi_j)\}\\
=g(T(T(\xi_i,\xi_j),X),\xi_j)=g(T(-[\xi_i,\xi_j]_H,X),\xi_j)=g([[\xi_i,\xi_j]_H,X],\xi_j)\\
=-d\eta_j([\xi_i,\xi_j]_H,X)=-2\omega_j([\xi_i,\xi_j],X).
\end{gather*}
Hence, the second equality in (\ref{seventy}) follows. Analogous calculations show the
validity of the first equality in \eqref{seventy}. Then, \eqref{rhotor} is a consequence
of \eqref{seventy} and \eqref{cornv}.
\end{proof}

The vertical derivative of the qc-scalar curvature is determined in the next
Proposition.

\begin{prop}
On a QC manifold we have
\begin{equation}\label{ricvert}
 \rho_i(\xi_i,\xi_j)+\rho_k(\xi_k,\xi_j)=\frac{1}{16n(n+2)}\xi_j(Scal).
\end{equation}
\end{prop}
\begin{proof}
Since $\nabla$ preserves the splitting $H\oplus V$, the first Bianchi identity and
\eqref{cornv} imply
\begin{multline*}
-2(\rho_i(\xi_i,\xi_j)+\rho_k(\xi_k,\xi_j))=g(\sigma_{\xi_i,\xi_j,\xi_k}
\{R(\xi_i,\xi_j)\xi_k\},\xi_j)\\
=\ g(\sigma_{\xi_i,\xi_j,\xi_k}\{(\nabla_{\xi_i}T)(\xi_j,\xi_k) +
T(T(\xi_i,\xi_j),\xi_k)\},\xi_j) =-\frac{1}{8n(n+2)}\xi_j(Scal)
\end{multline*}
\end{proof}

\subsection{The Bianchi identities}
In order to derive the essential information contained in the Bianchi identities we need
the next Lemma, which is an application of a standard result in differential geometry.
\begin{lemma}\label{norma}
 In a neighborhood of any point $p\in M^{4n+3}$
and an Q-orthonormal basis \\
$\{X_1(p),X_2(p)=I_1X_1(p)\dots,X_{4n}(p)=I_3X_{4n-3}(p),\xi_1(p),\xi_2(p),\xi_3(p)\}$ of the tangential space at p
there exists a $\mathbb Q$ - orthonormal frame field \\
\centerline{$\{X_1,X_2=I_1X_1,\dots,X_{4n}=I_3X_{4n-3},\xi_1,\xi_2,\xi_3\},
X_{a_|p}=X_a(p),
 \xi_{s_|p}=\xi_s(p), a=1,\dots,4n, i=1,2,3,$} such that the
connection 1-forms of the Biquard connection are all zero at the point p:
\begin{equation}\label{norm}
(\nabla_{X_a}X_b)_|p=(\nabla_{\xi_i}X_b)_|p=(\nabla_{X_a}\xi_t)_|p=(\nabla_{\xi_t}\xi_s)_|p=0,
\end{equation}
for $a,b =1,\dots,4n,  s,t,r=1,2,3.$
In particular,
$$((\nabla_{X_a}I_s)X_b){_|p}=((\nabla_{X_a}I_s)\xi_t){_|p}=((\nabla_{\xi_t}I_s)X_b){_|p}
=((\nabla_{\xi_t}I_s)\xi_r){_|p}=0.$$
\end{lemma}
\begin{proof}
Since $\nabla$ preserves the splitting $H\oplus V$ we can apply the standard arguments
for the existence of a normal frame with respect to a metric connection (see e.g.
\cite{Wu}). We sketch the proof for completeness.

Let $\{\tilde X_1,\dots,\tilde X_{4n},\tilde\xi_1,\tilde\xi_2,\tilde\xi_3\}$ be a
Q-orthonormal basis around p such that $\tilde X_{a_|p}=X_a(p),\quad
\tilde\xi_{i_|p}=\xi_i(p)$. We want to find a modified  frame
$X_a=o^b_a\tilde X_b, \quad \xi_i=o^j_i\tilde\xi_j,$ which
satisfies the normality conditions of the lemma.

Let $\varpi$ be the $sp(n)\oplus sp(1)$-valued connection 1-forms with respect to
$\{\tilde X_1,\dots,\tilde X_{4n},\tilde\xi_1,\tilde\xi_2,\tilde\xi_3\}$,
$$\nabla\tilde X_b=\varpi^c_b\tilde X_c, \quad
\nabla\tilde\xi_s=\varpi^t_s\tilde\xi_t, \quad B\in\{\tilde X_1,\dots,\tilde
X_{4n},\tilde\xi_1,\tilde\xi_2,\tilde\xi_3\}.$$ Let $\{x^1,\dots,x^{4n+3}\}$ be a
coordinate system around p such that $$\frac{\partial}{\partial x^a}(p)=X_a(p),\quad
\frac{\partial}{\partial x^{4n+t}}(p)=\xi_t(p), \quad a=1,\dots,4n, \quad t=1,2,3.$$ One
can easily check that the matrices
$$o^b_a=exp\left(-\sum_{c=1}^{4n+3}\varpi^b_a(\frac{\partial}{\partial x^c})_{|p}x^c\right)\in Sp(n),\quad
o^s_t=exp\left(-\sum_{c=1}^{4n+3}\varpi^s_t(\frac{\partial}{\partial
x^c})_{|p}x^c\right)\in Sp(1)$$ are the desired matrices making the identities
\eqref{norm} true.

Now, the last identity in the lemma is a consequence of the fact that the choice of the
orthonormal basis of $V$ does not depend on the action of $SO(3)$ on $V$ combined with
Corollary~\ref{fortyone} and Proposition~\ref{forty}.
\end{proof}
\begin{dfn}\label{d:normal frame}
We refer to the orthonormal frame constructed in Lemma~\ref{norma} as a \emph{qc-normal
frame}.
\end{dfn}
Let us fix a qc-normal frame
$\{e_1,\dots,e_{4n},\xi_1,\xi_2,\xi_3\}$. We shall denote with
$X,Y,Z$ horizontal vector fields $X,Y,Z\in H$ and keep the notation
for the torsion of type (0,3)  $T(B,C,D)=g(T(B,C),D), B,C,D\in
H\oplus V$.
\begin{prop}\label{bianci}
{On a  quaternionic contact manifold $(M^{4n+3},g ,\mathbb Q)$} the following identities
hold
\begin{align}\label{binc2}
& 2\sum_{a=1}^{4n}(\nabla_{e_a}Ric)(e_a,X)-X(Scal)\ =
\ 4\sum_{r=1}^3Ric(\xi_r,I_rX)-8n\sum_{r=1}^3\rho_r(\xi_r,X);\\
& \label{biric1} Ric(\xi_s,I_sX)\ =\
2\rho_q(I_tX,\xi_s)+2\rho_t(I_sX,\xi_q)+\sum_{a=1}^{4n}(\nabla_{e_a}T)(\xi_s,I_sX,e_a);
\\
& \label{birhozet} 4n(\rho_s(X,\xi_s)-\zeta_s(\xi_s,X))\ =\
2\rho_q(I_tX,\xi_s)+2\rho_t(I_sX,\xi_q)-\sum_{a=1}^{4n}(\nabla_{e_a}T)(\xi_s,X,I_se_a);\\
& \label{zet} \zeta_s(\xi_s,X)\ =\
-\frac{1}{4n}\sum_{a=1}^{4n}(\nabla_{e_a}T)(\xi_s,I_sX,e_a),
\end{align}
\noindent where $s\in\{1,2,3\}$ is fixed and $(s,t,q)$ is an even permutation of
$(1,2,3)$.
\end{prop}

\begin{proof}
The second Bianchi identity  implies
$$2\sum_{a=1}^{4n}(\nabla_{e_a}Ric)(e_a,X)-X(Scal)+2\sum_{a=1}^{4n}Ric(T(e_a,X),e_a)+
\sum_{a,b=1}^{4n}R(T(e_b,e_a),X,e_b,e_a)=0.$$ Apply \eqref{torha} in the last equality to
get \eqref{binc2}.

The first Bianchi identity combined with \eqref{torb0}, \eqref{torha} and the fact that
$\nabla$ preserves the orthogonal splitting $H\oplus V$ yield
\begin{multline*}
Ric(\xi_s,I_sX)=\sum_{a=1}^{4n}\left((\nabla_{e_a}T)(\xi_s,I_sX,e_a)+
2\sum_{r=1}^3\omega_r(I_sX,e_a)T(\xi_r,\xi_s,e_a)\right)=\\
=\sum_{a=1}^{4n}(\nabla_{e_a}T)(\xi_s,I_sX,e_a)
+2T(\xi_s,\xi_t,I_qX) + 2T(\xi_q,\xi_s,I_tX),
\end{multline*}
which, together with \eqref{rhotor}, completes the proof of
\eqref{biric1}.

In a similar fashion, from  the first Bianchi identity, \eqref{torb0}, \eqref{torha} and
the fact that $\nabla$ preserves the orthogonal splitting $H\oplus V$ we can obtain the
proof of \eqref{birhozet}. Finally, take \eqref{rizet} with $B=\xi_i$ and combine the
result with \eqref{biric1} to get \eqref{zet}.
\end{proof}

The following Theorem gives relations between $Sp(n)Sp(1)$-invariant tensors and is
crucial for the solution of the Yamabe problem, which we shall undertake in the last
Section. We define the horizontal divergence $\nabla^*P$ of a (0,2)-tensor field $P$ with
respect to Biquard connection to be the (0,1)-tensor defined by
$\nabla^*P(.)=\sum_{a=1}^{4n}(\nabla_{e_a}P)(e_a,.),$
where $e_a, a=1,\dots,4n$ is an orthonormal basis on $H$.

\begin{thrm}\label{t:horizontal system}
The horizontal divergences of the curvature and torsion tensors satisfy the system $B\,
b\ =\ 0$, where

\[
\mathbf{B}\ =\ \left(%
\begin{array}{ccccc}
 -1 & 6 & 4n-1 & \frac{3}{16n(n+2)} & 0\\
    -1& 0 & n+2 & \frac{3}{16(n+2)} & 0 \\ 1 & -3 & 4 & 0 & -1 \\
\end{array}%
\right),
\]
\[
\mathbf{b}\ =\ \left(%
\begin{array}{ccccc}
  \nabla^*\, T^0, &
  \nabla^*\, U ,&
  \mathbb A ,&
  d\,Scal\,|_{_H} ,&
 \sum_{j=1}^3 Ric\ (\xi_j,I_j \, . \,)
\end{array}%
\right)^t,
\]
\noindent with $T^0$ and $U$ defined in \eqref{tor}  and
\[
\mathbb A(X)=g(I_1[\xi_2,\xi_3]+I_2[\xi_3,\xi_1]+I_3[\xi_1,\xi_2],X).
\]
\end{thrm}

\begin{proof}
Throughout the proof of Theorem \ref{t:horizontal system}  $(s,t,q)$ will denote an even
permutation of $(1,2,3)$. Equations \eqref{sixtyeight} and \eqref{rhotor} yield
\begin{gather}\label{tr1}
\sum_{r=1}^3\rho_r(X,\xi_r)=-\frac{3}{32n(n+2)}X(Scal)-\frac12\mathbb A(X),
\\\label{tr2} \sum_{s=1}^3\rho_q(I_tX,\xi_s)=\mathbb A(X).
\end{gather}
\noindent Using the properties of the torsion described in Proposition~\ref{torbase} and
\eqref{tors1}, we obtain
\begin{align}\label{tr3}
& \sum_{s=1}^3\sum_{a=1}^{4n}(\nabla_{e_a}T)(\xi_s,I_sX,e_a)\ =\
\nabla^*T^0(X)-3\nabla^*U(X),\\\label{tr4} &
\sum_{s=1}^3\sum_{a=1}^{4n}(\nabla_{e_a}T)(\xi_s,X,I_se_a)\ =\
\nabla^*T^0(X)+3\nabla^*U(X).
\end{align}

\noindent  Substituting \eqref{tr2} and \eqref{tr3} in the sum of
\eqref{biric1} written for $s=1,2,3$, we obtain the  third row of
the system. The second row can be obtained by inserting \eqref{zet}
into \eqref{birhozet}, taking the sum over $s=1,2,3$ and applying
\eqref{tr1}, \eqref{tr2}, \eqref{tr3}, \eqref{tr4}.

The second Bianchi identity and applications of \eqref{torha} give
\begin{align}\notag
& \sum_{s=1}^3\left( \sum_{a=1}^{4n}   \Bigl[ \
 (\nabla_{e_a}Ric)(I_sX,I_s{e_a})+4n(\nabla_{e_a}\zeta_s)(I_sX,{e_a})\
\Bigr]  - 2Ric(\xi_s,I_sX) +8n\zeta_s(\xi_s,X)\right)
\\ \label{tr5}  & \hskip.8truein  +8n\sum_{s=1}^3\Bigl[\zeta_s(\xi_t,I_qX)-
\zeta_s(\xi_q,I_tX)-\rho_s(\xi_t,I_qX)+\rho_s(\xi_q,I_tX)\Bigr]=0.
\end{align}
\noindent Using \eqref{rizet}, \eqref{rhoh} as well as \eqref{tors1},\eqref{tors2} and
\eqref{New8} we obtain the next sequence of equalities
\begin{gather}
\sum_{s=1}^3\Bigl[Ric(I_sX,I_s{e_a})+4n\zeta_s(I_sX,{e_a})\Bigr]=2\sum_{s=1}^3\Bigl[\rho_s(I_qX,I_te_a)-
\rho_s(I_tX,I_qe_a)\Bigr] \nonumber\\ =\ -4T^0(X,e_a)+24U(X,e_a)+\frac{3\,Scal}{2n(n+2)}g(X,e_a),\nonumber\\
\label{tr7}
8n\sum_{s=1}^3\Bigl[\zeta_s(\xi_t,I_qX)-\zeta_s(\xi_q,I_tX)\Bigr]\\\nonumber
= \sum_{s=1}^3\Bigl[4Ric(\xi_s,I_sX)-8\rho_s(\xi_s,X)
  + 4\rho_s(\xi_t,I_qX)-4\rho_s(\xi_q,I_tX\Bigr].\\\nonumber
  \sum_{s=1}^3\Bigl[-2Ric(\xi_s,I_sX)
+8n\zeta_s(\xi_s,X)\Bigr]
= \sum_{s=1}^3\Bigl[-4Ric(\xi_s,I_sX) -4\rho_s(\xi_t,I_qX)+4\rho_s(\xi_q,I_tX)\Bigr],\\\nonumber
\end{gather}
\noindent  Substitute  \eqref{tr7} in \eqref{tr5},  and then use
\eqref{tr1} and \eqref{tr2} to get the first row of the system.
\end{proof}
We are ready to prove one of our main observations.

\begin{thrm}\label{maincon}
 The qc-scalar curvature of a qc-Einstein quaternionic contact
manifold of dimension bigger than seven is a global constant. In
addition, the vertical distribution $V$ of a qc-Einstein structure
is integrable and the Ricci tensors are given by
\begin{gather*}\rho_{t|H}=\tau_{t|H}=-2\zeta_{t|H}=-\frac{Scal}{8n(n+2)}
\omega_t \quad s,t=1,2,3., \\ Ric(\xi_s,X)=\rho_s(X,\xi_t)=\zeta_s(X,\xi_t)=0, \quad
s,t=1,2,3.
\end{gather*}
\end{thrm}

\begin{proof}
Suppose the quaternionic contact manifold is qc-Einstein.
According to Proposition~\ref{tor-ein}, the quaternionic contact
torsion vanishes, $T_{\xi}=0, \xi\in V$. Since $n>1$, Theorem
\ref{t:horizontal system} gives immediately that the horizontal
gradient  of the scalar curvature vanishes, i.e., $X(Scal)=0,\quad
X\in H$.  Notice that this fact implies also $\xi(Scal)=0, \xi\in
V$, taking into account that for any $p\in M$ one has
$[e_a,I_se_a]_{|p}=T(e_a,I_se_a)_{|p}=2\xi_{s|p}$.  Now,
\eqref{biric1}, \eqref{seventy}, \eqref{sixtythree}, \eqref{rhoh}
and \eqref{sixtytwo} complete the proof.
\end{proof}

\begin{rmrk}\label{r:open question}
It is an open question whether a seven dimensional qc-Einstein structure must have a constant qc scalar curvature.
\end{rmrk}

\subsection{Examples of qc-Einstein structures}
\begin{exm} The flat model. \end{exm} The quaternionic Heisenberg group $G(\mathbb H)$ with its standard
left invariant quaternionic contact structure (see
Section~\ref{5.2}) is the simplest example. The Biquard connection
coincides with the flat left-invariant  connection on $G(\mathbb
H)$. More precisely, we have the following Proposition.
\begin{prop}
{Any  quaternionic contact manifold $(M, g, \mathbb{Q})$} with flat Biquard connection is
locally isomorphic to $G(\mathbb H)$.
\end{prop}
\begin{proof}
Since the Biquard connection $\nabla$ is flat, there exists a local Q-orthonormal frame\\
$\{T_a,\, I_1T_a,\,  I_2T_a, \, I_3T_a,\,  \xi_1,\,  \xi_2,\,  \xi_3:\  a=1,\dots,n\}$
which is $\nabla$-parallel. Theorem~\ref{maincon} tells us that the quaternionic contact
torsion vanishes and the vertical distribution is integrable. In addition, \eqref{cornv}
and \eqref{torha} yield $[\xi_i,\xi_j]=0$ with the only non-zero commutators
$[I_iT_a,T_a]=2\xi_i, i,j=1,2,3$ (cf. \eqref{e:commutators}). Hence, the manifold has a
local Lie group structure which is locally isomorphic to $\QH$ by the Lie theorems. In
other words, there is a local diffeomorphism $\Phi:\, M \rightarrow \QH$ such that $\eta\
=\ \Phi^* \Theta$, where $\Theta$ is the standard contact form on $\QH$, see
\eqref{e:Heisenbegr ctct forms}.
\end{proof}
\begin{exm}\label{3sas1} {\it The 3-Sasakian Case.}
\end{exm}
Suppose $(M,g)$ is a (4n+3)-dimensional  Riemannian manifold with a given 3-Sasakian
structure, i.e., the cone metric on $M\times\mathbb R$ is a hyperk\"ahler metric, namely,
it has holonomy contained in $Sp(n+1)$ \cite{BGN}. Equivalently, there are
three Killing vector fields $\{\xi_1,\xi_2,\xi_3\}$, which satisfy
\par (i) $g(\xi_i,\xi_j)=\delta_{ij},\ i,j=1,2,3$
\par (ii) $[\xi_i,\xi_j]=-2\xi_k,$ for any cyclic permutation
$(i,j,k)$ of $(1,2,3)$
\par (iii) $(D_B \tilde{I_i})C=g(\xi_i,C)B-g(B,C)\xi_i,\ i=1,2,3,\ B,C\in
\Gamma(TM)$, where $\tilde{I_i}(B)=D_B\xi_i$ and $D$ denotes the Levi-Civita connection.

A 3-Sasakian manifold of dimension $(4n+3)$ is Einstein with positive Riemannian scalar curvature $(4n+2)(4n+3)$
\cite{Kas} and if complete it is compact with finite fundamental group due to Mayer's
theorem (see \cite{BG} for a nice overview of 3-Sasakian spaces).

\noindent Let $H={\{\xi_1,\xi_2,\xi_3\}}^{\perp}$. Then
\begin{gather*}
\tilde{I_i}(\xi_j)=\xi_k,\quad \tilde{I_i}\circ\tilde{I_j}(X)=\tilde{I_k}X, \quad
\tilde{I_i}\circ\tilde{I_i}(X)=-X,\ X\in H,\\
d\eta_i(X,Y)=2g(\tilde{I_i}X,Y),\ X,Y\in H.
\end{gather*}
Defining $V=span\{\xi_1,\xi_2,\xi_3\},\quad I_{i|H}=\tilde{I}_{i|H},\ I_{i|V}=0$ we
obtain a quaternionic contact structure  on $M$ \cite{Biq1}. It is easy to calculate that
\begin{gather*}
{{\xi_i}\lrcorner d\eta_{j}}_{|H}=0,\ d\eta_i(\xi_j,\xi_k)=2,\
d\eta_i(\xi_i,\xi_k)=d\eta_i(\xi_i,\xi_j)=0, \\ A_1=A_2=A_3=0 
\eqref{fiftyone},
\qquad \tilde{u}={1\over 2}Id_H \ \text{ cf. \eqref{thirq}}.
\end{gather*}
This quaternionic contact structure  satisfies the conditions \eqref{bi1} and therefore
it admits the Biquard connection $\nabla.$ More precisely, we have

\par (i) $\nabla_XI_i=0,X\in H$, \quad\quad $\nabla_{\xi_i}I_i=0$,
\quad\quad$\nabla_{\xi_i}I_j=-2I_k,\ \nabla_{\xi_j}I_i=2I_k,$
\par (ii) $T(\xi_i,\xi_j)=-2\xi_k$
\par (iii) $T(\xi_i,X)=0,\ X\in H$.

\noindent
>From  Proposition~\ref{tor-ein}, Theorem~\ref{maincon}, \eqref{coneform1} and
\eqref{sixtyone}, we obtain the following Corollary.
\begin{cor}\label{3sas}
 Any 3-Sasakian manifold  is
a qc-Einstein with positive qc-scalar curvature
$$Scal=16n(n+2).$$
For any  $s,t,r=1,2,3$, the Ricci-type tensors are given by
\begin{equation}
\begin{split}
\rho_{t|H} & \ =\ \tau_{t|H}=-2\zeta_{t|H}=-2\omega_t \\
 Ric(\xi_s,X)& \ =\ \rho_s(X,\xi_t)=\zeta_s(X,\xi_t)=0\\
\rho_s(\xi_t,\xi_r)& \ =\ 0.
\end{split}
\end{equation}
The nonzero parts of the curvature $R$ of the Biquard connection is expressed in terms of
the curvature of the Levi-Civita connection  $R^g$ as follows
\begin{enumerate}[i)]
\item $ R(X,Y,Z,W)  \ =\  R^g(X,Y,Z,W)\ \qquad\newline
\qquad \ +\  \sum_{s=1}^3\{\omega_s(Y,Z)\omega_s(X,W)\ -\
\omega_s(X,Z)\omega_s(Y,W)-2\omega_s(X,Y)\omega_s(Z,W)\};$\newline
\item $ R(\xi,Y,Z,W)\ =\ -R(Y,\xi,Z,W)\ =\ R^g(\xi,Y,Z,W);$\newline
\item $  R(\xi,\bar\xi,Z,W)\ = \ R^g(\xi,\bar\xi,Z,W;)$\newline
\item $R(X,Y,\xi,\bar\xi)\ =\ -4\{ \, \eta_1\wedge\eta_2(\xi,\bar\xi)\omega_3(X,Y)\ +\
\eta_2\wedge\eta_3(\xi,\bar\xi)\omega_1(X,Y)\ +\
\eta_3\wedge\eta_1(\xi,\bar\xi)\omega_2(X,Y) \}$,\newline
\end{enumerate}

\noindent where $X,Y,Z,W\in H$ and $\xi,\bar\xi\in V.$
\end{cor}
In fact, 3-Sasakian spaces are locally the only qc-Einstein manifolds (cf.
Theorem~\ref{Ein2MO}). Before we turn to the proof of this fact we shall consider some
special cases of QC-structures suggested by the above example. These structures will be
relevant in Chapter 6, see for ex. Theorem \ref{crfth1}. We recall that the Nijenhuis
tensor $N_{I_i}$ corresponding to $I_i$ on $H$   is defined as usual by
$N_{I_i}(X,Y)=[I_iX,I_iY]-[X,Y]-I_i[I_iX,Y]-I_i[X,I_iY],\quad X,Y\in H.$

\begin{dfn}\label{d:hyper cplx ctct}
{A quaternionic contact structure $(M, g, \mathbb{Q})$ }is said to be
\emph{hyperhermitian contact} (\emph{HC structure} for short) if the horizontal bundle
$H$ is formally integrable with respect to $I_1,I_2,I_3$ simultaneously, i.e. for
$i=1,2,3$ and  any $X,\, Y\in H$, we have
\begin{equation}\label{e:vertical Nijenhuis}
N_{I_i}(X,Y)=0 \quad {\text mod}\quad  V.
\end{equation}
\end{dfn}

\noindent In fact a QC structure is locally a HC structure exactly when  two of the
almost complex structures on $H$ are formally integrable due to the next identity
essentially established in \cite[(3.4.4)]{AM}
\begin{multline*}
2N_{I_3}(X,Y)-N_{I_1}(X,Y)+I_2N_{I_1}(I_2X,Y)+I_2N_{I_1}(X,I_2Y)-N_{I_1}(I_2X,I_2Y)- \\
N_{I_2}(X,Y)+I_1N_{I_2}(I_1X,Y)+I_1N_{I_2}(X,I_1Y)-N_{I_2}(I_1X,I_1Y) = 0 \quad {\text
mod}\quad V.
\end{multline*}
\noindent On the other hand, the  Nijenhuis tensor has the following expression in terms
of a connection $\nabla$ with torsion $T$ satisfying \eqref{der}(see e.g. \cite{Iv})
\begin{equation}\label{ntor}
N_{I_i}(X,Y)=
T_{Ii}^{0,2}(X,Y)+\beta_i(Y)I_jX-\beta_i(X)I_jY-I_i\beta_i(Y)I_kX+I_i\beta_i(X)I_kY,
\end{equation}
\noindent where  the 1-forms $\beta_i$ and
the (0,2)-part of the torsion $T_{I_i}^{0,2}$ with respect to the almost complex
structure $I_i$ are defined on $H$, correspondingly, by
\begin{gather}\label{ni1f}
\beta_i=\alpha_j+I_i\alpha_k,\\\label{torhypcom}
T_{I_i}^{0,2}(X,Y)=T(X,Y)-T(I_iX,I_iY)+I_iT(I_iX,Y)+I_iT(X,I_iY).
\end{gather}
Applying the above formulas to the Biquard connection and taking into account
\eqref{torha} one sees that \eqref{e:vertical Nijenhuis} is equivalent to
$(\beta_i)_{|_H}=0$. Hence we have the following proposition.

\begin{prop}\label{hypcom1}
{A quaternionic contact structure  $(M, g, \mathbb{Q})$ }is a hyperhermitian contact
structure if and only if the connection 1-forms satisfy the relations
\begin{equation}\label{intqc}
\alpha_j(X)=\alpha_k(I_iX), \quad X\in H
\end{equation}
The Nijenhuis tensors of a HC
structure satisfy $ N_{I_i}(X,Y)= T_{I_i}^{0,2}(X,Y),\quad X,Y\in H. $
\end{prop}

Given a QC structure $(M, g,\mathbb Q)$ let us consider
the three almost complex structures $(\eta_i,\tilde I_i)$
\begin{equation}\label{thcr}
\tilde I_iX = I_iX, \quad X\in H, \quad \tilde I_i(\xi_j)=\xi_k,\quad \tilde
I_i(\xi_i)=0.
\end{equation}
\noindent With these definitions $(\eta_i,\tilde I_i)$ are almost
CR structures (i.e. possibly non-integrable) exactly when the QC
structure is HC since the condition $d\eta_i(\tilde I_iX,\tilde
I_i\xi_j)=d\eta_i(X,\xi_j)$ is equivalent to
$\alpha_k(X)=-\alpha_j(I_iX)$ in view of \eqref{coneforms}. Hence,
$d\eta_i$ is a (1,1)-form with respect to $\tilde I_i$ on
$\xi_i^{\perp}=H\oplus\{\xi_j,\xi_k\}$ and a HC structure supports a
non integrable hyper CR-structure $(\eta_i,\tilde I_i)$.

A natural question is to examine when $\tilde I_i$ is  formally integrable, i.e
$N_{\tilde I_i}=0 \quad {\text mod} \quad \xi_i$.

\begin{prop}
{Let $(M, g, \mathbb{Q})$ be a hyperhermitian contact structure}. Then the CR structures
$(\eta_i,\tilde I_i)$ are integrable
if and only if the next two equalities hold
\begin{gather}\label{crint2}
d\eta_j(\xi_k,\xi_i)=d\eta_k(\xi_i,\xi_j),\hskip1truein
d\eta_j(\xi_j,\xi_i)-d\eta_k(\xi_k,\xi_i)=0.
\end{gather}
\end{prop}

\begin{proof} From \eqref{torha} it follows $T^{0,2}_{\tilde I_1}(X,Y)=0$ using also \eqref{torhypcom}.
Substituting the latter into \eqref{ntor} taken with respect to $\tilde I_i$  shows
$N_{\tilde I_i}|_H=0 \quad {\text mod} \quad \xi_i$ is equivalent to \eqref{intqc}.
Corollary~\ref{fortyone} implies
\begin{multline*}N_{\tilde I_i}(X,\xi_j)=(\alpha_j(I_iX)+\alpha_k(X))\xi_i+
(\alpha_j(\xi_k)+\alpha_k(\xi_j))I_kX + (\alpha_j(\xi_j)- \alpha_k(\xi_k))I_jX+\\
+T(\xi_k,I_iX)-I_iT(\xi_k,X)-T(\xi_j,X)-I_iT(\xi_j,I_iX).
\end{multline*}

\noindent Taking the trace part and the trace-free part in the
right-hand side allows us to conclude that $N_{\tilde
I_i}(X,\xi_j)=0 \quad {\text mod}\quad \xi_i$ is equivalent to the
system
\begin{gather*}
T(\xi_k,I_1X)-I_1T(\xi_k,X)-T(\xi_j,X)-I_1T(\xi_j,I_1X)=0,\\
\alpha_j(\xi_k)+\alpha_k(\xi_j)=0 \quad \alpha_j(\xi_j)-\alpha_k(\xi_k)=0.
\end{gather*}
An application of Proposition~\ref{torbase}, \eqref{tors1} and
\eqref{tors2} shows the first equality is trivially satisfied, while
\eqref{coneform1} tells us that the other equalities are equivalent
to  \eqref{crint2}.
\end{proof}

\subsection{Proof of Theorem~\ref{Ein2MO}}
\begin{proof}[Proof of Theorem~\ref{Ein2MO}]
The equivalence of a) and c) was proved in Proposition
\ref{tor-ein}. We are left with proving the implication a) implies
b). Let $(M,\tilde g, \mathbb{Q})$ be a qc-Einstein manifold with
qc-scalar curvature $\overline{Scal}$. According to
Theorem~\ref{maincon}, $\overline{Scal}$ is a global constant on
$M$. We define $\eta=\frac{\overline{Scal}}{16n(n+2)}\, \tilde\eta$.
Then $(M,g , \mathbb{Q})$ is a qc-Einstein manifold with qc-scalar
curvature $Scal=16n(n+2),$ horizontal distribution $H=Ker(\eta)$ and
involutive vertical distribution $V= span\{\xi_1,\xi_2,\xi_3\}$ (see
\eqref{New19},\eqref{New24} and \eqref{e:conf change scalar curv}).

We shall show that the Riemannian cone is a hyperk\"ahler manifold.
Consider the structures defined by \eqref{thcr}.
We have the relations
\begin{gather} \eta_i(\xi_j)=\delta_{ij},\quad \eta_i\tilde
I_j=-\eta_j\tilde I_i=\eta_k,\qquad \tilde I_i\xi_j=-\tilde
I_j\xi_i=\xi_k\nonumber\\\label{New11} \tilde I_i\tilde
I_j-\eta_j\otimes\xi_i=-\tilde I_j\tilde
I_i+\eta_i\otimes\xi_j=\tilde I_k\qquad \\\nonumber \tilde
I_i^2=-Id+\eta_i\otimes\xi_i,\qquad \eta_i\tilde I_i=0,\qquad \tilde
I_i\xi_i=0,\qquad g(\tilde I_i .,\tilde I_i
.)=g(.,.)-\eta_i(.)\,\eta_i(.).
\end{gather}

Let $D$ be the Levi-Civita connection of the metric $g$ on $M$
determined by the structure $(\eta,Q)$. The next step is to show
\begin{equation}\label{New9}\sideremark{New9}
D\tilde I_i=Id\otimes\eta_i-g\otimes\xi_i-\sigma_j\otimes\tilde I_k+
\sigma_k\otimes\tilde I_j,
\end{equation}
for some appropriate 1-forms $\sigma_s$ on $M$. We consider all
possible cases.

\noindent \textbf{Case 1}[$X,Y,Z\in H$] The well known formula
\begin{multline}\label{New14}
2\,g(D_AB,C)\ =\ Ag(B,C)\ +\ Bg(A,C)\ -\ Cg(A,B)\\
+g([A,B],C)-g([B,C],A)+g([C,A],B),\qquad A,B,C\in \Gamma(TM)
\end{multline}
yields
\begin{equation}\label{H_New9}
2\,g(\,(D_X\tilde I_i)\,Y,
Z)=d\omega_i(X,Y,Z)-d\omega_i(X,I_iY,I_iZ)+g(N_i(Y,Z),I_iX).
\end{equation}
\noindent  We compute $d\omega_i$ in terms of the Biquard
connection. Using \eqref{torha}, \eqref{der} and \eqref{ni1f}, we
calculate
\begin{multline}\label{Omega2Nab}
d\omega_i(X,Y,Z)-d\omega_i(X,I_iY,I_iZ)\ =\
=\ -2\,\alpha_j(X)\,\omega_k(Y,Z)+2\alpha_k(X)\,\omega_j(Y,Z)\\
-\beta_i(Y)\,\omega_k(Z,X)\ -\ I_i\beta_i(Y)\,\omega_j(Z,X)\ -\
\beta_i(Z)\,\omega_k(X,Y)\ -\ I_i\beta_i(Z)\, \omega_j(X,Y).
\end{multline}
\noindent A substitution of \eqref{ntor} and \eqref{Omega2Nab} in
\eqref{H_New9} gives
\begin{equation}\label{qkcon1}
g((D_X\tilde I_i)Y,
Z)=-\alpha_j(X)\omega_k(Y,Z)+\alpha_k(X)\omega_j(Y,Z).
\end{equation}
Letting $\sigma_i(X)=\alpha_i(X)$, we obtain equation \eqref{New9}.

\noindent \textbf{Case 2} [ $\xi_s, \xi_t \in V$ and $Z\in H$] Using
the integrability of the vertical distribution $V$ and
\eqref{New14}, we compute
\begin{multline*}
2g((D_{\xi_s}\tilde I_i)\xi_t, Z)=2g(D_{\xi_s}\tilde I_i\xi_t,
Z)+2g(D_{\xi_s}\xi_t, I_iZ)=\\ -g([\tilde
I_i\xi_t,Z],\xi_s)-g([\xi_s,Z],\tilde I_i\xi_t)-g([\xi_s,I_iZ],\xi_t)-g([\xi_t,I_iZ],\xi_s).
\end{multline*}
An application of \eqref{bi1} allows to conclude $g((D_{\xi_s}\tilde
I_i)\xi_t, Z)=0$ for any $i,s,t\in\{1,2,3\}$.

\noindent \textbf{Case 3} [ $X,Y\in H$ and $C\in V$] First, let
$C=\xi_1$. We have
\begin{multline*}
2g((D_X\tilde I_1)Y, \xi_1)=2g(D_X\tilde I_1Y,\xi_1)\\
=\ -\xi_1g(X,I_1Y)
+g([X,I_1Y],\xi_1)-g([X,\xi_1],I_1Y)-g([I_1Y,\xi_1],X)
\\=-(\mathcal{L}_{\xi_1}g)(X,I_1Y)+\eta_1([X,I_1Y])=-d\eta_1(X,I_1Y) = -2g(X,Y).
\end{multline*}
after using \eqref{tornew1}, $T_{\xi_s}=0,\ s=1,2,3,$ and
\eqref{con1}.

For $C=\xi_2$, we calculate applying  \eqref{New11} and
\eqref{New14} that
\begin{multline*}
2g((D_X\tilde I_1)Y, \xi_2)=2g(D_X\tilde I_1Y,\xi_2)+2g(D_XY, \xi_3)\\
=\ -\xi_1g(X,I_iY)-\xi_3g(X,Y) +g([X,I_iY],\xi_2)-g([X,\xi_2],I_iY)\\
\ -g([I_iY,\xi_2],X)+g([X,Y],\xi_3)-g([X,\xi_3],Y)-g([Y,\xi_3],X)\\=
-(\mathcal{L}_{\xi_2}g)(X,I_1Y)-(\mathcal{L}_{\xi_3}g)(X,Y)+
\eta_2([X,I_1Y])+\eta_3(]X,Y])=0.
\end{multline*}
The other possibilities in this case can be checked in a similar
way.

\noindent \textbf{Case 4} [$X\in H$ and $ A,B\in V].$ We verify
\eqref{New9} for $\tilde I_1, A=\xi_1,B=\xi_2$ and $\tilde I_1,
A=\xi_2,B=\xi_3$ since the other verifications  are similar. Using
the integrability of $V$, \eqref{New14}, \eqref{New11} and
\eqref{coneforms}, we find
\begin{multline*}
2g((D_X\tilde I_1)\xi_1,\xi_2)=2g(D_X\xi_1,\xi_3)=g([X,\xi_1],\xi_3)-g([X,\xi_3],\xi_1)=-2\alpha_2(X),\\
2g((D_X\tilde I_1)\xi_2,
\xi_3)=2g(D_X\xi_3,\xi_3)-2g(D_X\xi_2,\xi_2)=0.
\end{multline*}

\noindent \textbf{Case 5} [$A,B,C\in V$] Let us extend the
definition of the three 1-forms $\sigma_s$ on $V$ as follows
\begin{gather}\label{qkcon2}
\sigma_i(\xi_i)=1+\frac{1}{2}(d\eta_i(\xi_j,\xi_k)-d\eta_j(\xi_k,\xi_i)-d\eta_k(\xi_i,\xi_j))\\\nonumber
\sigma_i(\xi_j)=d\eta_j(\xi_j,\xi_k),\qquad
\sigma_i(\xi_k)=d\eta_k(\xi_j,\xi_k).
\end{gather}
A small calculation leads  to the formula
\begin{gather}\label{I_V}
g(\tilde I_iA, B)=(\eta_j\wedge\eta_k)(A,B).
\end{gather}
On the other hand, we have
\begin{multline}\label{Dxi_V}
2(D_A\eta_i)(B)=2g((D_A\xi_i,B)\\
=A\eta_i(B)+\xi_ig(A,B)-B\eta_i(A)+g([A,\xi_i],B)-\eta_i([A,B])-g([\xi_i,B],A)\\
=\sum_{s,t=1}^3\Bigl[\eta_s(A)\,\eta_t(B)\,d\eta_i(\xi_s,\xi_t)\ -\
\eta_s(A)\,\eta_t(B)\,d\eta_t(\xi_s,\xi_i)\ -\
\eta_s(B)\,\eta_t(A)\,d\eta_t(\xi_s,\xi_i)\Bigr]\\
=\ 2\,\eta_j\wedge\eta_k(A,B)\ -\ 2\, \sigma_j(A)\,\eta_k(B)\ +\
2\,\sigma_k(A)\,\eta_j(B).
\end{multline}
With the help of \eqref{I_V} and \eqref{Dxi_V} we see
\begin{multline}
g((D_A\tilde
I_i)B,C)=D_A(\eta_j\wedge\eta_k)(B,C)=[D_A(\eta_j)\wedge\eta_k+\eta_j\wedge
D_A(\eta_k)](B,C)\\
=((A\lrcorner(\eta_k\wedge\eta_i)-\sigma_k(A)\eta_i +\sigma_i(A)\eta_k)\wedge\eta_k)(B,C) \\
\hskip.5truein +
(\eta_j\wedge(A\lrcorner(\eta_i\wedge\eta_j)-\sigma_i(A)\eta_j+\sigma_j(A)\eta_i))(B,C)\\
=(\eta_k(A)\eta_i\wedge\eta_k(B,C)+\eta_j(A)\eta_i\wedge\eta_j(B,C))-\sigma_j(A)g(\tilde
I_kB,C)+\sigma_k(A)g(\tilde
I_jB,C)\\=\eta_i(B)g(A,C)-g(A,B)\eta_i(C)-\sigma_j(A)g(\tilde
I_kB,C)+\sigma_k(A)g(\tilde I_jB,C).
\end{multline}

\noindent \textbf{Case 6} [$A\in V$ and $ Y,Z\in H].$ Let $A=\xi_s,\
s\in\{1,2,3\}$. The right hand side of \eqref{New9} is equal to
$-\sigma_j(\xi_s)\omega_k(Y,Z)+\sigma_k(\xi_s)\omega_j(Y,Z)$. On the
left hand side of \eqref{New9}, we have
\begin{multline}
2g((D_{\xi_s}\tilde I_i)Y,Z)=2g(D_{\xi_s}(I_iY),Z)+2g(D_{\xi_s}Y,I_iZ)\\
=\{\xi_sg(I_iY,Z)+ g([\xi_s,I_iY],Z)-g([\xi_s,Z],I_iY)-g([I_iY,Z],\xi_s)\}\\
+\{\xi_sg(Y,I_iZ)+ g([\xi_s,Y],I_iZ)-g([\xi_s,I_iZ],Y)-g([Y,I_iZ],\xi_s)\}\\
=g((\mathcal L_{\xi_s}I_i)Y,Z)-g((\mathcal
L_{\xi_s}I_i)Z,Y)+\omega_i(I_sY,Z)+\omega_i(Y,I_sZ)
\end{multline}
Now, recall Lemma~\ref{twentytwo} to compute the skew symmetric part
of $g((\mathcal L_{\xi_s}I_i)Y,Z)$, and also use formulas
\eqref{qkcon2}, to get for $i=s$
\[
g(\,(D_{\xi_i}\tilde I_i)\,Y,Z)\ =\
d\eta_i(\xi_i,\xi_j)\,\omega_j(Y,Z) \ +\
d\eta_i(\xi_i,\xi_k)\,\omega_k(Y,Z)\ =\
-\sigma_j(\xi_i)\,\omega_k(Y,Z)\ +\ \sigma_k(\xi_i)\,\omega_j(Y,Z).
\]
Similarly, for $s=j$ we have
\begin{multline*}
g(\,(D_{\xi_j}\tilde I_i)\, Y,Z)\ =\
d\eta_j(\xi_i,\xi_j)\,\omega_j(Y,Z)\
-\ \frac{1}{2}\,\Bigl(-d\eta_i(\xi_j,\xi_k)\ +\ d\eta_j(\xi_k,\xi_i)\ -\ d\eta_k(\xi_i,\xi_j)\Bigr)\,\omega_k(Y,Z)\\
 -\ \omega_k(Y,Z)\
=\ -\sigma_j(\xi_j)\,\omega_k(Y,Z)\ +\
\sigma_k(\xi_j)\,\omega_j(Y,Z),
\end{multline*}
which completes the proof of \eqref{New9}.

At this point, consider the Riemannian cone $N=M\times \mathbb R^+$
with the cone metric $g_N=t^2g+dt^2$ and the almost complex
structures
$$\phi_i(E,f\frac{d}{dt})=(\tilde I_iE+\frac{f}{t}\xi_i, -t\eta_i(E)\frac{d}{dt}), \quad i=1,2,3,\quad E\in \Gamma(TM).$$
Using the O'Neill formulas for warped product  \cite[p.206]{On},
\eqref{New11} and the just proved \eqref{New9} we conclude (see also
\cite{MO}) that the Riemannian cone $(N,g_N,\phi_i, i=1,2,3)$ is a
quatrnionic K\"ahler manifold with connection 1-forms defined by
\eqref{qkcon1} and \eqref{qkcon2}. It is classical result (see e.g
\cite{Bes}) that a quaternionic K\"ahler manifolds of dimension bigger
than 4 are Einstein with non-negative scalar curvature. This fact implies
that the cone $N=M\times \mathbb R^+$ with the warped product metric $g_N$
must be Ricci flat (see e.g. \cite[p.267]{Bes}) and therefore it is
locally hyperk\"ahler (see e.g. \cite[p.397]{Bes}). This means that
locally there exists a $SO(3)$-matrix $\Psi$ with smooth entries such that
the triple $(\tilde\phi_1,\tilde\phi_2,\tilde\phi_3)= \Psi\cdot
(\phi_1,\phi_2,\phi_3)^t$ is $D$-parallel. Consequently
$(M,\Psi\cdot\eta)$ is locally 3-Sasakian. Example~\ref{3sas1} and
Proposition~\ref{tor-ein} complete the proof.
\end{proof}

\begin{cor}\label{3-sas}
{Let $(M,g,\mathbb Q)$ be a QC structure} on a (4n+3)-dimensional
manifold with positive qc-scalar curvature, $Scal>0$. The next
conditions are equivalent
\begin{enumerate}
\item[i)] {The structure $(M, \frac{16n(n+2)}{Scal} g,
\mathbb{Q})$ is} locally 3-Sasakian; \item[ii)] There exists a
(local) 1-form $\eta$ such that the connection 1-forms of the
Biquard connection vanish on $H$,
$\alpha_i(X)=-d\eta_j(\xi_k,X)=0, X\in H, \quad i,j,k =1,2,3$ and
$Scal$ is constant if $n=1$.
\end{enumerate}
\end{cor}
\begin{proof}
In view of Theorem~\ref{Ein2MO} and Example~\ref{3sas1} it is
sufficient to prove the following Lemma.
\end{proof}
\begin{lemma}
If a QC structure has zero connection one forms restricted to the horizontal space $H$
then it is qc-Einstein, or equivalently, it has zero torsion.
\end{lemma}
\begin{proof}
If $\alpha_i(X)=0$ for $i=1,2,3$ and $X\in H$ then \eqref{sixtyone} together with
\eqref{torha}  yield\\
\centerline{$2\rho_i(X,Y)=-\alpha_i([X,Y])=\alpha_i(T(X,Y))=2\sum_{s=1}^3\alpha_i(\xi_s)\omega_s(X,Y).$}\\
Substitute the latter  into \eqref{rhoh} to conclude considering the
$Sp(n)Sp(1)$-invariant parts of the obtained equalities that
$T^0(X,Y)=U(X,Y)=\alpha_i(\xi_j)=0,\quad
\alpha_i(\xi_i)=-\frac{Scal}{8n(n+2)}$.
\end{proof}

\begin{cor}\label{umb}
Any totally umbilical hypersurface $M$ in a quaternionic-K\"ahler manifold
admits a canonical  qc-Einstein structure with a non-zero scalar
curvature.

\end{cor}
\begin{proof}
Let $M$ be a  hypersurface in the quaternionic-K\"ahler manifold $(\tilde
M, \tilde g)$. With $N$ standing for  the unit normal vector field on $M$
the second fundamental form is given by $II(A,B)=-\tilde g(\tilde D_A
N,B),$ with $A,B\in TM$ and $\tilde D$ being the Levi-Civita connection of
$(\tilde M,\tilde g)$. Since $M$ is a totally umbilical hypersurface of an
Einstein manifold $\tilde M$, taking a suitable trace in the Codazzi
equation \cite[Proposition 4.3]{KoNo},  we find $II(A,B)=-Const\,\tilde
g(A,B)$. Thus, after a homothety  of $\tilde M$ we can assume
$II(A,B)=-\tilde g(A,B)$.

Consider  a local basis $\tilde J_1,\tilde J_2,\tilde J_3$ of the
quaternionic structure of $\tilde M$ satisfying the quaternionic
identities. We define the horizontal distribution $H$ of $M$ to be the
maximal subspace of $TM$ invariant under the action of $\tilde J_1,\tilde
J_2,\tilde J_3$, whose restriction to $M$ will be denoted with $ I_1,I_2,
I_3$. We claim that $(H, I_1, I_2,I_3)$ is a qc-structure on $M$. Defining
$\eta_i(A)= \tilde g(\tilde J_i(N),A)$ and $\xi_i=\tilde J_i N$, a small
calculation shows
\begin{equation*}
d\eta_i(X,Y)=II(X,I_iY)-II(I_iX,Y)=2\tilde g(I_iX,Y), \quad X,Y\in H,
\end{equation*}
Hence, the conditions in the definition of a qc-structure are satisfied.

Let $D$ be the Levi-Civita connection of the restriction $g$ of $\tilde g$
to $M$. Then we have
\begin{equation}\label{sec_fund}
\tilde D_AB=D_AB+II(A,B)N,\qquad A,B\in TM
\end{equation}

Define $\tilde I_i(A)= \tilde J_i(A)_{TM}$ the orthogonal projection on
$TM$, $\ A\in TM$. Since by assumption $\tilde M$ is a
quaternionic-K\"ahler manifold we have $\tilde D \tilde
J_i=-\sigma_j\otimes \tilde J_k + \sigma_k\otimes \tilde J_j$. This
together with \eqref{sec_fund}, after some computation gives (compare with
\eqref{New9})
\begin{equation}\label{der2}\sideremark{New9}
D\tilde I_i=Id\otimes\eta_i-g\otimes\xi_i-\sigma_j\otimes\tilde I_k+
\sigma_k\otimes\tilde I_j.
\end{equation}
Using \eqref{der2} we will show that the torsion of the qc-structure is
zero. The same computation as in the Case~3 in the proof of
Theorem~\ref{Ein2MO} gives
\begin{multline*}
-2g(X,Y) = 2g((D_X\tilde I_i)Y, \xi_i)=2g(D_X\tilde I_iY,\xi_i)\\
=\ -\xi_ig(X,I_iY)
+g([X,I_iY],\xi_i)-g([X,\xi_i],I_iY)-g([I_iY,\xi_i],X)=-(\mathcal{L}_{\xi_i}g)(X,I_iY)-2g(X,Y).
\end{multline*}
Hence $0=(\mathcal{L}_{\xi_i}g)(X,Y)=2T^0_{\xi_i}(X,Y).$  A computation
analogous to the the proof of Theorem~\ref{Ein2MO}, Case~6 gives
\begin{multline}
-2\sigma_3 (\xi_1)\omega_1 (Y,Z)+2\sigma_1 (\xi_1)\omega_3 (Y,Z)= 2g((D_{\xi_1}\tilde I_2)Y,Z)\\
=g((\mathcal L_{\xi_1}I_2)Y,Z)-g((\mathcal
L_{\xi_1}I_2)Z,Y)+\omega_2(I_1Y,Z)+\omega_2(Y,I_1Z).
\end{multline}
From Lemma~\ref{twentytwo} it follows
\begin{multline*}
-2\sigma_3 (\xi_1)\omega_1 (Y,Z)+2\sigma_1 (\xi_1)\omega_3
(Y,Z)=\\
=-4g(I_3\tilde{u}Y,Z)+2d\eta_1(\xi_2,\xi_1)\omega_1(Y,Z) +
(d\eta_1(\xi_2,\xi_3)-d\eta_2(\xi_3,\xi_1)-d\eta_3(\xi_1,\xi_2)-2)\omega_3(Y,Z).
\end{multline*}
Working similarly in the remaining cases we conclude that the traceless
part $u$ of the tensor $\tilde u$ vanishes.
\end{proof}

\section{Conformal transformations of a qc-structure}
Let $h$ be a positive smooth function on a QC manifold $(M, g, \mathbb{Q})$. Let
$\bar\eta=\frac{1}{2h}\eta$ be a conformal deformation of the QC structure $\eta$ (to be
precise we should let $\bar g =\frac{1}{2h} g$ on $H$ and consider $(M, \bar g,
\mathbb{Q})$). We denote the objects related to $\bar\eta$ by overlining the same object
corresponding to $\eta$. Thus,
 $d\bar\eta=-\frac{1}{2h^2}dh\wedge\eta+\frac{1}{2h}d\eta$
and  $\bar g=\frac{1}{2h}g$. The new  triple $\{\bar\xi_1,\bar\xi_2,\bar\xi_3\}$ is
determined by \eqref{bi1}. We have
\begin{equation}\label{New19}
\bar\xi_s=2h\xi_s+I_s\nabla h,\quad s=1,2,3,
\end{equation}
where $\nabla h$ is the horizontal gradient defined by $g(\nabla h,X)=dh(X), X\in
H$.

The horizontal sub-Laplacian and the norm of the horizontal gradient  are defined
respectively by
\begin{equation}\label{e:horizontal laplacian and gradient}
\triangle h\ =\ tr^g_H(\nabla dh)\  = \ \sum_{\alpha=1}^{4n}\nabla dh(e_\alpha,e_\alpha
),\qquad
 |\nabla h|^2\ =\ \sum_{\alpha=1}^{4n}dh(e_\alpha)^2.
\end{equation}
The Biquard connections $\nabla$ and $\bar\nabla$ are connected by a (1,2) tensor S,
$\bar\nabla_AB=\nabla_AB+S_AB, A,B\in\Gamma(TM).$
The condition \eqref{torha} yields\\
\centerline{$g(S_XY,Z)-g(S_YX,Z)=-h^{-1}\sum_{s=1}^3\omega_s(X,Y)dh(I_sZ), \quad X,Y,Z\in H.$}
 From $\bar\nabla\bar g=0$ we get
$g(S_XY,Z)+g(S_XZ,Y)=-h^{-1}dh(X)g(Y,Z), \quad X,Y,Z\in H.$

The last two equations determine $g(S_XY,Z)$ for $X,Y,Z\in H$ due
to the equality\\
\leftline{$g(S_XY,Z)=-(2h)^{-1}\{dh(X)g(Y,Z)-\sum_{s=1}^3dh(I_sX)\omega_s(Y,Z)$}\\
\rightline{$+dh(Y)g(Z,X)+\sum_{s=1}^3dh(I_sY)\omega_s(Z,X)-dh(Z)g(X,Y)+\sum_{s=1}^3dh(I_sZ)\omega_s(X,Y)\}.$}
Using Biquard's Theorem~\ref{biqcon}, we obtain after some calculations that
\begin{multline}\label{New20}
g(\bar T_{\bar\xi_1}X,Y)-2hg(T_{\xi_1}X,Y)-g(S_{\bar\xi_1}X,Y)=\\
-\nabla dh(X,I_1Y)+h^{-1}(dh(I_3X)dh(I_2Y)-dh(I_2X)dh(I_3Y)).
\end{multline}
The identity $d^2=0$ yields
$\nabla dh(X,Y)-\nabla dh(Y,X)=-dh(T(X,Y)).$
Applying \eqref{torha}, we can write
\begin{equation}\label{symdh}
\nabla dh(X,Y)=[\nabla dh]_{[sym]}(X,Y)-\sum_{s=1}^3 dh(\xi_s)\omega_s(X,Y),
\end{equation}
where $[.]_{[sym]}$ denotes the symmetric part of the correspondin (0,2)-tensor.

Decomposing \eqref{New20} into [3] and [-1] parts according to \eqref{New21},
using the properties of the torsion tensor $T_{\xi_i}$ and \eqref{tor} we come to the
next transformation formulas:
\begin{gather}\label{e:T^o conf change}
\overline T^0(X,Y) \ =\ T^0(X,Y)\ +\ h^{-1}[\nabla dh]_{[sym][-1]},\\\label{e:U conf
change}
 \bar U(X,Y) \ =\
U(X,Y)\ +\  (2h)^{-1}[\nabla dh-2h^{-1}dh\otimes dh]_{[3][0]},
\end{gather}
\begin{align*}
 g(S_{\bar\xi_1}X,Y)& \ =\ -\frac{1}{4}\Big[-\nabla
dh(X,I_1Y)+\nabla dh(I_1X,Y) -\nabla dh(I_2X,I_3Y)+\nabla
dh(I_3X,I_2Y)\Bigr]\\
 \hskip.4truein & -\
(2h)^{-1}\Bigl[dh(I_3X)dh(I_2Y)-dh(I_2X)dh(I_3Y)+
dh(I_1X)dh(Y)-dh(X)dh(I_1Y)\Bigr]\\
 \hskip.8truein & \ + \
\frac{1}{4n}\left(-\triangle h+2h^{-1}|\nabla h|^2\right)g(I_1X,Y)
-dh(\xi_3)g(I_2X,Y)+dh(\xi_2)g(I_3X,Y),
\end{align*}
 where $[.]_{[sym][-1]}$ and $[.]_{[3][0]}$ denote the symmetric
 $[-1]$-component and the traceless [3] part of the
 corresponding (0,2) tensors on $H$, respectively. Observe that for $n=1$
 \eqref{e:U conf change} is trivially satisfied.


\noindent Thus, using \eqref{sixtyfour}, we  proved the following Proposition.

\begin{prop}\label{trric}
Let $\overline{\eta}=\frac1{2h}\eta$ be a conformal transformation of a given QC
structure $\eta$. Then the trace-free parts of the corresponding qc-Ricci tensors are
related by the equation
\begin{multline}\label{New24}\sideremark{New24}
Ric_0(X,Y)-{\overline {Ric}}_0(X,Y)\\
= -(2n+2)h^{-1}[\nabla dh]_{[sym][-1]}(X,Y)-(2n+5)h^{-1}\Bigl[\nabla
dh-2h^{-1}dh\otimes dh\Bigr]_{[3][0]}(X,Y).
\end{multline}
 For $n=1,\qquad Ric_0(X,Y)-{\overline {Ric}}_0(X,Y) = - 4h^{-1}[\nabla
dh]_{[sym][-1]}(X,Y)$.
\end{prop}

\noindent In addition, the qc-scalar curvature transforms by the formula \cite{Biq1}
\begin{equation}\label{e:conf change scalar curv}
\overline {Scal} =2h(Scal) -8(n+2)^2h^{-1}|\nabla h|^2 +8(n+2)\triangle h.
\end{equation}

\subsection{Conformal transformations preserving the qc-Einstein condition}
In this section we investigate the question of conformal transformations, which preserve
the qc-Einstein condition. A straightforward consequence of \eqref{New24} is the
following
\begin{prop}\label{con0}\sideremark{con0}
Let $\bar\eta=\frac{1}{2h}\eta$ be a conformal deformation of a given qc-structure $(M,
g, \mathbb{Q})$. Then the trace-free part of the qc-Ricci tensor does not change if and
only if the function $h$ satisfies the differential equations
\begin{align}\label{con01}
& 3(\nabla_Xdh)Y-\sum_{s=1}^3(\nabla_{I_sX}dh)I_sY \ =\ -4\sum_{s=1}^3dh(\xi_s)\omega_s(X,Y),\\
\label{con03} &
(\nabla_Xdh)Y-2h^{-1}dh(X)dh(Y)+\sum_{s=1}^3\left[(\nabla_{I_sX}dh)I_sY-2h^{-1}dh(I_sX)dh(I_sY)\right]
 =\ \lambda g(X,Y), 
\end{align}
for some smooth function $\lambda$ and any $X,Y\in H$.
\end{prop}
Note that for $n=1$ \eqref{con03} is trivially satisfied. Let us fix a qc-normal frame,
cf. definition \ref{d:normal frame}, $\{T_{\alpha},
X_{\alpha}=I_1T_{\alpha},Y_{\alpha}=I_2T_{\alpha},Z_{\alpha}=I_3T_{\alpha},
\xi_1,\xi_2,\xi_3\},  \alpha =1\dots,n$ at a point $p\in M$.
\begin{lemma}\label{vcom01}
If $h$ satisfies \eqref{con01} then we have at $p\in M$ the relations
\begin{equation}\label{e:hHessian of h}
\begin{array}{c}  (I_j\, T_\alpha )\, T_\alpha \, h\ =\ -\  T_\alpha \, (I_j\,
T_\alpha ) \, h\ =\ \xi_j\, h\\
 (I_j\, T_\alpha )\,
(I_i\, T_\alpha ) \, h\ =\  -\ (I_i\, T_\alpha )\, (I_j\, T_\alpha ) \, h\ =\ \xi_k\, h.
\end{array}
\end{equation}
\end{lemma}
\begin{proof}
Working with the fixed qc-normal frame, equation \eqref{con01} gives
$$4T_{\alpha}X_{\alpha}h\, (p)\ -\ [T_{\alpha},X_{\alpha}]h\,(p)\ +\ [Y_{\alpha},Z_{\alpha}]h\, (p)\ =\ -4\xi_1h\, (p).$$
Lemma~\ref{norma} and \eqref{torha} yield $[T_{\alpha},X_{\alpha}]h\,
(p)\ -\ [Y_{\alpha},Z_{\alpha}]h\, (p)\ =\ 0$. Hence, \eqref{e:hHessian of h} follow.
\end{proof}

\subsection{Quaternionic Heisenberg group. Proof of
Theorem \ref{t:einstein preserving}}\label{5.2} 

The proof of Theorem~\ref{t:einstein preserving}  will be presented as separate
Propositions and Lemmas in the rest of the Section, see \eqref{e:final form of h} for the
final formula.
We  use the following model of the quaternionic Heisenberg group  $\QH$. Define $\QH
\ =\ \Hn\times\text {Im}\, \mathbb{H}$ with the group law given by
$
( q', \omega')\ =\ (q_o, \omega_o)\circ(q, \omega)\ =\ (q_o\ +\ q, \omega\ +\ \omega_o\ +
\ 2\ \text {Im}\  q_o\, \bar q),
$
\noindent where $q,\ q_o\in\Hn$ and $\omega, \omega_o\in \text {Im}\, \mathbb{H}$.
In coordinates, with the obvious notation, a basis of left invariant horizontal vector fields $T_{\alpha},
X_{\alpha}=I_1T_{\alpha},Y_{\alpha}=I_2T_{\alpha},Z_{\alpha}=I_3T_{\alpha}, \alpha
=1\dots,n$ is given by
\begin{equation*}\label{qHh}
\begin{aligned}
T_{\alpha}\  =\  \dta {}\ +2x^{\alpha}\dx {}+2y^{\alpha}\dy
{}+2z^{\alpha}\dz {} \,\qquad
  X_{\alpha}\  =\  \dxa {}\ -2t^{\alpha}\dx {}-2z^{\alpha}\dy
{}+2y^{\alpha}\dz {} \,\\
Y_{\alpha}\  =\  \dya {}\ +2z^{\alpha}\dx {}-2t^{\alpha}\dy
{}-2x^{\alpha}\dz {}\,\qquad
Z_{\alpha}\  =\  \dza {}\ -2y^{\alpha}\dx {}+2x^{\alpha}\dy {}-2t^{\alpha}\dz {}\,.
\end{aligned}
\end{equation*}
\noindent The central (vertical) vector fields $\xi_1,\xi_2,\xi_3$ are described as
follows\\
\centerline{$\xi_1=2\dx {}\,\quad \xi_2=2\dy {}\,\quad \xi_3=2\dz {}\,.$}
A small calculation shows the following commutator relations
\begin{equation}\label{e:commutators}
[I_j\, T_\alpha, T_\alpha]\ =\ 2\xi_j\quad\quad\quad  [I_j\, T_\alpha, I_i\,  T_\alpha]\
=\ 2\xi_k.
\end{equation}
The standart 3-contact form
$\tilde\Theta=(\tilde\Theta_1,\ \tilde\Theta_2, \ \tilde\Theta_3)$ is
\begin{equation}\label{e:Heisenbegr ctct forms}
2\tilde\Theta\  =\ \
d\omega \ - \ q' \cdot d\bar q' \ + \ dq'\, \cdot\bar q'.
\end{equation}
\noindent The described horizontal and vertical vector fields are parallel with respect to the
Biquard connection and constitute an orthonormal basis of the tangent space.

We turn to the proof of Theorem \ref{t:einstein preserving}. We start with a Proposition
in which we shall determine the vertical Hessian of $h$.
\begin{prop}\label{vcom01H}
If $h$ satisfies \eqref{con01} on $\QH$ then we have the relations
\begin{equation}\label{vcomHH}
 \xi_1^2(h)=\xi_2^2(h)=\xi_3^2(h)=8\mu_o, \quad
\xi_i\xi_j(h)=0, \quad i\not=j=1,2,3,
\end{equation}
where $\mu_o>0$ is a constant. In particular,
\begin{equation}\label{e:form of h}
h(q,\omega)\ =\ g(q)\ +\ \mu_o\, \bigl [ \ (x\ +\ x_o(q)\,)^2\ +\ (y\ +\ y_o(q)\,)^2\ +\
(z\ +\ z_o(q)\,)^2\ \bigr  ]
\end{equation}
for some real valued functions $g,\ x_o\, \ y_o$ and $z_o$ on $\Hn$. Furthermore we have
\begin{equation*}\label{vcomxx}
T_{\alpha}Z_{\alpha}X^2_{\alpha}(h)=T_{\alpha}Z_{\alpha}Y^2_{\alpha}(h)=0, \quad
T^2_{\alpha}\xi_j(h)=0.
\end{equation*}
\end{prop}
\begin{proof}
Equations \eqref{e:hHessian of h} and \eqref{e:commutators} yield the next sequence of equalities
\begin{multline*}
2\xi_i\xi_j\ h \ =\ -2T_{\alpha}\ (I_i\, T_{\alpha})\ \xi_j\ h \
= \ -2T_{\alpha}\ \xi_j \  (I_i \, T_{\alpha})\  h
= \  \ T_{\alpha}\ [T_{\alpha}, I_j \, T_{\alpha} ]\ (I_i T_{\alpha})\
 h \ \\= \  \ T^2_{\alpha}\,  ( I_j \, T_{\alpha} )\ (I_i\,  T_{\alpha} ) \ h \ -
\ 2T_{\alpha}\ ( I_j\,  T_{\alpha})\ T_{\alpha} \ ( I_i
T_{\alpha})\ h
= \ \ T^2_{\alpha}\ \xi_k \ h \ - \ \xi_i\ \xi_j\ h.
\end{multline*}
Hence, $ 3\xi_i\xi_j\ h=T^2_{\alpha}\ \xi_k\ h$. Similarly, interchanging the roles of
$i$ and $j$ together with $\{ I_iT_\alpha, I_j T_\alpha \}= 0$ we find $
3\xi_j\xi_i(h)=-T^2_{\alpha}\ \xi_k(h)$. Consequently
$\xi_i\xi_j\ h \ = \ T^2_{\alpha}\ \xi_k\ h \ = \ 0.$
An analogous calculation shows that $\xi_i\xi_k\ h \ = \ 0.$
 Furthermore, we have
\begin{multline*}
2\xi_1^2(h)=2X_{\alpha}T_{\alpha}\xi_1(h)=2X_{\alpha}\xi_1T_{\alpha}(h)=
-X_{\alpha}[Y_{\alpha}Z_{\alpha}]T_{\alpha}(h)\\=
X_{\alpha}Z_{\alpha}Y_{\alpha}T_{\alpha}(h)-
X_{\alpha}Y_{\alpha}Z_{\alpha}T_{\alpha}(h)=\xi_2^2(h)+\xi_3^2(h).\hskip.5truein
\end{multline*}
We derive similarly
$2\xi^2_2(h)=\xi_1^2(h)+\xi_3^2(h), \quad
2\xi_3^2(h)=\xi_2^2(h)+\xi_1^2(h).$ Therefore
$\xi_1^2(h)=\xi_2^2(h)=\xi_3^2(h), \quad
\xi_i^3(h)=\xi_i\xi_j^2(h)=0, \quad i\not=j=1,2,3$ which proves part of  \eqref{vcomHH}.

Next we prove that the common value of the second derivatives is a constant. For this we
differentiate the equation  $T^2_\alpha\, \xi_k h = 0$ with respect to $I_k T_\alpha$
from where, \eqref{e:hHessian of h} and \eqref{e:commutators}, we get
\begin{multline*}
0  \ =\ \xi_k \ (I_k \, T_\alpha) \ T^2_\alpha\ h\ =\ \xi_k\ T_\alpha\ (I_k \, T_\alpha)
\ T_\alpha\ h\ + \ \xi_k\ [I_k\,
T_\alpha,\ T_\alpha]\ h\\
 =\ T_\alpha\  \xi_k^2\ h \ +\ 2\ T_\alpha\ \xi_k^2 \ h\ = \ 3\ T_\alpha\  \xi_k^2\ h  .
\end{multline*}
In order to see the vanishing of $(I_i\, T_\alpha)\  \xi_k^2\ h$ we shall need
\begin{equation}\label{e:IiTa^2xih }
(I_i \, T_\alpha)^2\, \xi_k h = 0.
 \end{equation}
\noindent The latter can be seen by the following calculation.
\begin{multline*}
2\xi_i\xi_j\ h \ =\ 2\xi_i\ (I_i\, T_\alpha) \ (I_k\, T_{\alpha}) \ h \ = \ 2(I_i\,
T_\alpha)\ \xi_i \  (I_k \, T_{\alpha})\  h  \ = \  \ T_{\alpha}\ [T_{\alpha}, I_j
\, T_{\alpha} ]\ (I_i T_{\alpha})\
 h\ \\ = \ \ (I_i\, T_\alpha)^2\ T_\alpha  ( I_k \, T_{\alpha} ) \ h \
 -\ (I_i\, T_\alpha)\ T_\alpha\ (I_i\, T_\alpha)\ (I_k\,
 T_\alpha)\ h
= \ - \ (I_i\, T_\alpha)^2\ \xi_k \ h \ - \ \ \xi_i\ \xi_j\ h,
\end{multline*}
from where $0= 3\xi_i\xi_j h = -  (I_i\, T_\alpha)^2\xi_k h.$
Differentiate $(I_i \, T_\alpha)^2\, \xi_k h = 0$ with respect to $I_j T_\alpha$ to get
\begin{multline*}
0  \ =\ \xi_k \ (I_j \, T_\alpha)\ (I_i \, T_\alpha)^2\, \xi_k\ h =\ \xi_k\
(I_i\,T_\alpha)\ (I_k \, T_\alpha) \ T_\alpha\ h\ + \ \xi_k\ [I_j\,T_\alpha,
I_i\, T_\alpha]\ (I_i\,T_\alpha)\ h\\
 =\ (I_i\, T_\alpha)\  \xi_k^2\ h \ +\ 2\ (I_i\, T_\alpha)\ \xi_k^2 \ h\ = \ 3\ (I_i\, T_\alpha)\  \xi_k^2\ h
 .\hskip1truein
\end{multline*}
We proved the vanishing of all derivatives of the common value of $\xi^2_j h$, i.e., this
common value is a constant, which we denote by $8\mu_o$.  Let us note that $\mu_o>0$
follows easily from the fact that $h>0$ since $g$ is independent of $x,\ y$ and $z$.

The rest equalities of the proposition follow easily from \eqref{e:hHessian of h} and
\eqref{vcomHH}.
\end{proof}
In view of Proposition~\ref{vcom01H}, we  define $h\ =\ g\ +\ \mu_o\, f$, where
\begin{equation}\label{e:f def}
f\ =\ (x+x_o(q))^2\ +\  (y+y_o(q))^2\ +\  (z+z_o(q))^2.
\end{equation}
The following simple Lemma is one of the keys to integrating our system.
\begin{lemma}\label{l:key integrating}
Let $X$ and $Y$ be two parallel horizontal vectors
\begin{enumerate}[a)]
\item If $\omega_s(X, Y)\ =\ 0$  for $s=1,2,3$  then
\begin{equation}\label{e:omegaXY=0}
4 XY h - 2h^{-1} \Bigl[dh(X)\ dh (Y) +\sum_{s=1}^3 dh (I_sX)\ dh (I_sY)\Bigr]
=\ \lambda\, g(X, Y).
\end{equation}
\item If $g(X,Y)\ =\ 0$ then
\begin{equation}\label{e:gXY=0}
2\, XY h\ -\ h^{-1} \, \{\,  dh(X)\ dh (Y)\ +\ \sum_{s=1}^3 dh (I_sX)\ dh (I_sY)\, \} \
=\ 2\ \sum_{s=1}^3 \big \{ \ (\xi_s h)\ \omega_s(X, Y)\  \big \}.
\end{equation}
\item If $g(X,Y)\ =\ \omega_s(X, Y)\ =\ 0$  for $s=1,2,3$  we have for
any $j\in \{1,2,3\}$
\begin{gather}\label{e:XYh}
XY \, (\xi_j  h)\ =\ 0, \qquad 8\, XY h \  = \ \mu_o \, \big \{  (X\xi_j f) \,
(Y\xi_j f) \ + \ \sum_{s=1}^3 (I_sX\, \xi_j f ) \, (I_sY\, \xi_j f)  \, \big \} .
\end{gather}
\end{enumerate}
\end{lemma}
\begin{proof}
The equation of a) and b) are obtained by adding \eqref{con01} and \eqref{con03}. Let us
prove part c).  From \eqref{con01} and \eqref{con03} taking any two horizontal vectors
satisfying $g(X,Y)=\omega_s(X, Y)=0$, we obtain \hspace{4mm}
$2h\nabla dh (X, Y)\ =\ dh(X)\ dh (Y)\ +\ \sum_{s=1}^3 dh (I_sX)\ dh (I_sY).$

\noindent If $X, \ Y$ are also parallel, differentiate along $\xi_j$ twice to get consequently
\begin{multline}\label{e:xi main integrating equality1}
2\xi_j h\, XY h \  + \ 2h XY\xi_j h
 = \ (X\xi_j h)\, (Yh)\  + \
\sum_{s=1}^3 [ (I_sX\, \xi_j h)\, (I_sY\, h)\  + \ (I_s X\, \xi_j h)\, (I_sY\, h) ],\\
2\xi^2_j  h\, XY h \  + \ 4h XY\xi_j h\ = \ 2 \big\{\,  (X\xi_j h)\, (Y\xi_jh)\ + \
\sum_{s=1}^3  (I_sX\, \xi_j h)\, (I_sY\, \xi_j h) \, \big\}.
\end{multline}
\noindent Differentiate three times along $\xi_j$ and use $\xi_j ^2
h=$const, cf. \eqref{vcomHH} to get
$2\ (\xi_j^2 h)\ XY \, (\xi_j  h)\ =\ 0,$
from where the first equality in \eqref{e:XYh} follows. With this information the second line in
\eqref{e:xi main
integrating equality1} reduces to the second equality in \eqref{e:XYh}.
\end{proof}
In order to see that after a suitable translation the functions
$ x_o\, \ y_o$ and $z_o$ can be made equal to zero we prove the following proposition.
\begin{prop}\label{p:HHVh}\sideremark{p:HHVh}
If $h$ satisfies \eqref{con01} and \eqref{con03} on $\QH$ then we have
\begin{enumerate}[a)]
\item For $s\in \{1,2,3\}$ and $i,j,k$ a cyclic permutation of
$1,2,3$
\begin{equation}\label{e:HHVh}
\begin{aligned}
T_\alpha\ T_\beta \  (\xi_s h)\  & =\ (I_i T_\alpha)\ (I_iT_\beta)
\ (\xi_s h)\   =\  0 \quad \forall  \ \alpha,\ \beta\\
(I_i T_\alpha)\ T_\beta \ (\xi_s h) \ & =\ (I_i T_\alpha)\
(I_iT_\beta) \ (\xi_s h)\  =\ 0, \ \alpha\not=\beta\\
 (I_j\, T_\alpha )\, T_\alpha \ (\xi_s \, h) \ & =\ -\  T_\alpha \,
(I_j\, T_\alpha ) \ (\xi_s h) \ =\ 8\ \delta_{sj}\ \mu_o\\
 (I_j\, T_\alpha )\, (I_i\, T_\alpha )\ (\xi_s h) \  & =\  -\ (I_i\,
T_\alpha )\, (I_j\, T_\alpha ) \ (\xi_s h) \ =\ 8\ \delta_{sk} \ \mu_o\  ,
\end{aligned}
\end{equation}
\noindent i.e., the horizontal Hessian of a vertical derivative of $h$ is determined
completely. \item There is a point $(q_o, \omega_o)\in\QH$, $q_o=(q^1_o, \, q^2_o, \dots,
q^n_o)\in\Hn$ and $\omega=ix_o+jy_o+kz_o\in Im (\mathbb{H})$, such that, \\
\centerline{$
ix_o(q)\ +\ jy_o(q)\ +\ k z_o(q)\ =\ w_o\  +\ 2\ \text {Im}\  q_o \, \bar q.
$}
\end{enumerate}
\end{prop}
\begin{proof}
a) Taking $\alpha\not=\beta$ and using $X=T_\alpha$ and $Y=T_\beta$ in \eqref{e:XYh}, we obtain
$
T_\alpha T_\beta \xi_s h= 0, \quad \alpha\not=\beta.
$
When $\alpha=\beta$ the same equality holds by \eqref{e:IiTa^2xih }.

The vanishing of the other derivatives can be obtained similarly. Finally, the rest of
the second derivatives can be determined from \eqref{e:hHessian of h}.

b)  From the identities in \eqref{e:HHVh} all second derivatives of $x_o, \ y_o$ and
$z_o$ vanish. Thus  $x_o, \ y_o$ and $z_o$ are linear function.  The fact that the
coefficients are related as required amounts to the following system
\begin{align*}
&  T_\alpha \, x_o\ =\ Z_\alpha \, y_o\ =\ - Y_\alpha\, z_o, \qquad
 X_\alpha\, x_o\ =\ Y_\alpha\, y_o\ =\ Z_\alpha\, z_o\\
&  Y_\alpha\, x_o\ =\ X_\alpha\, y_o\ =\ - T_\alpha\, z_o,\qquad
  Z_\alpha\, x_o \ =\ -T_\alpha\, y_o\ =\ - X_\alpha\, z_o.
\end{align*}
>From \eqref{e:form of h} we have \hspace{4mm} $\xi_1 h\ =\ 4\mu_o
(x+x_o(q)), \quad  \xi_2 h\ =\ 4\mu_o (y+y_o(q)), \quad  \xi_3 h\ =\
4\mu_o (z+z_o(q)).$ Therefore, the above system is equivalent to
\begin{align*}
& T_\alpha \, \xi_1\ h\ =\ Z_\alpha \, \xi_2\ h\ =\ - Y_\alpha\, \xi_3\ h, \qquad
  X_\alpha\, \xi_1\ h\ =\ Y_\alpha\, \xi_2\ h\ =\ Z_\alpha\, \xi_3\ h\\
&  Y_\alpha\, \xi_1\ h \ =\ X_\alpha\, \xi_2\ h\ =\ - T_\alpha\, \xi_3\ h\qquad
 Z_\alpha\, \xi_1\ h \ =\ -T_\alpha\, \xi_2\ h\ =\ - X_\alpha\, \xi_3\ h.
\end{align*}
Let us prove the first line.  Denote $ a\ =\ T_\alpha \, \xi_1\ h,\quad b\ =\ Z_\alpha \,
\xi_2\ h,\quad c\ =\ - Y_\alpha\, \xi_3\ h.$ From \eqref{e:hHessian of h} and \eqref{e:commutators}
it follows
\begin{align*}
& a\ = \ T_\alpha Z_\alpha  Y_\alpha  \ h\ =\ Z_\alpha  T_\alpha Y_\alpha   \ h \ +
\ [T_\alpha ,Z_\alpha  ]Y_\alpha   \ h\ =\ -b\ +\ 2c\\
& b\ = \ Z_\alpha  Y_\alpha  T_\alpha \ h\ =\ Y_\alpha  Z_\alpha  T_\alpha  \ h \ +
\ [Z_\alpha  ,Y_\alpha  ] T_\alpha  \ h\ =\ -c\ +\ 2a\\
& c\ = \ Y_\alpha  T_\alpha Z_\alpha  \ h\ =\ T_\alpha Y_\alpha Z_\alpha   \ h \ +\
[Y_\alpha  ,T_\alpha ]Z_\alpha   \ h\ =\ -a\ +\ 2b,
\end{align*}
which implies $a=b=c$. The rest of the identities of the system can be obtained
analogously.
\end{proof}
So far we have proved that if $h$ satisfies the system \eqref{con01} and \eqref{con03} on
$\QH $ then, in view of the translation invariance of the system, after a suitable
translation we have
\[
h(q, \omega) \ =\ g(q)\ +\ \mu_o\, (x^2\ +\ y^2\ +\ z^2).
\]
\begin{prop}\label{p:find g}
If $h$ satisfies the system \eqref{con01} and \eqref{con03} on $\QH $ then after a
suitable translation we have $$g(q)\ =\ (b\ +\ 1\ +\ \sqrt \mu_o \ |q|^2)^2, \quad\quad
b\ + \ 1\ >\ 0.$$
\end{prop}
\begin{proof}
Notice that  $\xi_1 h\ =\ 4\mu_ox,  \xi_2 h\ =\ 4\mu_o y,
\xi_3 h\ =\ 4\mu_o z.$  With this equations \eqref{e:hHessian of h} become
\begin{equation}\label{e:vcomh final}
\begin{array}{c}
T_{\alpha}X_{\alpha}(h)=Y_{\alpha}Z_{\alpha}(h)=-X_{\alpha}T_{\alpha}(h)=-Z_{\alpha}Y_{\alpha}(h)=-4\mu_ox,\\
T_{\alpha}Y_{\alpha}(h)=Z_{\alpha}X_{\alpha}(h)=-Y_{\alpha}T_{\alpha}(h)=-X_{\alpha}Z_{\alpha}(h)=-4\mu_o y,\\
T_{\alpha}Z_{\alpha}(h)=X_{\alpha}Y_{\alpha}(h)=-Z_{\alpha}T_{\alpha}(h)=-Y_{\alpha}X_{\alpha}(h)=-4\mu_o
z.
\end{array}
\end{equation}
Let us also write explicitly some of the derivatives of $f$, which shall be used to
express the derivatives of $g$ by the derivatives of $h$. For all $\alpha$ and $\beta$ we
have
\begin{align*}
& T_\beta f\ =\ 4\, (x^\beta x \ +\ y^\beta y\ +\ z^\beta z),\quad\quad\quad \quad\quad\quad \qquad
X_\beta f\ =\ 4\, (-t^\beta x\ -\ z^\beta y\ +\ y^\beta
z),\\
& Y_\beta f\ =\ 4\, (z^\beta x \ -\ t^\beta y\ -\ x^\beta z),\quad\quad\quad \quad\quad\quad \qquad
Z_\beta f\ =\ 4\, (-y^\beta x\ +\ x^\beta y\ -\ t^\beta
z),\\
& T_\alpha T_\beta f\ =\ 8\, (x^\alpha x^\beta\ +\ y^\alpha y^\beta\ +\ z^\alpha
z^\beta),\quad\quad \quad \qquad X_\alpha X_\beta f\ =\ 8\, (t^\alpha t^\beta\ +\ z^\alpha z^\beta\
+\ y^\alpha y^\beta),\\
& Y_\alpha Y_\beta f\ =\ 8\, (z^\alpha z^\beta\ +\ t^\alpha t^\beta\ +\ x^\alpha
x^\beta),\quad\quad \quad \qquad Z_\alpha Z_\beta f\ =\ 8\, (y^\alpha y^\beta\ +\ x^\alpha x^\beta\
+\ t^\alpha t^\beta),\\
& T_\alpha X_\beta f = -4\delta_{\alpha\beta}x +8(-x^\alpha
t^\beta -y^\alpha z^\beta+ z^\alpha y^\beta), \quad T_\alpha Y_\beta f = -4\delta_{\alpha\beta}y+ 8(x^\alpha
z^\beta-y^\alpha t^\beta  - z^\alpha x^\beta),\\
& T_\alpha Z_\beta f = -4\delta_{\alpha\beta}z + 8(-x^\alpha
y^\beta +  y^\alpha x^\beta  -  z^\alpha t^\beta), \quad
 X_\alpha T_\beta f = 4\delta_{\alpha\beta}x + 8(-t^\alpha
x^\beta - z^\alpha y^\beta +y^\alpha z^\beta),\\
& X_\alpha Y_\beta f = -4\delta_{\alpha\beta}z + 8(-t^\alpha
z^\beta + z^\alpha t^\beta - y^\alpha x^\beta),\quad
 X_\alpha Z_\beta f = 4\delta_{\alpha\beta}y+ 8(t^\alpha y^\beta- z^\alpha
x^\beta + y^\alpha t^\beta).
\end{align*}
\noindent From the above formulas we see that the fifth order horizontal derivatives of
$f$ vanish. In particular the fifth order derivatives of $h$ and $g$ coincide.

Taking $X=Y=T_\alpha$ in \eqref{e:omegaXY=0} we obtain
\begin{equation}\label{e:TTh}
4T^2_\alpha h\ -\ 2h^{-1} \{ (T_\alpha h)^2+(X_\alpha h)^2+(Y_\alpha h)^2+(Z_\alpha
h)^2\}=\lambda.
\end{equation}
\noindent Using in the same manner $X_\alpha$, $Y_\alpha$ and $Z_\alpha$ we see the
equality of the second derivatives
\begin{equation}\label{e:T^2h etc}
T^2_\alpha h\ =\ X^2_\alpha h\ =\ Y^2_\alpha h\ =\ Z^2_\alpha h.
\end{equation}
Therefore, using \eqref{e:hHessian of h} and \eqref{vcomHH}, we have
$T^3_\alpha h  =T_\alpha X^2_\alpha h =X_\alpha T_\alpha
 X_\alpha h +[T_\alpha, X_\alpha]X_\alpha h =-3 X_\alpha \xi_1 h
 = 24\mu_o t^\alpha$
and thus
$ T^4_\alpha h\ =\ 24\mu_o.$
In the same fashion we conclude
\begin{multline}\label{e:T^3h etc}
T^3_\alpha h \ =\ 24\mu_o t^\alpha \qquad X^3_\alpha h \ =\ 24\mu_o x^\alpha \qquad
Y^3_\alpha h \ =\ 24\mu_o y^\alpha \qquad Z^3_\alpha h \ =\ 24\mu_o z^\alpha.
\end{multline}
\noindent Similarly, taking $X=T_\alpha$, $Y=X_\beta$ and $j=1$ in \eqref{e:XYh} we find
\[
\ T_\alpha X_\beta h \  =\ 8\mu_o\ ( -x^\alpha t^\beta \ +\ t^\alpha x^\beta\ -\ y^\alpha
z^\beta\ +\ z^\alpha y^\beta  ), \quad\quad \alpha\not=\beta.
\]
Plugging $X=T_\alpha$, $Y=T_\beta$ with $\alpha\not=\beta$ in \eqref{e:XYh} we obtain
\[
T_\alpha T_\beta h \ =\ X_\alpha h X_\beta h\ =\ Y_\alpha h Y_\beta h \ =\ Z_\alpha h
Z_\beta h\ =\ 8\mu_o(t^\alpha t^\beta\ +\ x^\alpha x^\beta \ +\ y^\alpha y^\beta \ +\
z^\alpha z^\beta).
\]
\noindent The other mixed second order derivatives when $\alpha\not=\beta$ can be
obtained by taking suitable $X$ and $Y$. In view of the formulas for the derivatives of
$f$ and \eqref{e:vcomh final}, we conclude
\begin{equation}\label{e:TalphaTbeta g}
\begin{aligned}
T_\alpha X_\beta\, g \ =\ 8\mu_o t^\alpha x^\beta, \quad T_\alpha Y_\beta\, g \ =\ 8\mu_o
t^\alpha y^\beta, \quad T_\alpha Z_\beta\, g
\ =\ 8\mu_o t^\alpha z^\beta\\
 X_\alpha Y_\beta\, g \ =\ 8\mu_o x^\alpha y^\beta, \quad
X_\alpha Z_\beta\, g \ =\ 8\mu_o x^\alpha z^\beta, \quad Y_\alpha Z_\beta\, g \ =\ 8\mu_o
y^\alpha z^\beta,\\
T_\alpha T_\beta g\ = 8\mu_o t^\alpha t^\beta, \quad X_\alpha
X_\beta g\ = 8\mu_o x^\alpha x^\beta,\quad
 Y_\alpha Y_\beta g\ = 8\mu_o y^\alpha y^\beta, \quad
Z_\alpha Z_\beta g\ = 8\mu_o z^\alpha z^\beta,\\
T_\alpha X_\alpha\, g \ =\ 8\mu_o t^\alpha x^\alpha, \quad T_\alpha Y_\alpha\, g \ =\
8\mu_o t^\alpha y^\alpha, \quad T_\alpha
Z_\alpha\, g \ =\ 8\mu_o t^\alpha z^\alpha\\
 X_\alpha Y_\alpha\, g \ =\ 8\mu_o x^\alpha y^\alpha, \quad
X_\alpha Z_\alpha\, g \ =\ 8\mu_o x^\alpha z^\alpha, \quad Y_\alpha Z_\alpha\, g \ =\
8\mu_o y^\alpha z^\alpha.
\end{aligned}
\end{equation}
%
%
\noindent A consequences of the considerations so far is the
fact that all second order derivative are quadratic functions
  of the variables from the first layer, except the pure (unmixed) second
derivatives, in which case we know \eqref{e:T^2h etc}
  and \eqref{e:T^3h etc}.
   It is easy to see then that
the fifth order horizontal derivatives of $h$ vanish. With the information so far after a
small argument we can assert that $g$ is a polynomial of degree 4 without terms of degree
3, and of the form
\[
g=\mu_o \sum_{\alpha=1}^n (t^4_\alpha\ +\ x^4_\alpha\ +\ y^4_\alpha\ +\ z^4_\alpha)\ +
p_2,
\]
where $p_2$ is  a polynomial of degree two. Furthermore, the mixed second order
derivatives of $g$ are determined, while the pure second order derivatives are equal. The
latter follows from \eqref{e:TTh} taking $q=0, \ \omega=0$.  Let us see that there are no
terms of degree one on $p_2$. Taking $X=T_\alpha$, $Y=T_\beta$, $\alpha\not=\beta$ and
$j=1$ in \eqref{e:xi main integrating equality1} we find
\begin{multline*}
(4\mu_o x)\ \big \{ \ 4T_\alpha T_\beta g \ +\ 32\mu_o(x^\alpha
x^\beta \ +\ y^\alpha y^\beta\ +\ z^\alpha z^\beta ) \ \big \}\\
 =\ 2 \ \Big \{  (8x^\alpha)( T_\beta g+ 4\mu_o (x^\beta x \ +\ y^\beta y\ + \ z^\beta z))
 \ +\
 (-8t^\alpha)( X_\beta g+ 4\mu_o (-t^\beta x \  - \ z^\beta y\ + \ y^\beta z)) \\
 +\ (8z^\alpha)( Y_\beta g+ 4\mu_o (z^\beta x \ -\ t^\beta y\ - \
x^\beta z)) \
 +\ (-8y^\alpha)( Z_\beta g+ 4\mu_o (-y^\beta x \ +\
x^\beta y\ - \ t^\beta z)) \\
+\ (8x^\beta)( T_\alpha g+ 4\mu_o (x^\alpha x \ +\ y^\alpha y\ + \ z^\alpha z))
 \
  +\ (-8t^\beta)( X_\alpha g+ 4\mu_o (-t^\alpha x \  - \ z^\alpha y\ + \
y^\alpha z)) \\
 +\  (8z^\beta)( Y_\alpha g+ 4\mu_o (z^\alpha x \ -\
t^\alpha y\ - \ x^\alpha z)) +\  (-8y^\beta)( Z_\alpha g+ 4\mu_o (-y^\alpha x \ +\
x^\alpha y\ - \ t^\alpha z))
  \Big \}\\
  =\ 16\ \big (x^\alpha T_\beta g \ +\ x^\beta T_\alpha g
  \  -\ t^\alpha X_\beta g \ -\ t^\beta X_\alpha g
  \ +\ z^\alpha Y_\beta g \ +\ z^\beta Y_\alpha g
  \  -\ y^\alpha Z_\beta g \ -\ y^\beta Z_\alpha g
\big )
\\
 +\ 128 \mu_o\ x\ (t^\alpha t^\beta\ +\ x^\alpha
  x^\beta \ +\ t^\alpha t^\beta\ +\ x^\alpha
  x^\beta \ )
\end{multline*}
Taking into account \eqref{e:TalphaTbeta g} we proved
\begin{multline*}
x^\alpha T_\beta g  + x^\beta T_\alpha g
    - t^\alpha X_\beta g  - t^\beta X_\alpha g
   + z^\alpha Y_\beta g  +z^\beta Y_\alpha g
    - y^\alpha Z_\beta g  - y^\beta Z_\alpha g  = 0,  \quad \alpha\not=\beta.
\end{multline*}
\noindent Comparing coefficients  in front of the linear terms implies that $g$ has no
first order terms.
Thus, we can assert that $g$ can be written in the following form
\begin{equation}\label{e:near end form of g}
g\ =\ \big ( 1\ +\ \sqrt {\mu_o}\, |q|^2 \big )^2 \ +\ 2a\, |q|^2 \ +\ b.
\end{equation}
Hence,
$
h\ =\ \big ( 1\ +\ \sqrt {\mu_o}\, |q|^2 \big )^2 \ +\ a\, |q|^2 \ +\ b \ +\ \mu_o\,
(x^2\ +\ y^2\ +\ z^2).
$  Taking $X=T_\alpha$, $Y=T_\beta$ in
\eqref{e:gXY=0} we obtain
$
16\mu_0\, (1\ +\ b) \ =\ 4\, (a\ +\ 2\,  \sqrt{\mu_o})^2.
$
Therefore,
\begin{multline*}
g\ =\ \mu_o \, |q|^4 \ + \ (a\ +\ 2\sqrt {\mu_o})\ \sqrt{\mu_o}\, |q|^2\ +\ b\ + \ 1\ =\
2\sqrt
{b\ +\ 1}\, |q|^2\ +\  b\ +\ 1
=\ \big ( b\ +\ 1\ +\ \sqrt{\mu_o} \, |q|^2 \big )^2.
\end{multline*}
 In turn
the formula for $h$ becomes
\begin{equation}\label{e:final form of h} h \ =\ \big ( b\
+\ 1\ +\ \sqrt{\mu_o} \, |q|^2 \big )^2\ +\ \mu_o\, (x^2\ +\ y^2\ +\ z^2). \end{equation}
Setting
$c\ =\ (b\ +\ 1)^2 \quad\text{ and }\quad \nu\ =\ \frac {\sqrt{\mu_o}
}{1+b}\ >\ 0$ the solution takes the  form \\
$h \ =\ c\ \Big [  \big ( 1\ +\
\nu\, |q|^2 \big )^2\  +\ \nu^2\, (x^2\ +\ y^2\ +\ z^2)\Big ],$ which completes the
proof of Theorem~\ref{t:einstein preserving}.
\end{proof}
 Let us note that the final conclusion can be reached also using the fact
that a qc-Einstein structure has necessarily constant scalar curvature by Theorem
\ref{maincon}, together with the result of \cite{GV} identifying all partially symmetric
solutions of the Yamabe equation on $\QH$, i.e., of the equation
\[
\sum_{\alpha = 1}^n \big ( T^2_\alpha u \ +\ X^2_\alpha u \ +\ Y^2_\alpha u \ +\
Z^2_\alpha u \big )\ =\ -u^{\frac {Q+2}{Q-2}}.
\]
The fact that we are dealing with such a solution follows from \eqref{e:near end form of
g}. The current solution depends on one more parameter as the scalar curvature can be an
arbitrary constant. This constant will appear in the argument of \cite{GV} by first using
scalings  to reduce to a fixed scalar curvature one for example.
\section{Special functions and pseudo-Einstein quaternionic contact
  structures}

Considering only the [3]-component of the Einstein tensor of the Biquard connection due
to Theorem~\ref{sixtyseven} and by analogy with the CR-case \cite{L1} , it seems useful
to give the following Definition.
\begin{dfn}\label{d:qc-pseudo_einstein}
Let $(M, g, \mathbb{Q})$ be a quaternionic contact manifold of dimension bigger than 7.
We call $M$ \emph{qc-pseudo-Einstein} if the trace-free part of the [3]-component of the
qc-Einstein tensor vanishes.
\end{dfn}
\noindent Observe that for $n=1$ any QC structure is qc-pseudo-Einstein. According to
Theorem~\ref{sixtyseven} $(M,g,\mathbb{Q})$ is quaternionic qc-pseudo-Einstein exactly
when the trace-free part of the [3]-component of the torsion vanishes, $U=0$.
Proposition~\ref{trric} yields the following claim.
\begin{prop}\label{pseudoeinst}
Let $\bar\eta=u\eta$ be a conformal deformation of a given qc-structure. Then the
trace-free part of the $[3]$ component of the qc-Ricci tensor (i.e. $U$ ) is preserved if
and only if the function $u$ satisfies the differential equations
\begin{equation}\label{vonps1}
(\nabla_Xdu)Y+(\nabla_{I_1X}du)I_1Y+(\nabla_{I_2X}du)I_2Y+(\nabla_{I_3X}du)I_3Y=\frac1n\Delta
u\,g(X,Y).
\end{equation}
In particular, the qc-pseudo-Einstein condition persists under conformal deformation
$\tilde\eta=u\eta$ exactly when the function $u$ satisfies \eqref{vonps1}.
\end{prop}
 \begin{proof}
Definning $h=u^{-1}$ a small calculation shows
\begin{equation}\label{transf}
\nabla dh-2h^{-1}dh\otimes dh=u^{-1}\nabla du.
\end{equation}
Inserting \eqref{transf} into \eqref{con03} shows \eqref{vonps1}.
\end{proof}

\noindent Our next goal is to investigate solutions to
\eqref{vonps1}. We shall find geometrically defined functions, which
are solutions of \eqref{vonps1}.

\subsection{Quaternionic pluriharmonic functions}

We start with some analysis on the quaternion space $\mathbb{H}^n$.
\subsubsection{Pluriharmonic functions in $\Hn$} \label{ss:Pluriharmonic functions in
Hn} Let $\mathbb{H}$ be the four-dimensional real associative algebra of the quaternions.
The elements of $\mathbb{H}$ are of the form $q=t+ix+jy+kz$, where $t, x, y,
z\in\mathbb{R}$ and $i,\ j, \ k$ are the basic quaternions satisfying the multiplication
rules
$
i^2=j^2=k^2 = -1 \text { and } ijk =-1.
$
For a quaternion $q$ we define its conjugate $\bar q= t-ix-jy-kz$, and real and imaginary
parts, correspondingly, by
$\Re {q}\ =\ t \ \text { and }\ \Im q\ =xi+yj+zk.
$
 The most important operator for us is
the Dirac-Feuter operator $\dbar \ =\ \dt {}\ + \ i\dx {}\ + \ j\dy {}\ +\ k\dz {}$,
i.e.,
$
\dbar  F = \dt    F+ \ i\dx {F}\ + \ j\dy {F}\ +\ k\dz {F}
$
and in addition
\[
\mathcal{D}  F\ =\ \dt {F}\ - \ i\dx {F}\ - \ j\dy {F}\ -\ k\dz {F}.
\]
Note that if $F$ is a  quaternionic valued function due to the non-commutativity of the
multiplication the above expression is not the same as $F\dbar {}\ \overset {def}{=}\ \dt
{F}\ + \ \dx {F}i\ + \ \dy {F}j\ +\ \dz {F}k$. Also, when conjugating $\overline {\dbar
F}\neq \mathcal{D}\overline F$.
\begin{dfn}
A function $F\ :\ \mathbb{H}\rightarrow\ \mathbb{H}$, which is continuously
differentiable when regarded as a function of $\mathbb{R}^4$ into $ \mathbb{R}^4$ is
called \emph{quaternionic anti-regular} (\emph{quaternionic regular}), or just
\emph{anti-regular} (\emph{regular}) for short, if $\mathcal{D} F=0$ ($\dbar F\ =\ 0$).
\end{dfn}
These functions were introduced by Fueter \cite{F}. The reader can consult the paper of
A. Sudbery \cite{S} for the basics of the quaternionic analysis on $\mathbb{H}$.  Let us
note explicitly one of the most striking differences between complex and quaternionic
analysis. As it is well known the theory of functions of a complex variable $z$ is
equivalent to the theory of power series of $z$. In the quaternionic case, each of the
coordinates $t,\  x,\  y$ and $z$ can be written as a polynomial in $q$, see eq. (3.1) of
\cite{S}, and hence the theory of power series of $q$ is just the theory of real analytic
functions.  Our goal here is to consider functions of several quaternionic variables in
$\mathbb{H}^n$ and on manifolds with quaternionic structure and present some applications
in geometry.

For a point $q\in \Hn$ we shall write $q=(q^1, \dots, q^n)$ with $q^\alpha\in
\mathbb{H}$, $q^\alpha=t^\alpha+i x^\alpha+jy^\alpha+kz^\alpha$ for $\alpha = 1,\dots,
n$. Furthermore, $q^{\overline\alpha}=\overline{q^\alpha}$ , i.e., $q^{\overline\alpha} =
t^\alpha-i x^\alpha-jy^\alpha-kz^\alpha$.

We recall that a function $F\ :\ \Hn\rightarrow\ \mathbb{H}$, which is continuously
differentiable when regarded as a function of $\mathbb{R}^{4n}$ into $ \mathbb{R}^4$ is
called quaternionic regular, or just regular for short, if\\
\centerline{$
\dbara F\ =\ \dta {F}\ + \ i\dxa {F}\ + \ j\dya {F}\ +\ k\dza {F}\ =\
 0, \quad \alpha = 1,\dots, n.$}
In other words, a real-differentiable function of several quaternionic variables is
regular if it is regular in each of the variables (see \cite{Per1,Per2,Joy}.
The
condition that $F=f+iw+ju+kv$ is regular is equivalent to the following
Cauchy-Riemann-Feuter equations
\begin{equation}\label{reg}
\begin{aligned}
\dta f -\dxa w  -\dya u - \dza v \ =\ 0,\quad
\dta w + \dxa f + \dya v - \dza u \ =\ 0,\\
\dta u - \dxa v + \dya f + \dza w \ =\ 0, \quad
\dta v + \dxa u - \dya w + \dza f  \ = \ 0.
\end{aligned}
\end{equation}
\begin{dfn}
A real-differentiable function $f: \Hn\mapsto \mathbb{R}$ is called $\bar
Q$-\emph{pluriharmonic} if it is the real part of a regular function.
\end{dfn}
\begin{prop}\label{pluriharmonic}\sideremark{pluriharmonic}
Let $f$ be a real-differentiable function $f: \Hn\mapsto \mathbb{R}$. The following
conditions are equivalent
\begin{enumerate}[i)]
\item $f$ is $\bar Q    $-pluriharmonic, i.e., it is the real part of a regular
function;
\item  $\dbarb\dira f\ =\ 0$ for every $\alpha, \beta\in \{1,\dots, n\}$,   where
$\dira\ =\ \dta {}\ - \ i\dxa {}\ - \ j\dya {}\ -\ k\dza {}$;
\item $f$ satisfies the following system of PDEs
\begin{equation}\label{qreg}
\begin{aligned}
 \dtbta f +  \dxbxa f + \dybya f + \dzbza f \  & =\ 0, \quad
 \dxbta f - \dtbxa f - \dybza f + \dzbya f \ & =\ 0 \\
 -\dtbya f + \dxbza f +\dybta f - \dzbxa f \ & =\ 0, \quad
 -\dtbza f - \dxbya f + \dybxa f + \dzbta f\ & = \ 0.
\end{aligned}
\end{equation}
\end{enumerate}
\end{prop}
\begin{proof}
It is easy  to check that  $\dbarb\dira f =0$ is equivalent to
\eqref{qreg}.

We turn to the proof of ii) implies i). Let $f$ be real valued function on
$\Hn$, such that, $\dbarb \dira f\ = \ 0$. We shall construct a
real-differentiable regular function $F: \Hn \mapsto \mathbb{H}$. In fact,
for $q\in\Hn$ we define
\begin{equation*}
F(q)\ =\ f(q)\ +\ \Im\, \int_0 ^1 s^2\, (\dira f)(sq)\,  q^\alpha \, ds.
\end{equation*}
\noindent In order to rewrite the imaginary part in a different way we compute
\begin{multline*}
\Re \, \int_o^1 s^2\, (\dira f)(sq)\,  q_\alpha \, ds
=\Re \int_0^1 s^2 \Bigl ( \dta {f}- i\dxa {f} - j\dya {f} -  k\dza {f}\Bigr )\,
(sq) \,
(t_\alpha + i x_\alpha + j y_\alpha + k z_\alpha)\, ds\\
   =\ \int_0^1 s^2 \, \Bigl ( \dta f (sq)t_\alpha + \dxa f (sq)x_\alpha +
\dya f (sq)y_\alpha +\dza f (sq)z_\alpha \Bigr ) \,
  ds\\
   =\ \int_0^1 s^2 \ds  \Bigl ( f(sq)\Bigr)\, ds\ =\ s^2f(sq)|_0^1\
  -\ 2\int_0^1 sf(sq)\, ds
   =\ f(q)\ -\ 2 \int_0^1 sf(sq)\, ds.
\end{multline*}
\noindent Therefore we have
$\Im\, \int_o^1 s^2\, (\dira f)(sq)\, q^\alpha \, ds\ =\ \int_o^1 s^2\, (\dira f)(sq)\,
q^\alpha \, ds \ -\ f(q) \ +\ 2 \int_0^1 sf(sq)\, ds.
$

\noindent In turn, the formula for $F(q)$ becomes
$F(q)\ =\ \int_o^1 s^2\, (\dira f)(sq)\,  q^\alpha \, ds \
 +\ 2 \int_0^1 sf(sq)\, ds.
$ 
\noindent  Hence
$\dbarb F (q)\ =\ \int_o^1 s^2\, \dbarb \bigl [(\dira f)(sq)\, q^\alpha \bigr ] \, ds \
 +\ 2 \int_0^1 s\, \dbarb \bigl [ f(sq) \bigr ]\, ds.
$ 
We compute the first term,
\begin{multline*}
\dbarb \bigl[(\dira f)(sq)\, q^\alpha \bigr] = (\dtb + i\dxb+j\dyb+k\dzb)
\bigl[(\dira f)(sq)\, q^\alpha \bigr]\\
   =\ \dbarb \bigl [(\dira f)(sq)\, \bigr ] \,q^\alpha +\ \dira f (sq)\  \dtb {q_\alpha}+\ i \dira f (sq)\ \dxb {q_\alpha} \\
    +\ j\dira f (sq)\ \dyb {q_\alpha} \ + \ k\dira f (sq)\ \dzb {q_\alpha}\\
    =\ \dbarb \bigl [(\dira f)(sq)\, \bigr ] \,q^\alpha +\ \delta_{\alpha\beta}
    \{ \dira f (sq) \ +\ i \dira f (sq)i \ +\j\dira f (sq)j \ + \ k\dira f (sq)k  \}.
\end{multline*}
The  last term can be simplified, using the fundamental property that the coordinates of
a quaternion can be expressed by the quaternion only, as follows
\begin{multline*}
 \dirb f (sq) \ +\ i \dirb f(sq)i \ +\ j\dirb f (sq)j \ + \ k\dirb f (sq)k\\
   =\ \diracb f \ + \ i \diracb f i\\
   + \ j \diracb f j\ + \ k \diracb f k \\
   =\ -2\dtb f \ -\ i\dxb f\ -\ j\dyb f \ -\ k\dzb f\ +\ i\dxb f \ -\  j\dyb f\  -\  k\dzb f\\
     -\ i \dxb f \ + \ j \dyb f \ -\ k\dzb f\  -\  i \dxb f\ - \ j \dyb f\ + \ k \dzb f\\
       =\ -2\dtb f \ -\ 2i\dxb f\ -\  2j\dyb f\ -\ 2 k\dzb f\ =\ -2\dbarb f (sq).
\end{multline*}
Going back to the computation of $\dbarb F (q)$, we find
\begin{multline*}
\dbarb F (q)\ =\ \int_0^1 \dbarb \bigl [(\dira f)(sq)\, \bigr ]
\,q^\alpha \, ds\ -\ 2\int_0^1 s^2 \dbarb f (sq)\, ds\ +\ 2\int_0^1 s^2 \dbarb f (sq)\, ds\\
   =\ \int_0^1 \dbarb \bigl [(\dira f)(sq)\, \bigr ]
\,q^\alpha \, ds.
\end{multline*}
Hence, if $\dbarb \dira f\ =\ 0$ for every $\alpha$ and $ \beta$ we have $\dbarb F (q)\
=\ 0$.

Next we show that i) implies ii). Using \eqref{reg}, we have
\begin{multline*}
\dxbta f - \dtbxa f \  + \dzbya f \ - \dybza f= \dxbxa w + \dxbya u + \dxbza v + \dtbta w + \dtbya v -\dtbza u\\
-\dzbta u + \dzbxa v - \dzbza w + \dybta v + \dybxa u - \dybya w.
\end{multline*}
Both sides must be equal to zero by noticing that the left hand side is antisymmetric
while on the right we have an expression symmetric
 with respect to exchanging $\alpha$ with $\beta$.  The other
 identities can be obtained similarly.
\end{proof}
According to \cite{Sti} there are exactly two kinds of Cauchy-Riemann equations for
functions of several quaternionic variables. The second one turns out to be most suitable
for the geometric purposes considered in this paper.
\begin{dfn}\label{d:anti-regular functions}
A function $F\ :\ \Hn\rightarrow\ \mathbb{H}$, which is continuously differentiable when
regarded as a function of $\mathbb{R}^{4n}$ into $ \mathbb{R}^4$ is called
\emph{quaternionic anti-regular} ( also \emph{anti-regular}), if
\[
\mathcal D F\ =\ \dta {F}\ - \ i\dxa {F}\ - \ j\dya {F}\ -\ k\dza {F}\ =\
 0, \quad \alpha = 1,\dots, n.
\]
\end{dfn}
The condition that $F=f +iw+ju+kv$ is anti-regular function is equivalent to the
following Cauchy-Riemann-Feuter  equations
\begin{equation}\label{ahol}
\begin{aligned}
\dta f +\dxa w  +\dya u + \dza v \ =\ 0,\quad \dta w - \dxa f - \dya v + \dza u \ =\ 0,\\
\dta u + \dxa v - \dya f - \dza w \ =\ 0,\quad \dta v - \dxa u + \dya w - \dza
          f \ =\ 0.
\end{aligned}
\end{equation}
See also   \eqref{ahol1} for an equivalent form of the above system.

Anti-regular functions on hyperk\"ahler and quaternionic K\"ahler
manifolds are studied in \cite{CL1,CL2,LZ}, under the name quaternionic
maps, in connection with minimal surfaces and  maps between quaternionic
K\"ahler manifolds  preserving the sphere of almost complex structures.
Thus, the anti-regular functions considered here are quaternionic maps
between $\mathbb{H}^n$ and $\mathbb{H}$ with a suitable choice of the
coordinates.
\begin{dfn}
A real-differentiable function $f: \Hn\mapsto \mathbb{R}$ is called  \emph{quaternionic
pluriharmonic} (\emph{ Q-pluriharmonic} for short) if it is the real part of an
anti-regular function.
\end{dfn}
The anti-regular functions and their real part play a significant role in
the theory of hypercomplex manifold as well as in the theory of
quaternionic contact (hypercomplex contact) manifolds as we shall see
further in the paper. We need a real expression of the second order
differential operator $\dira\dbarb f$ acting on a real function $f$.

We  use the standard hypercomplex structure on $\Hn$ determined by the
action of the imaginary quaternions
\begin{align*}
I_1\, \diffta\ =\ \diffxa, & \quad I_1\, \diffya\  =\ \diffza, \quad\quad
I_2\diffta\ =\ \diffya, \quad I_2\,\diffxa\ =\ -\diffza\\
& I_3\diffta\ =\ \diffza, \quad\quad I_3\,\diffxa\ =\ \diffya.
\end{align*}
We recall a convention. For any p-form $\psi$ we consider the p-form $I_s\psi$ and three
(p+1)-forms $d_s\psi, \quad s=1,2,3$ defined by
$I_s\psi(X_1,\dots,X_p):=(-1)^p\psi(I_sX_1,\dots,I_sX_p),\qquad
d_s\psi:=(-1)^pI_sdI_s\psi.
$ 
Consider the  second order differential operators  $DD_{I_s}$ acting on the exterior
algebra  defined by \cite{HP}
\begin{equation}\label{DD1} DD_{I_i}:=
dd_i+d_jd_k=dd_i-I_jdd_i=dd_i-I_kdd_i.
\end{equation}
\begin{prop}\label{antiplu}
Let $f$ be a real-differentiable function $f: \Hn\rightarrow \mathbb{R}$. The following
conditions are equivalent
\begin{enumerate}[i)]
\item $f$ is {\bf Q}-pluriharmonic, i.e. it is the real part of an anti-regular function; \item
$DD_{I_s}\, f\ =\ 0,\quad s=1,2,3.$; \item  $\dira \dbarb f\ =\ 0$  for every $\alpha,
\beta\in \{1,\dots, n\}$, where
$\dira\ =\ \dta {}\ - \ i\dxa {}\ - \ j\dya {}\ -\ k\dza {}$,
\item  $f$ satisfies the following system of PDEs
\begin{align*}
 \dtbta f +  \dxbxa f + \dybya f + \dzbza f \   =\ 0,\quad
 -\dxbta f + \dtbxa f - \dybza f + \dzbya f \  =\ 0, \\
 \dtbya f + \dxbza f - \dybta f - \dzbxa f\  =\ 0, \quad
 \dtbza f - \dxbya f + \dybxa f - \dzbta f\  = \ 0.
\end{align*}
\end{enumerate}
\end{prop}
\begin{proof}
A simple calculation of $\dirb \dbara f$ gives the equivalence between
iii) and iv).

Next, we shall show that ii) is equivalent to iii). As $df = \dta f\diffta
+ \dxa f \diffxa  + \dya f \diffya + \dza f\diffza$ we have
$I_1 df =\dta f\diffxa - \dxa f\diffta + \dya f\diffza - \dza f \diffya.
$ 
A routine calculation gives the following formula
 \begin{multline}
  D\,D_{I_1}\, f\
  =\ \sum_{\alpha, \beta}\Re(\dirb\dbara f)\, [\diffta\wedge\diffxb +\diffya\wedge\diffzb ]\\
  +\sum_{\alpha < \beta}\Re(i\dirb\dbara f)\, [- \diffta\wedge\difftb
 - \diffxa\wedge\diffxb
 + \diffya\wedge\diffyb
 + \diffza\wedge\diffzb ]\\
 + \sum_{\alpha, \beta}\Re(j\dirb\dbara f)\,[\diffta\wedge\diffzb - \diffxa\wedge\diffyb]\ - \ \sum_{\alpha, \beta}\Re(k\dirb\dbara f)\,[ \diffta\wedge\diffyb + \diffxa\wedge\diffzb].
\end{multline}
Similar formulas hold for $DD_{I_2}$ and $DD_{I_3}$. Hence, the
equivalence of ii) and iii) follows.

The proof of the implication iii) implies i) is analogous to the proof of the
corresponding implication in Proposition~\ref{pluriharmonic}. Define
$F(q)\ =\ f(q)\ +\ \Im\, \int_o^1 s^2\, (\dbarb f)(sq)\, q^{\overline\beta} \, ds,
$ 
and a small calculation shows that this defines an anti-regular function, i.e., $\dira F\
=\ 0$ for every $\alpha$.

In order to see that iii) follows from i) we can proceed as in Proposition
\ref{pluriharmonic} and hence we skip the details. See also another proof in
Proposition~\ref{antipluman}
\end{proof}
\begin{rmrk}
We note that Proposition~\ref{pluriharmonic} and Proposition~\ref{antiplu} imply that the
real part of a regular function is not in the kernel of the operators $DD_{I_s}$ which is
one of the main difference between regular and anti-regular function.
\end{rmrk}

\subsection{Quaternionic pluriharmonic functions on hypercomplex manifold}

We recall that a hypercomplex manifold is a smooth $4n$-dimensional manifold $M$ together
with a triple $(I_1,I_2,I_3)$ of integrable almost complex structures satisfying the
quaternionic relations $I_1I_2=-I_2I_1=I_3$.  The second order differential operators
$DD_{I_i}$ defined in \cite{HP} by \eqref{ddi} having  the origin in the papers
\cite{Sal,Sal1,CSal} play an important r\'ole in the
 theory of quaternionic  plurisubharmonic functions (i.e. a real function for which
 $DD_{I_s}(.,I_s.)$ is positive definite) on hypercomplex manifold
 \cite{A1,A2,V,AV,A3} as well as the potential theory of
 HKT-manifolds. We recall that  Riemannian metric $g$ on a
 hypercomplex manifold compatible with the
three complex structures is said to be HKT-metric \cite{HP} if the three corresponding
K\"ahler forms $\Omega_s=g(I_s.,.)$ satisfy $d_1\Omega_1=d_2\Omega_2=d_3\Omega_3$. A
smooth real function is a  HKT-potential if locally it  generates the three K\"ahler
forms, $\Omega_s=DD_{I_s}f$ \cite{MS,GP}, in particular such a function is quaternionic
plurisubharmonic. The existence of a HKT potential on any HKT metric on $\Hn$ is proved
in \cite{MS} and for any HKT metric in \cite{BS}.

Regular functions on hypercomplex manifold are studied from analytical point
\cite{Per1,Per2}, from algebraic point \cite{Joy,Q}. However, as we have already
mentioned, regular functions are not the appropriate functions for our purposes mainly
because they have no direct connection with the second order differential operator
$DD_{I_s}$.

Here we consider anti-regular functions and their real parts on hypercomplex manifold.
\begin{dfn}
Let $(M,I_1,I_2,I_3)$ be a hypercomplex manifold. A quaternionic valued function $F:M
\longrightarrow f+iw+ju+kv\in \mathbb H$ is said to be  \emph{anti-regular} if  any one
of the following relations between the differentials of the coordinates  hold
\begin{equation}\label{ahol1}
\begin{aligned}
df=d_1w+d_2u+d_3v\qquad d_1f=-dw+d_3u-d_2v\\
d_2f=-d_3w-du+d_1v \qquad d_3f=d_2w-d_1u -dv.
\end{aligned}
\end{equation}
A real valued function $f:M\longrightarrow \mathbb R$ is said to be
\emph{quaternionic pluriharmonic} ( or \emph{Q-pluriharmonic}) if it is
the real part of an anti-regular function.
\end{dfn}

Observe that the system \eqref{ahol} is equivalent to \eqref{ahol1}. We have the
hypercomplex manifold analogue of Proposition~\ref{antiplu}

\begin{prop}\label{antipluman}
Let $(M,I_1,I_2,I_3)$ be a hypercomplex manifold and let $f$ be a real-differentiable
function on $M$, $f: M\longrightarrow \mathbb{R}$. The following conditions are
equivalent
\begin{enumerate}[i)]
\item $f$ is Q-pluriharmonic, i.e. it is the real part of an anti-regular function; \item
$DD_{I_s}\, f\ =\ 0,\quad s=1,2,3.$;
\end{enumerate}
\end{prop}

\begin{proof}
It is easy to verify that if each $I_s$ is integrable almost complex structure then we
have the identities \cite{HP}
\begin{equation}\label{comj} dd_s+d_sd=0, \qquad d_sd_r+d_rd_s=0,
  s,r=1,2,3.
\end{equation}
Using the commutation relations \eqref{comj}, we  get readily that i) implies ii). For
example, \eqref{ahol1} yields
\begin{gather*}
dd_1f+d_2d_3f +d^2w-d_2^2w-dd_3u+d_2d_1u+dd_2v+d_2dv\\
d_1df+d_3d_2f -d_1^2w+d_3^2w-d_1d_2u+d_3du-d_1d_3v-d_3d_1v.
\end{gather*}
Subtracting the two equations and using the commutation relations \eqref{comj} we get
$DD_{I_1}f=0$.

For the converse, observe that $DD_{I_1}f=0 \Leftrightarrow dd_1f=I_2dd_1f$. The
$\partial\bar\partial$-lemma for $I_2$ gives the existense of a smooth function $A_1$
such that $dd_1f=dd_2A_1$. Similarly, using the Poincare lemma, we obtain
$d_1f-d_2A_1-dB_1=0, \quad d_2f-d_3A_2-dB_2=0, \quad
d_3f-d_1A_3-dB_3=0$ for a smooth functions $A_1,A_2,A_3,B_1,B_2, B_3$.
The latter implies
$df+d_1(A_2+B_1)+d_2(B_2-A_3)+d_3(A_1-B_3)=0.$
Set $w=-A_2-B_1, u=A_3-B_2, v=B_3-A_1$ to get the equivalence between i) and ii).
\end{proof}

\subsubsection{Restriction on hyper-surfaces}\label{6.2.1} In this section  we
shall denote with $\langle .,. \rangle$ the Euclidean scalar product in
$\mathbb{R}^{4n+4} \cong \Hnn$ and with $\tilde I_j$, $j=1,2,3$, the standard almost
complex structures on $\Hnn$. Let $M$ be a smooth hyper-surface in $\Hnn$ with a defining
function $\rho$, $M\ =\ \{ \rho=0 \}$, $d \rho\not= 0$, and $i : M\hookrightarrow \Hnn$
be the embedding. It is not hard to see that at every point $p\in M$ the subspace
$H_p=\bigcap_{j=1}^3 \, \tilde I_j\, (T_pM)$ of the tangent space $T_pM$ of $M$ at $p$ is
the largest subspace invariant under the almost complex structures and dim$H_p = 4n$. We
shall call $H_p$ the horizontal space at $p$. Thus on the horizontal space $H$ the almost
complex structures $I_j$, $j=1,2,3$, are the restrictions of the standard almost complex
structures on $\Hnn$. In particular, for a horizontal vector $X$ we have
\begin{equation}\label{acsrestr}
\tilde I_j\, i_{*}X\ =\ i_{*}\,(I_jX).
\end{equation}
  Let $\tilde\theta^j=\tilde I_j\ \frac{d\rho}{|d\rho|}$. We drop
the tilda in the notation of the almost complex structures when there is no ambiguity.

We define three one-forms on $M$ by setting $\theta^j=i^*\tilde\theta^j=i^*(\tilde I_j\
\frac{d\rho}{|d\rho|})$, i.e.,
$\theta_j\ (\ .\ )\ =\ - \frac {d\rho\, (\tilde I_j\, .\ )}{|d\rho|}\ =\ \langle \ .\, ,
\tilde I_j N \rangle,
$ 
where $N=\frac{D\rho}{|D\rho|}$ is the unit normal vector to $M$. 
We  describe the hypersurfaces which inherit a natural quaternionic contact structure
from the standard structures on $\Hnn$ (see also \cite{D1}) in the next
\begin{prop}\label{p:QRhypersurface}
If $M$ is a smooth hypersurface of $\Hnn$ then we have
\begin{equation}\label{e:horizontal2forms}
d\theta_1(I_1X,Y)=d\theta_2(I_2X,Y)=d\theta_3(I_3X,Y) \quad (X,\, Y\in H)
\end{equation}
if and only if the restriction of the second fundamental form of $M$ to the horizontal
space is invariant with respect to the almost complex structures, i.e. if $X$ and $Y$ are
two horizontal vectors we have  $II(I_j X, I_j Y)=II(X, Y)$. Furthermore, if the
restriction of the second fundamental form of $M$ to the horizontal space is positive
definite,  $II(X,X)> 0$  for any non-zero horizontal vector $X$, then $(M,\theta,
I_1,I_2)$ is a quaternionic contact manifold.
\end{prop}
\begin{proof}
Let $D$ be the Levi-Civita connection on $\mathbb{R}^{4n+4}$
and $X,\ Y$ be two horizontal vectors. As the
horizontal space is the intersection of the kernels of the one forms $\theta_j$ we have
\begin{multline}
d\theta_1(I_1X, Y)\ =\ -\theta_1([\tilde I_1X,Y])\ =\ -\langle [\tilde I_1X,Y],
\tilde I_1N\rangle\\
=\ -\langle D_{\tilde I_1X}Y - D_{Y}(\tilde I_1X), \tilde I_1N\rangle\ =\ -\langle
D_{\tilde I_1X}Y , \tilde I_1N \rangle + \langle D_{Y}(\tilde I_1X), \tilde I_1N\rangle\\
=\ \langle D_{\tilde I_1X}(\tilde I_1Y) , N \rangle + \langle D_{Y}X, N\rangle\
 =\ II(\tilde I_1X, \tilde I_1Y)+II(X,Y).
\end{multline}
Therefore $d\theta_1(I_1X,Y)=d\theta_2(I_2X,Y)$ iff $II(I_j X, I_j Y)=II(X, Y)$.

The last claim of the proposition is clear from the above formula.
In particular, $g_H(X, Y)=II(X,Y)$ is a metric on the horizontal
space when the second fundamental form is positive definite on the
horizontal space and we have $ d\theta_1(I_1X, Y)=2g_H(X, Y). $
Hence, $(M, \theta, I, J)$ becomes a quaternionic contact structure.
We denote the corresponding horizontal forms with $\omega_j$, i.e.,
$\omega_j(X,Y)=g_H(I_jX, Y)$.
\end{proof}
\noindent Let us note also that in the situation as above
$g=g_H+\Sigma_{j=1}^3\theta_j\otimes\theta_j$ is a Riemannian metric
on $M$. In view of the above observations we define a
QC-hypersurface of $\mathbb H^{n+1}$ as follows.
\begin{dfn}\label{d:QRhypersurface}
We say that a smooth embedded hypersurface of $\Hnn$ is a \emph{$QC$-hypersurface} if the
restriction of the second fundamental form of $M$ to the horizontal space is a definite
symmetric form, which is invariant with respect to the almost complex structures.
\end{dfn}
Clearly every sphere in $\Hnn$ is a QC-hypersurface and this is true also
for the ellipsoids $\sum_a\frac {|q^a|^2}{b_a}=1$. In fact, a hypersurface
of $\Hnn$ is a QC-hypersurface if and only if the (Euclidean) Hessian of
the defining function $\rho$, considered as a quadratic form on the
horizontal space, is a symmetric definite matrix from $GL(n, \mathbb{H})$,
the latter being the linear group of invertible matrices which commute
with the standard complex structures on $\mathbb{H}^n$. The same statement
holds for hypersurfaces in quaternionic K\"ahler and hyperk\"ahler
manifolds.

\begin{prop}\label{Qplurharm}
Let $i:M\rightarrow \Hn$ be a QC hypersurface in $\Hn$, $f$ a real-valued function on
$M$. If $f=i^*F$ is the restriction to $M$ of a  Q-pluriharmonic function $F$ defined on
$\Hn$, i.e. $F$  is the real part of an anti-regular function $F+iW+jU+kV$, then:
\begin{gather}\label{resant}
df=d(i^*F)=d_1(i^*W)+d_2(i^*U)+d(i^*V) \qquad \text{mod} \ \eta,\\\label{restplu}
DD_{I_1} f (X, I_1Y)\ =\ -4dF(D \rho) \, g_H( X, Y )\ - \ 4 (\xi_2 f)\, \omega_2(X, Y)
\end{gather}
for any horizontal vector fields $X, Y\in H$.
\end{prop}
\begin{proof}
Let us prove first \eqref{resant}. Denote with small letters the restrictions of the
functions defined on $\mathbb{H}^n$. For $X\in H$ from \eqref{acsrestr} we have
\begin{multline*}
(i^* \, \tilde I_1\, dW)(X)= (\tilde I_1 dW)(i_{*}\, X)= -dW
(\tilde I_1 \,i_{*}\, X)
 =-dW (i_{*}\,  (I_1 \ X))= -dw(I_1\, X)=d_1 w\, (X).
\end{multline*}
Applying the same argument to the functions $U$ and $V$ we see the validity of
\eqref{resant}.

Our goal is to write the equation for $f$
on $M$, using the fact that $f=i^*F$. Let us consider the function $\lambda$,
$
\lambda\ =\ \frac {dF(D\rho)}{|D\rho|^2},
$
and the one-form $d_M F$,
$
d_M F\ =\ dF\ -\ \lambda\, d\rho.
$
Thus the one-form $df$ satisfies the equation
$df\ =\ i^* (d_M F\ +\ \lambda\,d\rho)\ =\ i^*(d_M F),
$ 
taking into account that $(\lambda\circ i )\,d (\rho\circ i)\ =\ 0$ as $\rho$ is constant
on $M$. From Proposition \ref{antiplu}, the assumption on $F$ is equivalent to
$DD_{\tilde I_j}F=0$. Therefore, we have
\begin{align*}
0\ =\ DD_{{\tilde I_1}}F\  =\ d{\tilde I_1} dF\ -\ {\tilde I_2} d{\tilde I_1} dF
                     =\ d\,({\tilde I_1} d_MF\ +\ \lambda {\tilde I_1} d\rho)\
                                -\ {\tilde I_2} d\, ({\tilde I_1} d_MF\ +\ \lambda {\tilde I_1} d\rho)\\
                     =\ d\,{\tilde I_1} d_MF\ +\ d\lambda \wedge{\tilde I_1}
                              d\rho\ +\ \lambda d{\tilde I_1} d\rho
                             -\  {\tilde I_2} d\,{\tilde I_1} d_MF\ - \ {\tilde I_2} (d\lambda \wedge{\tilde I_1}
                              d\rho)\ -\ \lambda {\tilde I_2}{\tilde I_1}
                              d\rho.
\end{align*}
Restricting to $M$, and in fact, to the horizontal space $H$ we find
\begin{align}\label{restrictioneqn1}
0\  =\ i^* \, (DD_{{\tilde I_1}}F)\vert_H
  =\ i^*d\ ({\tilde I_1} d_M
F)\vert_H \ +\ d(\lambda\circ i) \wedge i^* ({\tilde I_1}
                              d\rho)\vert_H \ +\ (\lambda\circ i ) d i^* ({\tilde I_1} d\rho)\vert_H\\
                                     \quad \quad -\   i^* ({\tilde I_2} d\,{\tilde I_1} d_MF)\vert_H \ -
                                    \ i^* ({\tilde I_2} (d\lambda \wedge{\tilde I_1}
                                        d\rho)) \vert_H \ -\ (\lambda\circ i)\ i^*( {\tilde I_2}d\,{\tilde I_1}d\rho)
                                        \vert_H.\notag
\end{align}
Since the horizontal space is in the kernel of the one-forms $\theta_j\ =\ \tilde I_j\,
d\rho\vert_H$ it follows that
$i^*\, (\tilde I_j\, d\rho )\vert_H\ =\ 0.
$ 
Hence, two of the terms in \eqref{restrictioneqn1} are equal to zero, and we have
\begin{align}\label{restrictioneqn2}
0\ =\ i^* \, (DD_{{\tilde I_1}}F)_{\vert_H} \  =\ i^* ( d\, {\tilde I_1} d_M
F \ -\  {\tilde I_2} d\,{\tilde I_1} d_MF)\ _{\vert_H}
                                     +\  (\lambda\circ i )\ i^*  \bigl ( d\, {\tilde I_1} d\rho
                                            \ -\ {\tilde I_2}d\,{\tilde I_1}d\rho \bigr )
                                            _{\vert_H}. \notag
\end{align}
In other words for horizontal $X$ and $Y$ we have
\begin{equation}\label{e:restriction}
i^* ( d\, {\tilde I_1} d_M F \ -\  {\tilde I_2} d\,{\tilde I_1} d_MF)(X, IY)\ =\
-(\lambda\circ i )\ i^*  \bigl ( d\, {\tilde I_1} d\rho
                                            \ -\ {\tilde I_2}d\,{\tilde I_1}d\rho \bigr
                                            )(X, IY)
\end{equation}
\noindent The right-hand side is proportional to the metric. Indeed,
recalling
\begin{equation*}
i^*\ ( \tilde I_j\,d\, \rho) (X)\ =\ |d \rho|\ \theta_j(X)  d\, \theta_j
(X, Y)\ =\ 2g(I_j X, Y)\end{equation*}  we obtain the identity
\begin{multline*}
i^*  \bigl ( d\, {\tilde I_1} d\rho \ -\ {\tilde I_2}d\,{\tilde I_1}d\rho \bigr)
(X, Y)\ =\ 2|d \rho|\ g(I_1X, Y)\ -\ 2 |d \rho|\ g(I_1I_2X, I_2Y)\\
         =\ 2|d \rho|\ g(I_1X, Y)\ -\ 2g(I_3X, I_2Y)\ =\ 4|d \rho|\ g(I_1X, Y).
        \end{multline*}
Let us consider now the term in the left-hand side of \eqref{e:restriction}. Decomposing
$d_MF$ into horizontal and vertical parts we write
$ d_M
F\ =\ d_Hf\ +\ F_j\tilde\theta^j.
$
>From the definitions of the forms $\tilde\theta^j$ we have
$
\tilde I_1 \tilde\theta^1 = \frac {d \rho}{|d \rho|}, \quad \tilde I_1 \tilde\theta^2
=\tilde\theta^3 ,\quad \tilde I_1 \tilde\theta^1 \tilde\theta^3= -\tilde\theta^2.
$
Therefore
\begin{multline*}
d{\tilde I_1} d_M F\ =\ d{\tilde I_1} d_H F\ +\ dF_j\wedge{\tilde I_1}\theta^j\ +\ F_1 d(
\frac {d \rho}{|d \rho|})\ +\ F_2d\tilde\theta^3 \ -\
F_3d\tilde\theta^2\\
=\ d{\tilde I_1} d_H F\ +\ dF_j\wedge{\tilde I_1}\theta^j\ -\ |d\rho|^{-2}d |d\rho|\wedge
d\rho\ + \  F_2d\tilde\theta^3 \ -\ F_3d\tilde\theta^2\\
{\tilde I_2} d{\tilde I_1} d_M F\ =\ {\tilde I_2}  d{\tilde I_1} d_H F\ +\ {\tilde I_2} dF_j\wedge{\tilde I_2} {\tilde I_1}\theta^j\\
-\ |d\rho|^{-2}{\tilde I_2} d|d\rho|\wedge {\tilde I_2} d\rho\ + \ F_2{\tilde I_2}
d\tilde\theta^3 \ -\ F_3 {\tilde I_2}  d\tilde\theta^2.
\end{multline*}
\noindent From $I_2 d\theta^3\ =\ -d\theta^3$, $I_2 d\theta^2\ =\ d\theta^2$ and the
above it follows
\begin{multline*}
i^* ( d\, {\tilde I_1} d_M F -  {\tilde I_2} d\,{\tilde I_1} d_MF)_{|_H} = DD_{I_1}f
+F_2 d\theta^3-  F_3 d\theta^2 +F_2 d\theta^3
+ F_3 d\theta^
=DD_{I_1}f +4F_2\ \omega_3.
\end{multline*}
\noindent In conclusion, we proved
$
DD_{I_1} f (X, Y)\ =\ -4(\lambda\circ i)\ |\nabla\rho|\  g(I_1X, Y)\ -\ 4F_2\, \omega_3
(X, Y)
$
from where the claim of the Proposition.
\end{proof}

\subsection{Anti-CRF functions on Quaternionic contact manifold}

Let $(M,\eta,\mathbb Q)$ be a (4n+3)-dimensional quaternionic
contact manifold and $\nabla$ denote the Biquard connection on $M$.
The equation \eqref{resant} suggests the following
\begin{dfn}\label{d:anti-CRF functions}
A smooth $\mathbb H$-valued function $F:M\longrightarrow \mathbb H, \quad F=f+iw+ju+kv,$
is said to be an \emph{anti-CRF function} if the smooth real valued functions $f,w,u,v$
satisfy
\begin{equation}\label{crfhyp}
df\ =\ d_1 w\ +\ d_2 u\ +\ d_3 v \ \text{  mod  }\ \eta,
\end{equation}
where $d_i\ =\ I_i\circ d.$
\end{dfn}
Choosing a local frame $\{T_a,X_a=I_1T_a,Y_a=I_2T_a,Z_a=I_3T_a,\xi_1,\xi_2,\xi_3\}, \quad
a=1,\dots,n$ it is easy to check that a $\mathbb H$-valued function $F=f+iw+ju+kv$ is an
anti-CRF function if it belongs  to the kernel of the operators
\begin{equation}\label{crf1}
D_{T_{\alpha}}\ =\ T_{\alpha}-iX_{\alpha}-jY_{\alpha}-kZ_{\alpha},
\qquad D_{T_{\alpha}}F\ =\ 0, \qquad \alpha =1,\dots n.
\end{equation}
\begin{rmrk}We note that anti-CRF functions have  different properties than the
CRF functions \cite{Per1,Per2}which are defined to be in the kernel of the operator
\begin{equation*}
\overline{D}_{T_{\alpha}}\ =\ T_{\alpha}+iX_{\alpha}+jY_{\alpha}+kZ_{\alpha}, \qquad
\overline D_{T_{\alpha}}F\ =\ 0, \qquad \alpha =1,\dots n.
\end{equation*}
\end{rmrk}
Equation \eqref{crfhyp} and a small calculation give the following Proposition.
\begin{prop}\label{crf}
A $\mathbb H$-valued function $F=f+iw+ju+kv$ is an anti-CRF function if and only if the
smooth functions $f,w,u,v$ satisfy the horizontal Cauchy-Riemann-Fueter equations
\begin{equation}\label{crf2}
\sideremark{crf2}
\begin{array}{c}
T_{\alpha}f=-X_{\alpha}w-Y_{\alpha}u-Z_{\alpha}v, \quad
X_{\alpha}f=T_{\alpha}w+Z_{\alpha}u-Y_{\alpha}v,\\
Y_{\alpha}f=-Z_{\alpha}w+T_{\alpha}u+X_{\alpha}v,\quad
Z_{\alpha}f=Y_{\alpha}w-X_{\alpha}u+T_{\alpha}v.
\end{array}
\end{equation}
\end{prop}
Having the quaternionic contact form $\eta$ fixed, we may extend the definitions
\eqref{DD1} of $DD_{I_i}$ to the second order differential operator $DD_{I_i}$ acting on
the real-differentiable functions $f:M\rightarrow \mathbb R$ by
\begin{equation}\label{Ivan8} DD_{I_i}f:=dd_if+d_jd_kf=
dd_if-I_jdd_if=d(I_idf)-I_j(d(I_idf)).
\end{equation}
The following proposition provides some formulas, which shall be used later.
\begin{prop}\label{comddi1}
On a QC-manifold, for $X,Y\in H$, we have the following commutation relations
\begin{gather}
DD_{I_i}f(X,I_iY)-DD_{I_k}f(X,I_kY)=-I_iN_{I_j}(X,I_iY)(f)-N_{I_k}(I_jX,I_jY)(f),\nonumber\\\label{comddi2}
d_id_jf(X,Y)+d_jd_if(X,Y)=-N_{I_k}(I_jX,I_iY)(f),\\\nonumber
dd_if(X,Y)+d_idf(X,Y)=N_{I_i}(I_iX,Y)f,\\\label{ddi2}
d_i^2f(X,Y)=-2\xi_i(f)\omega_i(X,Y)+2\xi_j(f)\omega_j(X,Y)+2\xi_k(f)\omega_k(X,Y),
\end{gather}
\noindent  In particular, on a hyperhermitian contact manifold we have
\begin{gather}
DD_{I_i}f(X,I_iY)-DD_{I_k}(X,I_kY)=4\xi_i(f)\omega_i(X,Y)-4\xi_j(f)\omega_j(X,Y),\nonumber\\\label{comddi3}
d_id_jf(X,Y)+d_jd_if(X,Y)=-4(\xi_i(f)\omega_j+\xi_j(f)\omega_i), \\\nonumber
dd_if+d_idf=4(\xi_k(f)\omega_j-\xi_j(f)\omega_k).
\end{gather}
\end{prop}

\begin{proof}
By the  definition \eqref{Ivan8} we obtain the second  and the third formulas in
\eqref{comddi2} as well as $DD_i(X,Y)+(dd_k-d_jd_i)(X,I_jY)=-I_iN_{I_j}(X,Y)$. The first
equality in \eqref{comddi2} is a consequence of the latter and the second equality in
\eqref{comddi2}. We have
$$d_i^2f(X,Y)=-I_id(I_i^2df)(X,Y)=d(df-\sum_{s=1}^3\xi_s(f)\eta_s)(I_iX,I_iY)$$
which is exactly  \eqref{ddi2}.
If $H$ is formally integrable then the formula \eqref{ntor} reduces to
$N_i(X,Y)=T_i^{0,2}(X,Y)$. The equation \eqref{comddi3} is an easy consequences of the
latter equality, \eqref{comddi2} and  \eqref{torha}
\end{proof}

Let us make the conformal change $\bar\eta=\frac{1}{2h}\eta$. The endomorphisms $\bar
I_i$ will coincide with $I_i$ on the horizontal distribution $H$ but they will have a
different kernel - the new vertical space $span\{\bar\xi_1,\bar\xi_2,\bar\xi_3\}$, where
$\bar\xi_s=2h\xi_s+I_s(\nabla h)$ (see \eqref{New19}). Hence,
for any $P\in \Gamma(TM)$ we have
\begin{multline}\label{Ivan5}
\bar I_i(P)=\bar I_i(P - \sum_{s=1}^3\tilde\eta_s(P) \bar\xi_s)=I_i(P
- \frac{1}{2h}\sum_{s=1}^3\eta_s(P)(2h\xi_s+I_s(\nabla h))) \\
=I_i(P)+\frac{1}{2h}\{\eta_i(P)\nabla h-\eta_j(P)I_k\nabla h+\eta_k(P)I_j\nabla h\}.
\end{multline}
\begin{prop}\label{p:DDi conf inv}
Suppose $\bar\eta=\frac{1}{2h}\eta$ are two conformal to each other structures.

a) The second order differential operator $DD_{I_i}$ (restricted on functions) transforms
as follows:
\begin{equation*} DD_{\bar
I_i}f-DD_{I_i}f=-2h^{-1}df(\nabla h)\omega_i-2h^{-1}df(I_j\nabla h)\omega_k\qquad
\text{mod} \ \eta.
\end{equation*}
b) If $f$ is the real part of the anti-CRF function $f+iw+ju+kv$ then the two forms
$$\Omega_i\ =\ DD_{I_i}\, f \ - \ \lambda \, \omega_i\ +\ 4(\xi_j f)\, \omega_k \qquad
\text{mod} \ \eta$$  are conformally invariant, where $\lambda\ =\ 4\big ( \xi_1 w+\xi_2u
+ \xi_3v \big ).$
\end{prop}

\begin{proof}
a) For any $X,Y\in H$, we compute
\begin{gather*}
d(\bar I_idf)(X,Y)=X(I_idf(Y))-Y(I_idf(X))-\bar I_idf[X,Y]
=d(I_idf)(X,Y)+df(\bar I_i[X,Y]-I_i[X,Y])
\end{gather*}
Here, we apply \eqref{Ivan5} to get
\begin{multline*}
d(\bar I_idf)(X,Y)= d(I_idf)(X,Y) +
\frac{1}{h}\{-df(\nabla h)\omega_i(X,Y) +  df(I_k\nabla h)\omega_j(X,Y) -
df(I_j\nabla h)\omega_k(X,Y)\}.
\end{multline*}
Now, apply the defining equation \eqref{Ivan8} to get the proof of a).

b) Assuming that $f$ is the real part of an anti-CRF function, from part a) we have 
\begin{align*}
\bar\Omega_i\ -\ \Omega_i\ & =\  D D_{\bar I_i}\, f \ - \ D D_{I_i}\, f \ - \
\bar\lambda\, \omega_i\ + \  \lambda\, \omega_i\ +\ 4\bar\xi_j f\, \omega_k \ -\ 4\xi_j
f\, \omega_k  \qquad
\text{mod} \ \eta \\
& = -4 (\xi_1 w  +\xi_2 u + \xi_3 v) \omega_i- 4\big ( (I_1
dh) w  +(I_2 dh)u+(I_3 dh) v \big ) \frac {\omega_i}{2h}\\
 & -\frac {2}{h}g(df, dh) \omega_i - \frac {2}{h} g(df,
d_j h) \omega_k  + 4(\xi_1 w+ \xi_2 u  + \xi_3 v) \omega_i + 4g(df, d_j h)
\frac {\omega_k}{2h}\ =\ 0 \quad \text{mod} \ \eta,
\end{align*}
taking into account \eqref{crfhyp}.
\end{proof}

We restrict our considerations to hyperhermitian contact manifolds.

\begin{thrm}\label{crfth1}
 If $f:M\rightarrow \mathbb R$ is the real part of an anti-CRF function $f+iw+ju+kv$
 on a (4n+3)-dimensional $(n>1)$ hyperhermitian contact
manifold $(M,\eta,Q)$.

Then the following equivalent conditions hold true:
\begin{enumerate}
\item[i)] 
The next equalities hold
\begin{equation}\label{crf03}
DD_{I_i}f=\lambda\omega_i- 4\xi_j(f)\omega_k \qquad \text{mod} \ \eta.
\end{equation}
\item[ii)] For
any $X,Y\in H$ we have the equality
\begin{align}\notag
(\nabla_Xdf)(Y) & +(\nabla_{I_1X}df)(I_1Y)+
(\nabla_{I_2X}df)(I_2Y)+(\nabla_{I_3X}df)(I_3Y)\\ \notag & =\
4\lambda g(X,Y)+
df(X)\alpha_3(I_3Y)+df(I_1X)\alpha_3(I_2Y)-df(I_2X)\alpha_3(I_1Y)-df(I_3X)\alpha_3(Y)\\\notag
& \quad
+df(Y)\alpha_3(I_3X)+df(I_1Y)\alpha_3(I_2X)-df(I_2Y)\alpha_3(I_1X)-df(I_3Y)\alpha_3(X).
\end{align}
\item[iii)] The function $f$ satisfies the  second order differential equations
\begin{align}\notag
& \Re(D_{T_{\beta}}\overline{D}_{T_{\alpha}}f) =\lambda
g(T_{\beta},T_{\alpha})\\\label{dd12} & \quad +
df(\nabla_{T_{\beta}}T_{\alpha})+df(\nabla_{I_1T_{\beta}}I_1T_{\alpha})+df(\nabla_{I_2T_{\beta}}I_2T_{\alpha})+
df(\nabla_{I_3T_{\beta}}I_3T_{\alpha})\\ \notag & \quad
+df(T_{\beta})\alpha_3(I_3T_{\alpha})+df(I_1T_{\beta})\alpha_3(I_2T_{\alpha})-
df(I_2T_{\beta})\alpha_3(I_1T_{\alpha})- df(I_3T_{\beta})\alpha_3(T_{\alpha})\\\nonumber
& \quad +df(T_{\alpha})\alpha_3(I_3T_{\beta})+df(I_1T_{\alpha})\alpha_3(I_2T_{\beta})-
df(I_2T_{\alpha})\alpha_3I_1(T_{\beta})- df(I_3T_{\alpha})\alpha_3(T_{\beta})\\\notag &
\\ \label{dd123} &
\Re(iD_{T_{\beta}}\overline{D}_{T_{\alpha}}f)=\Re(D_{I_1T_{\beta}}\overline{D}_{T_{\alpha}}f),
\quad
\Re(jD_{T_{\beta}}\overline{D}_{T_{\alpha}}f)=\Re(D_{I_2T_{\beta}}\overline{D}_{T_{\alpha}}f),\\\nonumber
& \hskip1truein
\Re(jD_{T_{\beta}}\overline{D}_{T_{\alpha}}f)=\Re(D_{I_3T_{\beta}}\overline{D}_{T_{\alpha}}f).
\end{align}
\end{enumerate}
The function $\lambda$ is determined by
\begin{equation}\label{lamb1}\lambda\ = \ 4(\xi_1(w)+\xi_2(u)+\xi_3(v)).\end{equation}
\end{thrm}
\begin{proof}
The proof includes a number of steps and occupies the rest of the section.

\emph{ i) } Suppose that there exists a smooth functions $w,u,v$ such that $F=f+iw+ju+kv$
is an anti-CRF function. The defining equation \eqref{crfhyp} yields
\begin{equation}\label{crfhyp1}
df=d_1w+d_2u+d_3v+\sum_{s=1}^3\xi_s(f)\eta_s,
\end{equation}
Since $d_s\eta_t(X,Y)=0$, for $s,t\in \{1,2,3\}, X,Y\in H$, applying \eqref{ddi2}
and \eqref{con1}, we obtain from \eqref{crfhyp1}
\begin{equation*}
\begin{aligned}
(dd_1f-dd_3u+dd_2v-2\xi_1(w)\omega_1-2\xi_2(w)\omega_2-2\xi_3(w)\omega_3)(X,Y)=0,\\
(d_1df-d_1d_2u-d_1d_3v+2\xi_1(w)\omega_1-2\xi_2(w)\omega_2-2\xi_3(w)\omega_3)(X,Y)=0,\\
(d_2d_3f+d_2d_1u+d_2dv-2\xi_1(w)\omega_1+2\xi_2(w)\omega_2-2\xi_3(w)\omega_3)(X,Y)=0,\\
(d_3d_2f+d_3du-d_3d_1v+2\xi_1(w)\omega_1+2\xi_2(w)\omega_2-2\xi_3(w)\omega_3)(X,Y)=0.
\end{aligned}
\end{equation*}
Summing  the first and the third equations, subtracting the second and the fourth and
using the commutation relations \eqref{comddi3} we obtain \eqref{crf03} with the
condition \eqref{lamb1}  which proves $i)$.

$iii)$ Equations \eqref{crf1} and \eqref{crf2} yield
\begin{align}\notag
2\Re(D_{T_{\beta}}& \overline{D}_{T_{\alpha}}f)\ = \
2(T_{\beta}T_{\alpha}f+X_{\beta}X_{\alpha}f+Y_{\beta}Y_{\alpha}f +Z_{\beta}Z_{\alpha}f)\\
\notag & =(\Re(D_{T_{\beta}}\overline{D}_{T_{\alpha}}f)+
\Re(D_{T_{\alpha}}\overline{D}_{T_{\beta}}f))+(\Re(D_{T_{\beta}}\overline{D}_{T_{\alpha}}f)-
\Re(D_{T_{\alpha}}\overline{D}_{T_{\beta}}f))=\\ \notag
& ([T_{\beta},T_{\alpha}]+[X_{\beta},X_{\alpha}]+[Y_{\beta},Y_{\alpha}]+[Z_{\beta},Z_{\alpha}])f
-([T_{\beta},X_{\alpha}]-[X_{\beta},T_{\alpha}]+[Y_{\beta},Z_{\alpha}]-[Z_{\beta},Y_{\alpha}])w\\
\notag & -
([T_{\beta},Y_{\alpha}]-[X_{\beta},Z_{\alpha}]-[Y_{\beta},T_{\alpha}]+[Z_{\beta},X_{\alpha})u-
([T_{\beta},Z_{\alpha}]+[X_{\beta},Y_{\alpha}]-[Y_{\beta},X_{\alpha}]-[Z_{\beta},T_{\alpha}])v.
\end{align}
Expanding the commutators and applying \eqref{crfhyp}, \eqref{torha}, \eqref{der} and
\eqref{intqc} gives \eqref{dd12}. Similarly, one can check the validity of \eqref{dd123}

$i)\Leftrightarrow ii)\Leftrightarrow iii)$ The next lemma establishes the equivalence
between i), ii) and iii).

\begin{lemma}\label{ddi}\sideremark{ddi} For any $X,Y\in H$ on a
  quaternionic contact manifold we have the identity
\begin{multline*}
DD_{I_1} f (X, I_1Y)= (\nabla_X df)Y  +  (\nabla_{I_1X} df)I_1Y +  (\nabla_{I_2X}
df)I_2Y  + (\nabla_ {I_3X}df)I_3X-4\xi_2(f)\omega_2(X,Y)\\
- df(X)\alpha_3(I_3Y)+df(I_1X)\alpha_2(I_3Y)+df(I_2X)\alpha_3(I_1Y)-df(I_3X)\alpha_2(I_1Y)\\
-df(Y)\alpha_2(I_2X)-df(I_1Y)\alpha_3(I_2X)+df(I_2Y)\alpha_2(X)+df(I_3Y)\alpha_3(X).
\end{multline*}
\end{lemma}
\begin{proof}[Proof of Lemma~\ref{ddi}]
Using the definition and also   \eqref{der}, \eqref{torha} and \eqref{symdh} we derive the
next sequence of equalities
\begin{multline*}
(dd_{I_1} f)(X, Y)
 = -(\nabla_X df)(I_1Y) + (\nabla_Y df)(I_1X)
- df(\nabla_X (I_1Y)  -\nabla_Y (I_1X) - I_1[X, Y])\\
= -(\nabla_X df)(I_1Y) + (\nabla_Y df)(I_1X)
+\alpha_2(X)df(I_3Y)-\alpha_3(X)df(I_2Y)-\alpha_2(Y)df(I_3X)+\alpha_3(Y)df(I_2X)\\
 = \ -(\nabla_Xdf)I_1Y\  +\  (\nabla_{I_1X} df)Y\ -
\ df (T(Y, \,I_1X))\\
+\alpha_2(X)df(I_3Y)-\alpha_3(X)df(I_2Y)-\alpha_2(Y)df(I_3X)+\alpha_3(Y)df(I_2X).
\end{multline*}
\begin{multline}\label{ddd1}
DD_{I_1}f(X,I_1Y)=(dd_{I_1}-I_2dd_{I_1})f \, (X, I_1Y)\\ = (\nabla_X df)Y  +
(\nabla_{I_1X} df)I_1Y  +  (\nabla_{I_2X} df)I_2Y + (\nabla_ {I_3X}df)I_3X -  df
(T(I_1Y,\, I_1X)) - df(T(I_3Y,\, I_3X))
\\-
df(X)\alpha_3(I_3Y)+df(I_1X)\alpha_2(I_3Y)+df(I_2X)\alpha_3(I_1Y)-df(I_3X)\alpha_2(I_1Y)\\
-df(Y)\alpha_2(I_2X)-df(I_1Y)\alpha_3(I_2X)+df(I_2Y)\alpha_2(X)+df(I_3Y)\alpha_3(X).
\end{multline}

\noindent A short calculation using \eqref{torha} gives
$ \ df (T(I_1Y, \, I_1X))\
+\ df(T(I_3Y, \, I_3X))\ =4\xi_2(f)\omega_2(X,Y).
 $
Inserting the last equality in \eqref{ddd1} completes the proof of
Lemma~\ref{ddi}.
\end{proof}
Since the structure is hyperhermitian contact, with the help of
\eqref{intqc} of Lemma~\ref{ddi} the proof of Theorem \ref{crfth1}
follows.
\end{proof}

We conjecture that the converse of the claim of Theorem~\ref{crfth1} is true. At this point we can
prove Lemma \ref{QR-lemma}, which supports the conjecture.  First we prove a useful
technical result.

\begin{lemma}\label{key1}\sideremark{key1}
Suppose $M$ is a quaternionic contact manifold of dimension $(4n+3)>7$. If $\psi$ is a smooth closed two-form
whose restriction to $H$ vanishes, then $\psi$ vanishes identically.
\end{lemma}
\begin{proof}[Proof of Lemma~\ref{key1}]
The hypothesis on $\psi$ show that $\psi$ is of the form $\psi\ =\
\sum_{s=1}^3\sigma_s\wedge\eta_s, $ where $\sigma_s$ are 1-forms.
 Taking the exterior differential and using \eqref{fifteen}, we obtain for $X\in H$

\begin{gather*}
0 =\sum_{a=1}^{4n}d\psi(e_a,I_ie_a,X) = -2\sum_{a=1}^{4n}\sum_{s=1}^3\sigma_s\wedge\omega_s(e_a,I_ie_a,X)=
(4n-2)\sigma_i(X)+2\sigma_j(I_kX)-2\sigma_k(I_jX),
\end{gather*}
where $e_1,\dots,e_{4n}$ is an orthonormal basis on $H$. For $n>1$
the latter implies $\sigma_{s_{|_H}}=0, s=1,2,3$. Hence, we have
$\psi=\sum_{1 \le s<t\le}A_{st}\,\eta_s\wedge\eta_t, $ where
$A_{st}$ are smooth functions on $M$. Now, the exterior derivative
gives $ 0=d\psi(e_a,I_se_a,\xi_t)\ =\ 2A_{st}.$
\end{proof}

The assumption in the next Lemma is a kind of $\partial\bar\partial_H$-lemma result,
which we do not know how to prove at the moment, but we believe that it is true. We show
how it implies the converse of Theorem \ref{crfth1}.

\begin{lemma}\label{QR-lemma}  Suppose, for  $i=1,2,3$,
$DD_{I_i}f\equiv dd_if+d_jd_kf=\sum_{s=1}^3p^i_s\omega_s \quad
\text{mod} \ \eta$ implies
\begin{equation}\label{dddd1}
dd_if-dd_jA_i\ =\ 2\sum_{s=1}^3r^i_s\, \omega_s \qquad \text{mod} \ \eta
\end{equation}
\noindent for some function $A_i$ on a QC manifold of dimension
$(4n+3)>7$. With this assumption, if $DD_{I_i}f =
\sum_{s=1}^3p^i_s\omega_s \quad \text{mod} \ \eta$, $i=1,2,3$, then
$f$ is a real part of an anti-CRF-function.
\end{lemma}
\begin{proof}[Proof of Lemma~\ref{QR-lemma}]
Consider the closed 2-forms\\
\centerline{$\Omega_i=d(d_if-d_jA_i-\sum_{s=1}^3r^i_s\eta_s).$}
We have $d\Omega_i=0$ and $\Omega_{i_{|_H}}=0$ due to \eqref{dddd1} and \eqref{thirteen}. Applying
Lemma~\ref{key1} we conclude $\Omega_i=0$, after which the Poincare lemma yields\\
\centerline{$d_if-d_jA_i-dB_i=0 \quad mod \quad \eta$}
\noindent for some smooth functions $A_1,A_2,A_3,B_1,B_2, B_3$.
%
The latter implies\\
\centerline{$df+d_1(A_2+B_1)+d_2(B_2-A_3)+d_3(A_1-B_3)=0\quad mod \quad \eta.$}
Setting $w=-A_2-B_1,\quad u=A_3-B_2,\quad v=B_3-A_1$ proves the claim.
\end{proof}

\begin{cor}\label{crfth13s}
Let $f:M\rightarrow \mathbb R$ be a smooth real function on a (4n+3)-dimensional $(n>1)$
3-Sasakian manifold $(M,\eta)$. If $f$ is the real part of an anti-CRF function
$f+iw+ju+kv$ then: 
\begin{enumerate}
\item[i)] The equation \eqref{crf03} holds.
\item[ii)] For
any $X,Y\in H$ we have the  equality
\begin{gather*}
(\nabla_Xdf)(Y)+(\nabla_{I_1X}df)(I_1Y)+
(\nabla_{I_2X}df)(I_2Y)+(\nabla_{I_3X}df)(I_3Y)=4\lambda g(X,Y).
\end{gather*}
\end{enumerate}
The function $\lambda$ is determined in \eqref{lamb1}.
\end{cor}

\begin{cor}\label{crfth13h}
Let $f:G(\mathbb H)\rightarrow \mathbb R$ be a smooth real function on the
(4n+3)-dimensional $(n>1)$ quaternionic Heisenberg group endowed with the standard flat
quaternionic contact structure and $\{T_a, X_a,Y_a,Z_a, \quad a=1,\dots,4n\}$
 be  $\nabla$-parallel basis on $G(\mathbb H)$.  If $f$ is the real part of an anti-CRF
 function $f+iw+ju+kv$ then the following
equivalent conditions hold true:
\begin{enumerate}
\item[i)] The equation \eqref{crf03} holds.
\item[ii)] The horizontal Hessian of $f$ is given by
\begin{gather*}\label{crf02h}
T_bT_af+X_bX_af+Y_bY_af+Z_aZ_b=4\lambda g(T_b,T_a);
\end{gather*}
\item[iii)] The function $f$ satisfies the following second order differential equation
\begin{equation*}
D_{T_b}\overline D_{T_a}f=\lambda (g-i\omega_1-j\omega_2-k\omega_3)(T_b,T_a);
\end{equation*}
\end{enumerate}
The function $\lambda$ is given by \eqref{lamb1}.
\end{cor}
Proposition~\ref{pseudoeinst}, Corollary~\ref{crfth13s} and
Example~\ref{3sas1} imply the next Corollary.
\begin{cor}\label{333s}
Let $(M,\eta)$ be a (4n+3)-dimensional $(n>1)$ 3-Sasakian manifold, $f:M\rightarrow
\mathbb R$ a positive smooth real function. Then the conformally 3-Sasakian QC structrure
$\bar\eta=f\eta$ is qc-pseudo Einstein if and only if the operators $DD_{I_s}f, s=1,2,3$
satisfy \eqref{crf03}. In particular, if $f$ is real part of anti CRF function then the
conformally 3-Sasakian qc structure $\bar\eta=f\eta$ is qc-pseudo Einstein.
\end{cor}

\section{Infinitesimal Automorphisms}
\subsection{3-contact manifolds}
We start with the more general notion of 3-contact manifold $(M,H)$, where $H$ is an
orientable codimension three distribution on $M$. Let $E\subset TM^*$ be the canonical
bundle determined by $H$, i.e. the bundle of 1-forms with kernel $H$. Hence, M is
orientable if and only if $E$ is also orientable, i.e. $E$ has a global non-vanishing
section  $vol_E$ locally given by $vol_E\ =\ \eta_1\wedge\eta_2\wedge\eta_3$. Denote by
$\eta=(\eta_1,\eta_2,\eta_3)$ the local 1 -form with values in $\mathbb R^3$. Clearly
$H=Ker~\eta$.

\begin{dfn}A $(4n+3)$-dimensional orientable smooth manifold $(M,\eta, H=Ker~\eta)$ is said to be
a \emph{3-contact manifold} if the restriction of each 2-form $d\eta_i,\  i=1,2,3$ to $H$
is non-degenerate, i.e.,
\begin{equation}\label{3connon}
d\eta_i^{2n}\wedge\eta_1\wedge\eta_2\wedge\eta_3 \ = \ u_i\ vol_{M}, \quad u_i>0,\quad
i=1,2,3,
\end{equation}
and the following compatibility conditions hold
\begin{equation}\label{3concom}
d\eta_1^p\wedge d\eta_2^q\wedge d\eta_3^r\wedge\eta_1\wedge\eta_2\wedge\eta_3=0, \quad p+q+r=2n, \quad 0<p,q,r,<2n.
\end{equation}
\end{dfn}

Denote the restriction of $d\eta_i$ on $H$ by $\Omega_i, \quad \Omega_i=(d\eta)_{|_H}, i=1,2,3$. The condition
\eqref{3connon} is equivalent to \\
\centerline{$\Omega_i^{2n}\not=0, \quad i=1,2,3, \qquad \Omega_1^p\wedge\Omega_2^q\wedge\Omega_3^r=0,
\quad p+q+r=2n, \quad 0<p,q,r<0.$}
We remark that the notion of 3-contact structure is slightly more general than the notion
of QC structure. For example, any real hypersurface $M$ in $\mathbb H^{n+1}$ with
non-degenerate second fundamental form carries 3-contact structure defined in the
beginning of Section~\ref{6.2.1}
(cf. Proposition~\ref{p:QRhypersurface} and Definition~\ref{d:QRhypersurface} where this
structure is QC if and only if \eqref{e:horizontal2forms} holds, or equivalently, the
second fundamental form is, in addition, invariant with respect to the hypercomplex
structure on $\mathbb H^{n+1}$). Another examples of 3-contact structure is the so called
quaternionic CR structure introduced in \cite{KN} and the so called weak QC structures
considered in \cite{D1}.  Note that in these examples the 1-form $\eta=(\eta_1,\eta_2,
\eta_3)$ are globally defined.

On any 3-contact manifold $(M,\eta, H)$ there exists a unique triple
$(\xi_1,\xi_2,\xi_3)$ of vector fields transversal to $H$ determined by the conditions\\
\centerline{$\eta_i(\xi_j)=\delta_{ij}, \quad (\xi_i\lrcorner d\eta_i)_{|_H}=0.$}
We refer to such a triple as fundamental vector fields or Reeb vector fields and denote
$V=span\{\xi_1,\xi_2,\xi_3\}$.
Hence, we have the splitting $TM=H\oplus V$.

The 3-contact structure $(\eta,H)$ and the vertical space $V$ are determined up to an action of $GL(3,\mathbb R)$, namely
for any $GL(3,\mathbb R)$ matrix $\Phi$ with smooth entries the structure $\Phi\cdot\eta$ is again a 3-contact structure.
Indeed, it is an easy algebraic fact that the conditions \eqref{3connon} and \eqref{3concom} also hold
for $\Phi\cdot\eta$. The Reeb vector field are transformed with the matrix with entries the adjunction
quantities of $\Phi$, i.e. with the inverse matrix $\Phi^{-1}$. This leads to the next
\begin{dfn}\label{d:3-ctct auto3}
A diffeomorphism $\phi$ of a 3-contact manifold $(M,\eta, H)$ is called a \emph{3-contact
automorphism} if $\phi$ preserves the 3-contact structure $\eta$, i.e.,
\begin{equation}\label{aut13c}
\phi^*\eta=\Phi\cdot\eta,
\end{equation}
for some matrix $\Phi\in GL(3,\mathbb R)$ with smooth
functions as entries and $\eta=(\eta_1,\eta_2,\eta_3)^t$ is considered as an element of
$\mathbb R^3$. \end{dfn}
The infinitesimal versions of
these notions lead to the following definition.

\begin{dfn}\label{d:3-ctct v field3}
A vector field $Q$ on a 3-contact manifold $(M,\eta, H)$ is an
infinitesimal generator of a 3-contact automorphism (3-contact
vector field) if its flow preserves the 3-contact structure, i.e.
\begin{equation*}
\LieQ\eta=\phi\cdot\eta,\qquad \phi\in gl(3,\mathbb R).
\end{equation*}
\end{dfn}
We show that any 3-contact vector field on a 3-contact manifold depend on 3-functions which satisfy
certain differential relations. We begin with describing infinitesimal automorfisms of
the 3-contact structure $\eta$ i.e. vector field $Q$ whose flow satisfies
\eqref{aut13c}. Our main observation is that 3-contact vector fields on a 3-contact
manifold are completely determined by their vertical components in the sense of
the following

\begin{prop}\label{qaut13c}
Let ($M,\eta,H)$ be a 3-contact manifold. A smoth vector field $Q$ on $M$ is 3-contact vector field
if and only if the functions  $f_i=\eta_i(Q), i=1,2,3$ satisfy the next compatibility conditions on $H$
\begin{equation}\label{compaut3c}
\begin{aligned}
& u_i(df_i+f_j(\xi_j\lrcorner d\eta_i)+f_k(\xi_k\lrcorner d\eta_i))_{|_H}\wedge\Omega_i^{(2n-1)})=\\
&u_j(df_j+f_k(\xi_k\lrcorner d\eta_j)+f_i(\xi_i\lrcorner d\eta_j))_{|_H}\wedge\Omega_j^{(2n-1)})
\quad {\text on}\quad H.
\end{aligned}
\end{equation}
The 3-contact vector field $Q$ has the form
\begin{equation*}
Q=Q_h+f_1\xi_1+f_2\xi_2+f_3\xi_3,
\end{equation*}
where $Q_h$ is the horizontal 3-contact hamiltonian field of $(f_1,f_2,f_3)$ defined on $H$ by
\begin{equation}\label{con3ham3c}
Q_h\lrcorner\eta_i=0, \quad Q_h\lrcorner(\Omega_i)=
-df_i-f_j(\xi_j\lrcorner d\eta_i)-f_k(\xi_k\lrcorner d\eta_i),\quad i=1,2,3, \quad {\text on} \quad H.
\end{equation}
\end{prop}

\begin{proof}
For a vector field $Q\in\Gamma(TM)$ we write $Q=Q_H+\sum_{s=1}^3\eta_s(Q)\xi_s$ where
$Q_H\in H$ is the horizontal part of $Q$. Applying \eqref{fifteen}, we
calculate
\begin{gather}\label{lieaut13c}
\LieQ\eta_i=Q\lrcorner d\eta_i+d(Q\lrcorner \eta_i)=\\\nonumber Q_H\lrcorner\Omega_i+[d(\eta_i(Q))
+\eta_i(Q)\xi_i\lrcorner d\eta_i+\eta_j(Q)\xi_j\lrcorner d\eta_i+\eta_k(Q)\xi_k\lrcorner
d\eta_i]_{|_H}\nonumber
\\\nonumber+[\xi_i(\eta_i(Q))-\eta_j(Q)d\eta_i(\xi_i,\xi_j)-\eta_k(Q)d\eta_i(\xi_i,\xi_k)]\eta_i
\\\nonumber+
[\xi_j(\eta_i(Q))+d\eta_i(Q,\xi_j)]\eta_j + [\xi_k(\eta_i(Q))+d\eta_i(Q,\xi_k)]\eta_k.
\end{gather}
Suppose $Q$ is a 3-contact vector field. Then \eqref{lieaut13c} and
the compatibility conditions \eqref{3concom} imply that $f_i$ and
$Q_H$ necessarily satisfy \eqref{compaut3c} and \eqref{con3ham3c},
respectively. Therefore $Q_H=Q_h$. The converse follows from
\eqref{lieaut13c} and the conditions of the proposition.
\end{proof}

\noindent The last Proposition implies that the space of 3-contact vector fields is
isomorphic to the space of triples consisting of smooth function $f_1,f_2,f_3$ satisfying
the compatibility conditions \eqref{compaut3c}.
\begin{cor}\label{connforms1} Let  $(M,\eta)$  be a 3-contact manifold. Then
\begin{enumerate}
\item[a)]
If $Q$ is a horizontal 3-contact vector field on $M$ then $Q$ vanishes identically.
\item[b)] The vector fields $\xi_i, i=1,2,3$ are 3-contact vector fields if and only if \\
\centerline{$ \xi_i\lrcorner d\eta_j  {|_H}=0, \quad i,j=1,2,3. $}
\end{enumerate}
\end{cor}

\subsection{QC vector fields}

{Suppose $(M, g, \mathbb{Q})$ is a quaternionic contact manifold.}

\begin{dfn}\label{d:3-ctct auto}
A diffeomorphism $\phi$ of a QC manifold $(M,[g],\mathbb Q)$ is
called a \emph{conformal quaternionic contact automorphism
(conformal qc-automorphism)} if $\phi$ preserves the QC structure,
i.e.
\begin{equation*}
\phi^*\eta=\mu\Psi\cdot\eta,
\end{equation*}
for some positive smooth function $\mu$ and some matrix $\Psi\in
SO(3)$ with smooth functions as entries and
$\eta=(\eta_1,\eta_2,\eta_3)^t$ is a local 1-form considered as an
element of $\mathbb R^3$.
\end{dfn}

In view of the uniqueness of the possible associated almost complex
structures, a quaternionic contact automorpism will preserve also
the associated (if any) almost complex structures, $\phi^* \mathbb
Q=\mathbb Q$ and consequently, it will preserve the conformal class
[g] on $H$. We note that {conformal} QC diffeomorphisms on
$S^{4n+3}$ are considered in \cite{Kam}. The infinitesimal versions
of these notions lead to the following definition.

\begin{dfn}\label{d:3-ctct v field}
{A vector field $Q$ on a QC manifold $(M, [g], \mathbb{Q})$} is an infinitesimal
generator of a conformal quaternionic contact automorphism (\emph{QC vector field} for
short) if its flow preserves the QC structure, i.e.
\begin{equation}\label{autvf}
\LieQ\eta=(\nu I+O)\cdot\eta,
\end{equation}
where $\nu$ is a smooth function and $O\in so(3)$.
\end{dfn}
In view of the discussion above a QC vector field on a QC manifold $(M,\eta,\mathbb Q)$
satisfies the conditions.
\begin{gather}\label{infaut1}
\LieQ g =\nu g,\\\label{qautf1}
\LieQ I=O\cdot I, \qquad O\in so(3), \quad I=(I_1,I_2,I_3)^t,
\end{gather}
If the flow of a vector field $Q$ is a conformal diffeomorphism of the horizontal metric $g$,
i.e. \eqref{infaut1} holds,  we shall call it \emph{ infinitesimal conformal
isometry}. If the function $\nu=0$ then $Q$ is said to be \emph{infinitesimal isometry}.

A QC vector field on a QC manifold is a 3-contact vector field of special type. Indeed,
let $\sharp$ be the musical isomorphism between $T^*M$ and $TM$ with respect to the fixed
Riemannian metric $g$ on $TM$ and recall that the forms $\alpha_j$ were defined in
\eqref{coneforms}. We have
\begin{prop}\label{qaut1}
{Let $(M, g, \mathbb{Q})$ be a quaternionic contact manifold}.  The vector field $Q$ is
a QC vector field if and only if
\begin{equation}\label{e:form of Q}
Q\ =\ \frac 12 \big ( f_j\, I_i \, \alpha^\sharp_k \ -\ f_k\, I_i \, \alpha^\sharp_j \ -\
I_i\,(df_i)^\sharp \big)\ +\ \sum_{s=1}^3 f_s\, \xi_s,
\end{equation}
for some functions $f_1, f_2$  and $f_3$ such that for any positive permutation $(i,j,k)$
of $(1,2,3)$ we have
\begin{gather}\label{e:o_ii=0}
f_j\, d\eta_i(\xi_j, \xi_i) +\ f_k\, d\eta_i(\xi_k, \xi_i)\ +\ \xi_i\, f_i\ =\ f_k\,
d\eta_j(\xi_k, \xi_j) +\ f_i\, d\eta_j(\xi_i, \xi_j)\ +\ \xi_j\, f_j\,
\\\label{e:o_ij=-o_ji}
 f_i\, d\eta_i(\xi_i, \xi_j) +\ f_k\, d\eta_i(\xi_k, \xi_j)\ +\ \xi_j\, f_i
 =\  -\  f_j\, d\eta_j(\xi_j, \xi_i)\ - \ f_k\,
d\eta_j(\xi_k, \xi_i)\  - \ \xi_i\, f_j,
\\\label{e:Q_H independ of i}
 f_j\, I_i(\alpha_k)^\sharp\ -\ f_k\, I_i
(\alpha_j)^\sharp\ -\ I_i(df_i)^\sharp\ =\ f_i\, I_k(\alpha_j)^\sharp\ -\ f_j\, I_k
(\alpha_i)^\sharp\ -\ I_k(df_k)^\sharp
\end{gather}
\end{prop}
\begin{proof}
Notice that \eqref{e:form of Q} implies $f_i=\eta_i(Q)$. By Cartan's formula
\eqref{autvf} is equivalent to
\[
Q  \lrcorner d\eta_i\ +\ df_i\ =\ \nu \, \eta_i\ + \ o_{is}\, \eta_s.
\]
In other words, both sides must be the same when evaluated on $\xi_t$, $t=1,2,3$ and also
when restricted to the horizontal bundle. Let $Q=Q_H + \sum_{s=1}^3f_s \xi_s$. Consider first the
action on the vertical vector fields.  Pairing with $\xi_t$ and taking successively $t=i,
j, k$ gives
\begin{gather}
f_j\, d\eta_i(\xi_j, \xi_i) +\ f_k\, d\eta_i(\xi_k, \xi_i)\ +\ \xi_i\, f_i\ =\ \nu \ +\
o_{ii}\nonumber\\\label{e:relation between Q and o}
\alpha_k (Q_H)\ +\ f_i\, d\eta_i(\xi_i, \xi_j) +\ f_k\, d\eta_i(\xi_k, \xi_j)\ +\ \xi_j\, f_i\ =\ o_{ij}
\\\nonumber
-\alpha_j(Q_H)\ +\ f_i\, d\eta_i(\xi_i, \xi_k) +\ f_j\, d\eta_i(\xi_j, \xi_k)\ +\ \xi_k\,
f_i\ =\ o_{ik}.
\end{gather}
 Equating the restrictions to the horizontal bundle, i.e.,
$d\eta_i(Q, .)|_H\ +\ df_i|_H\ =\ 0$, gives
\[
\Big (f_j\, d\eta_i(\xi_j, .)\ +\ f_k\, d\eta_i(\xi_k, .)\ +\ d\eta_i(Q_H, .)\ +\ df_i
\Big )|_H\ =\ 0.
\]
Since $g(A,.)|_H\ = \ 0 \Leftrightarrow A\ =\ \sum_{s=1}^3\eta_s(A)\xi_s$, the last equation is
equivalent to
\begin{equation}\label{e:horizontal part}
-f_j\alpha^\sharp_k\ +\ f_k \alpha^\sharp_j\ -\ 2I_i Q_H \ +\ (df_i)^\sharp \ =\ \sum_{s=1}^3\big (
-f_j\alpha_k (\xi_s)\ +\ f_k \alpha_j(\xi_s)\ +\ \xi_s f_i \big )\, \xi_s.
\end{equation}
Acting with $I_i$ determines $2Q_H \ = \ f_j\, I_i(\alpha_k)^\sharp\ -\ f_k\, I_i
(\alpha_j)^\sharp\ -\ I_i(df_i)^\sharp$, which implies \eqref{e:form of Q}. In addition
we have
\[
\alpha_j(Q_H)\ =\ \frac 12 \Big ( f_j\, \alpha_j(I_i(\alpha_k)^\sharp)\ -\ f_k\, \alpha_j
(I_i (\alpha_j)^\sharp)\ -\ \alpha_j(I_i(df_i)^\sharp) \Big )
\]
On the other hand, $o\in so(3)$ is equivalent to $o$ being a skew symmetric which is
equivalent to \eqref{e:o_ii=0} and \eqref{e:o_ij=-o_ji}, by the above computations.
Therefore, if we are given three functions $f_1, \ f_2, \ f_3$ satisfying
\eqref{e:o_ii=0}, \eqref{e:o_ij=-o_ji} and \eqref{e:Q_H independ of i}, then we define
$Q$ by \eqref{e:form of Q}. Using \eqref{e:relation between Q and o} we define $\nu$ and
$o$ with $o\in so(3)$. With these definitions $Q$ is a QC vector field.
\end{proof}

Using the formulas in Example~\ref{3sas} we obtain from Proposition~\ref{qaut1} the
following '3-hamiltonian' form of a QC vector field on 3-Sasakian manifold.
\begin{cor}\label{connforms}
{Let $(M,\eta)$ be a 3-Sasakian manifold}. Then any QC vector field $Q$ has the form
\begin{equation*}
Q=Q_h+f_1\xi_1+f_2\xi_2+f_3\xi_3,
\end{equation*}
where the smooth functions $f_1,f_2,f_3$  satisfy the conditions
\begin{equation*}
d_if_i=d_jf_j, \quad \xi_i(f_i)=\xi_j(f_j), \quad \xi_i(f_j)=-\xi_j(f_i),\quad i,j=1,2,3,
\end{equation*}
and the horizontal part $Q_h\in H$ is determined by
$$ Q_h\lrcorner d\eta_i=d_if_i, \quad  i\in\{1,2,3\}.
$$
The matrix in \eqref{autvf} has the form
\begin{gather*}
\nu I_{d_3}+O=\left(\begin{array}{ccc} \xi_1(\eta_1(Q)) & -\xi_1(\eta_2(Q))-2\eta_3(Q)
& -\xi_1(\eta_3(Q))+2\eta_2(Q)\\
\xi_1(\eta_2(Q))+2\eta_3(Q) & \xi_1(\eta_1(Q)) &
 -\xi_2(\eta_3(Q))-2\eta_1(Q)\\
\xi_1(\eta_3(Q))-2\eta_2(Q) & \xi_2(\eta_3(Q))+2\eta_1(Q) & \xi_1(\eta_1(Q))
\end{array}\right).
\end{gather*}
In particular, the Reeb vector fields $\xi_1,\xi_2,\xi_3$ are 3-contact vector
fields.
\end{cor}
Corollary~\ref{connforms1} tells us that on a QC manifold the Reeb vector fields
$\xi_1,\xi_2,\xi_3$ are 3-contact exactly when the connection 1-forms vanish on $H$. This
combined with Corollary~\ref{3-sas} gives 3-Sasakian structure compatible with the given
3-contact structure $H$, if the qc-scalar curvature is not zero (see Corollary~\ref{3con3sas} below).

First we shall investigate some useful properties of a QC vector field.
\begin{prop}\label{3coninf} {Let  $(M,[g],\mathbb Q)$ be QC manifold} and $Q$ be a QC vector field
determined by \eqref{autvf} and \eqref{e:relation between Q and o}.
The next equality hold
$$d\eta_i([Q,I_iX]^{\perp},Y)+d\eta_i(I_iX,[Q,Y]^{\perp})=0$$
\end{prop}
\begin{proof}
We have using \eqref{autvf} that
\begin{multline}\label{lieq?}
\LieQ d\eta_i(I_iX,Y)
=2(\LieQ \omega_i)(I_iX,Y)-d\eta_i([Q,I_iX]^{\perp},Y)-d\eta_i(I_iX,[Q,Y]^{\perp})\\
=-2(\LieQ g)(X,Y)+2g((\LieQ I_i)I_iX,Y) -d\eta_i([Q,I_iX]^{\perp},Y)-d\eta_i(I_iX,[Q,Y]^{\perp})\\
=(d\LieQ\eta_i)(I_iX,Y)=(d\nu\wedge\eta_i+\nu d\eta_i+do_{ij}\wedge\eta_j+o_{ij}d\eta_j+
do_{ik}\wedge\eta_k+o_{ik} d\eta_k)(I_iX,Y)\\= -2\nu g(X,Y)-2
o_{ij}\omega_k(X,Y)+2o_{ik}\omega_j(X,Y),
\end{multline}
\noindent where $o_{st}$ are the entries of the matrix $O$ given by
\eqref{e:relation between Q and o}.
Apply \eqref{infaut1} and \eqref{qautf1}
to \eqref{lieq?} to get the assertion.
\end{proof}

We are going to characterize the vanishing of the torsion of the Biquard connection in
terms of the existence of some special vertical vector fields. More precisely, we have
the following Theorem.

\begin{thrm}\label{flatqc}
{Let $(M,g,\mathbb Q)$ be a QC manifold} with non zero qc-scalar curvature. The following
conditions are equivalent
\begin{enumerate}
\item[i)] Each Reeb vector field is a QC vector field;
\item[ii)] The torsion of the Biquard connection is identically zero;
\item[iii)] Each Reeb vector field preserves the horizontal metric and the quaternionic structure
simultaneously, i.e. \eqref{infaut1} with $\nu=0$ and \eqref{qautf1} hold for $Q=\xi_i, i=1,2,3$;
\item[iv)] There exists a local 3-Sasakian structure in the sense of Theorem~\ref{Ein2MO}
\end{enumerate}
\end{thrm}

\begin{proof}
In the course of the proof we shall prove two Lemmas of independent
interest. Given a vector field $Q$, we define the symmetric tensor
$T^0_Q$ and the skew-symmetric tensor $u_Q$
\begin{equation}\label{e:u_Q}
T^0_Q \ = \
\sum_{s=1}^3\eta_s(Q)\,T^0_{\xi_s} , \qquad
u_Q\ =\ \sum_{s=1}^3\eta_s(Q)\, I_su, \end{equation} respectively,
such that, \hspace{4mm} $T(Q,X,Y)\ =\ g(T^0_Q
X, Y)\ +\ g(u_Q X, Y),$

\begin{lemma}\label{l:13part}
The tensors $T^0_Q$ and $u_Q$ lie in the  $[-1]$ component associated to the operator
$\dag$ cf. \eqref{e:cross} and \eqref{New21}.
\end{lemma}
\begin{proof}[Proof of Lemma~\ref{l:13part}]
By the definition of $u_Q$, we have
$
g(u_Q I_1X, I_1Y)\ =\ \sum_{s=1}^3\eta_s(Q)\,g(I_s uX,Y)\
$
and after summing we find\\
\centerline{$\sum_{j=1}^3 g(u_Q I_jX, I_jY)\
=\sum_{j=1}^3\-\eta_j(Q)\,g(I_j uX, Y)\ =\ -g(u_Q X, Y). $}
We turn to the second claim.  Recall that $T^0_{\xi_j}$ anti-commutes with $I_j$, see
\eqref{tors1}. Hence,
\begin{multline}
g(T^0_Q I_1X, I_1Y)\ =\ -\eta_1(Q)\,g(T^0_{\xi_1}X, Y)
    -\  \eta_2(Q)\,[g(T^{0\ --+}_{\xi_2}X, Y)\ -g(T^{0\ +--}_{\xi_2}X,
    Y)]\\ \notag
    -\ \eta_3(Q)\,[g(T^{0\ -+-}_{\xi_2}X, Y)\ -g(T^{0\ +--}_{\xi_3}X,
    Y)],\notag\\
g(T^0_Q I_2X, I_2Y)\ =\ -\eta_2(Q)\,g(T^0_{\xi_2}X, Y)
    -\  \eta_1(Q)\,[g(T^{0\ --+}_{\xi_1}X, Y)\ -g(T^{0\ -+-}_{\xi_1}X,
    Y)]\\ \notag
    -\ \eta_3(Q)\,[g(T^{0\ +--}_{\xi_3}X, Y)\ -g(T^{0\ -+-}_{\xi_3}X,
    Y)],\notag\\
g(T^0_Q I_3X, I_3Y)\ =\ -\eta_3(Q)\,g(T^0_{\xi_3}X, Y)
    -\  \eta_1(Q)\,[g(T^{0\ -+-}_{\xi_1}X, Y)\ -g(T^{0\ --+}_{\xi_2}X,
    Y)]\\ \notag
    -\ \eta_2(Q)\,[g(T^{0\ +--}_{\xi_2}X, Y)\ -g(T^{0\ --+}_{\xi_2}X,
    Y)].\notag
\end{multline}
\noindent Summing the above three equations we come to
\[
\sum_{j=1}^3g(T^0_Q I_jX, I_jY)\ =\ -\sum_{j=1}^3g(Q,\xi_j)\,g(T^0_{\xi_j}X, Y)\ =\ -
g(T^0_Q X, Y),
\]
which finishes the proof of Lemma~\ref{l:13part}.
\end{proof}

\begin{lemma}\label{p:Q infinit vf}
If $Q$ is an infinitesimal conformal isometry whose flow preserves
the quaternionic structure then
the next two equalities hold
\begin{gather}\label{cau2}
g(\nabla_XQ,Y)\ +\ g(\nabla_YQ,X)\ + \ 2g(T^0_Q X,Y)=\nu\, g(X, Y),
\\\label{cau3} 3g(\nabla_X Q, Y)-\sum_{s=1}^3g(\nabla_{I_sX} Q,I_sY) +
     4g(T^0_QX,Y) +4g(u_QX,Y)=-2\sum_{(ijk)}L_{ij}(Q)\omega_k(X,Y),
\end{gather}
where the sum is over all even permutation of (1,2,3) and
\begin{multline}\label{newlij}
L_{ij}(Q)=-L_{ji}(Q)=\xi_j(\eta_i(Q))-\eta_j(Q)d\eta_j(\xi_i,\xi_j)
\\+\frac12\eta_k(Q)\left(\frac{Scal}{8n(n+2)}
+d\eta_j(\xi_k,\xi_i)-d\eta_i(\xi_j,\xi_k)-d\eta_k(\xi_i,\xi_j)\right).
\end{multline}
\end{lemma}

\begin{proof}[Proof of Lemma \ref{p:Q infinit vf}]
In terms of the Biquard connection \eqref{infaut1}  reads exactly as \eqref{cau2}.
Furthermore, from \eqref{qautf1}, \eqref{e:relation between Q and o} and \eqref{der} it follows
\begin{multline}\label{infqua}
o_{ij}I_jX+o_{ik}I_kX =\ (\LieQ I_i)(X) =\\
    =\ -\nabla_{I_iX} Q\ +I_i\nabla_XQ-\alpha_j(Q)I_kX+\alpha_k(Q)I_jX
    -\ T(Q, I_i X)\ +\ I_iT(Q, X).
\end{multline}
\noindent A use of \eqref{e:relation between Q and o},
\eqref{coneforms} and \eqref{coneform1} allows us to write the last
equation in the form
\begin{multline}\nonumber
g(\nabla_XQ,Y)\ -\ g(\nabla_{I_iX} Q, I_i Y)
      +\ T(Q,X,Y)\ -\ T(Q,I_iX, I_iY)\ \\=(o_{ij}-\alpha_k(Q))\omega_k(X,Y)-(o_{ik}+\alpha_j(Q))\omega_j(X,Y)\\
=-L_{ij}(Q)\omega_k(X,Y)+L_{ik}(Q)\omega_j(X,Y),
\end{multline}
where $L_{ij}(Q)$ satisfy \eqref{newlij}. Summing the above
identities for the three almost complex structures and applying
Lemma~\ref{l:13part},  we obtain \eqref{cau3}, which completes the
proof of Lemma~\ref{p:Q infinit vf}.
\end{proof}

We are ready to  finish the proof of Theorem~\ref{flatqc}. Let $\xi_i, i=1,2,3$ be QC vector fields.
Then \eqref{cau2} for $Q=\xi_i$ yields $T_{\xi_i}=0, i=1,2,3,\quad \nu=0$ since $T_{\xi_i}$ is trace-free.
Consequently, for any cyclic permutation $(i,j,k)$ of $(1,2,3)$, \eqref{cau3} and \eqref{newlij} imply\\
\centerline{$u_{\xi_i}=0,\quad d\eta_j(\xi_i,\xi_j)=0,  \quad d\eta_i(\xi_j,\xi_k)=\frac{Scal}{8n(n+2)}$}
by comparing the trace and the trace-free part. Hence ii) follow.

Conversely, if the torsion of the Biquard connection vanishes, then \eqref{cau2} is trivially satisfied
for $\nu=0$ and \eqref{infqua} yields \eqref{qautf1} with $o_{ij}=\alpha_k$.
This establishes the equivalence between ii) and iii).

The other equivalences in the theorem follow from Theorem~\ref{Ein2MO},
Example~\ref{3sas1}, Corollary~\ref{connforms1} and Corollary~\ref{3-sas}.
\end{proof}

\begin{cor}\label{3con3sas}
{Let $(M,g,\mathbb Q)$ be a QC manifold} with non zero qc-scalar curvature. The following
conditions are equivalent
\begin{enumerate}[i)]
\item There exists a local 3-Sasakian structure compatible with $H=Ker~\eta$;
\item There are three linearly independent transversal QC-vector fields.
\end{enumerate}
\end{cor}
\begin{proof}
Let $\gamma_1,\gamma_2,\gamma_3$ be
linearly independent transversal QC-vector fields. Then there exist 1-forms
$\eta_{\gamma_1},\eta_{\gamma_2},\eta_{\gamma_3}$ satisfying
$\eta_{\gamma_i}(\gamma_j)=\delta_{ij}$, where $\delta_{ij}$ is the Kroneker symbol.
In view of the proof of Theorem~\ref{flatqc} it is sufficient to show $\gamma_1,\gamma_2,\gamma_3$
are the Reeb vector field for $\eta_{\gamma}$, i.e. we have to show
 that the compatibility conditions \eqref{bi1} are satisfied.
Indeed, the fact that $\gamma_i, i=1,2,3$ are  QC vector fields
means that \eqref{autvf} hold with respect to $\eta_{\gamma}$. Then \eqref{e:relation between Q and o}
gives $\nu=0$ and the second line of \eqref{lieaut13c}, for $\eta_{\gamma}$ and $Q=\gamma_i, i=1,2,3$,
imply \eqref{bi1} for the structure $\eta_{\gamma}$. Theorem~\ref{flatqc} completes the proof.
\end{proof}

 In the particular case when the vector field $Q$ is the
gradient of a function defined on the manifold $M$, we have the
following formulas.

\begin{cor}\label{c:gradient infinit vf}
If $h$ is a smooth real valued function on $M$ and $Q\ =\ \nabla h$
is a QC  vector field,
then for any horizontal vector fields $X$ and $Y$ we have
\begin{enumerate}[a)]
    \item $ [\ (\nabla dh)]_{[3][0]}(X,Y)\ =\ 0$
    \item $ [\ \nabla dh\ ]_{[sym][-1]}(X,Y)\ =\ -T^0_{Q}(X,Y)\ \ ( \text{ cf. \eqref{e:u_Q} })\ $
    \item $\ \  \  u_{Q}(X,Y)\ =\ 0 \ \ ( \text{ cf. \eqref{e:u_Q} }),\qquad  L_{ij}(\nabla h)=0.$
\end{enumerate}
\end{cor}

\begin{proof}
Use \eqref{cau2} and \eqref{symdh} to get
 $
 2\nabla dh (X,Y)\ +\ 2dh(\xi_j)\, \omega_j(X,Y)\ +\ 2g(T^0_Q X, Y)\ =\ \nu
g(X,Y). $ Decomposing in the $[-1]$ and   $[3]$ components completes
the proof of a) and b), taking  into account \eqref{cau3} and
Lemma~\ref{l:13part}. The skew-symmetric part of \eqref{cau3} gives
$2u_{Q}\ +\ \sum_{(ijk)}L_{ij}(\nabla h)\omega_k\ =\ 0,$ where the
sum is over all even permutations of $(1,2,3)$. Hence, c) follows by
comparing the trace  and trace-free parts of the last equality.
\end{proof}

\section{Quaternionic contact Yamabe problem}
\subsection{The Divergence Formula}

{Let $(M, \eta)$ be a quaternionic contact manifold with a fixed globally defined contact
form $\eta$. }For a fixed $j\in \{1,2,3\}$ the form
\begin{equation}\label{e:volumeform}
Vol_{\eta}=\eta_1\wedge\eta_2\wedge\eta_3\wedge\omega_j^{2n}
\end{equation}
 is a volume form. Note that $Vol_{\eta}$ is independent
of $j$.
We define the (horizontal) divergence of a horizontal vector
field/one-form $\sigma\in\Lambda^1\, (H)$ by
\begin{equation}\label{e:divergence}
\nabla^*\, \sigma\ =tr|_{H}\nabla\sigma=\
\sum_{a=1}^{4n}(\nabla_{e_\alpha}\sigma)(e_\alpha).
\end{equation}
Clearly the horizontal divergence does not depend on the basis and
is an $Sp(n)Sp(1)$-invariant. For any horizontal 1-form $\sigma\in\Lambda^1\, (H)$ we denote
with $\sigma^{\#}$ the corresponding horizontal vector field via
the horizontal metric defined with the equality
$\sigma(X)=g(\sigma^{\#},X)$. It is justified to call the function $\nabla^*\, \sigma$
divergence of $\sigma$ in view of the following Proposition.
\begin{prop}\label{p:divergence lemma} Let $(M,\eta)$ be a quaternionic contact manifold of dimension
(4n+3) and $\eta\wedge\omega_s^{2n-1}\ \overset {def}{=}\
\eta_1\wedge\eta_2\wedge\eta_3\wedge\omega_s^{2n-1}.$ For any horizontal 1-form $\sigma\in\Lambda^1\, (H)$ we have
$$d(\sigma^{\#}\lrcorner (Vol_{\eta}))\ =\
-(\nabla^*\sigma)\,\eta\wedge\omega^{2n}\ .$$ Therefore, if $M$ is
compact,
\[
\int_M (\nabla^*\sigma)\,\eta\wedge\omega^{2n}\ =\ 0.
\]
\end{prop}
\begin{proof}
We work in a qc-normal frame at a point $p\in M$ constructed in
Lemma~\ref{norma}. Since $\sigma$ is horizontal, we have
$\sigma^{\#}=g(\sigma^{\#},e_a)e_a$. Therefore, we calculate
$$\sigma^{\#}\lrcorner
(Vol_{\eta})=\sum_{a=1}^{4n}(-1)^ag(\sigma^{\#},e_a)\eta\wedge
e^{\#}_{a_1}\wedge\dots\wedge\hat{e^{\#}_a}\wedge\dots\wedge e^{\#}_{4n},$$ where
$\hat{e^{\#}_a}$ means that the 1-form $e^{\#}_a$ is missing in the above wedge
product. The exterior derivative of the above expression gives
\begin{equation*}d(\sigma^{\#}\lrcorner
(Vol_{\eta}))=-\sum_{a=1}^{4n}e_ag(\sigma^{\#},e_a)Vol_{\eta}=-(\nabla^*\sigma)\,Vol_{\eta}.
\end{equation*}
Indeed, since the Biquard connection preserves the metric,
the middle term is calculated as follows
$e_ag(X,e_a)=g(\nabla_{e_a}X,e_a)+g(X,\nabla_{a_a}e_a)$ which
evaluated at the point $p$ gives
$$e_ag(X,e_a)_{|_p}=g(\nabla_{e_a}X,e_a)_{|_p})=((\nabla_{e_a}\sigma)e_a)_{|_p}.$$
In order to obtain the last equality we also used the definition of the Reeb vector fields, \eqref{bi1}, and the following sequence of
identities
$$de^{\#}_b(e_b,e_a)_{|_p}=e^{\#}_b([e_b,e_a])_{|_p}=
e^{\#}_b(\nabla_{e_b}e_a-\nabla_{e_a}e_b-T(e_b,e_a))_{|_p}=0$$
since $T(e_b,a_a)$ is a vertical vector field.
This proves the first formula. If the manifold is compact, then
Stoke's theorem completes the proof.
\end{proof}
We note that the integral formula of the above theorem was
essentially proved in [\cite{Wei}, Proposition 2.1].

\subsection{Partial solutions of the QC-Yamabe problem}
In this Section we shall present a partial solution of the Yamabe problem on the
quaternionic sphere.  Equivalently, using the Cayley transform this provides a partial
solution of the Yamabe problem on the quaternionic Heisenberg group. The extra assumption
under which we classify the solutions of the Yamabe equation consists of assuming that
the "new" quaternionic structure has an integrable vertical space. The change of the
vertical space is given by \eqref{New19}. Of course, the standard quaternionic contact
structure has an integrable vertical distribution.  A note about the Cayley transform is
in order. We shall define below the explicit Cayley transform for the considered case,
but one should keep in mind the more general setting of groups of Heisenberg type
\cite{CDKR}. In that respect, the solutions of the Yamabe equation on the quaternionic
Hesenberg group, which we describe, coincide with the solutions on the groups of
Heisenebrg type \cite{GV}.

As in Section 5 we are considering a conformal transformation
$\tilde\eta=\frac{1}{2h}\eta$, where $\tilde\eta$ represents a fixed
quaternionic contact structure and $\eta$ is the "new" structure conformal
to the original one. In fact, eventually, $\tilde\eta$ will stand for the
standard quaternionic contact structure on the quaternionic sphere.  The
Yamabe problem in this case is to find all structures $\eta$, which are
conformal to $\tilde\eta$ and have constant scalar curvature equal to
$16n(n+2)$, see Corollary \ref{3sas}. The Yamabe equation is given by
\eqref{e:conf change scalar curv} and the problem is to find all solutions
of this equation.

\begin{prop}\label{qcYamab}
{Let $(M,\bar\eta)$ be a compact} qc-Einstein manifold of dimension $(4n+3)$. Let
$\bar\eta=\frac{1}{2h}\eta$ be a conformal deformation of the qc-structure $\bar\eta$ on
$M$. Suppose $\eta$ has constant scalar curvature.
\begin{enumerate}
\item[a)] If $n>1$, then any one of the following two
conditions
\begin{enumerate}
\item[i)] the
vertical space of $\eta$ is integrable;
\item[ii)] the QC structure $\eta$ is qc-pseudo Einstein;
\end{enumerate} implies that  $\eta$ is a qc-Einstein structure.
\item[b)]{ If $n=1$ and the vertical space of $\eta$ is integrable than
$\eta$ is a qc-Einstein structure.}
\end{enumerate}
\end{prop}
\begin{proof}
The proof follows the steps of the solution of the Riemannian Yamabe problem on the
standard unit sphere, see \cite{LP}. Theorem~\ref{Ein2MO} tells us that $\tilde\eta$ is a
qc-Einstein structure. Theorem \ref{sixtyseven} and equations \eqref{New24}, \eqref{e:T^o
conf change}, and \eqref{e:U conf change} imply
\begin{align}\label{ric-1}
[Ric_0]_{[-1]}(X,Y)\ & =\  (2n+2)T^0(X,Y)\ =\ -(2n+2)h^{-1}[\nabla
dh]_{[sym][-1]}(X,Y)\\\label{ric-3} [Ric_0]_{[3]}(X,Y)\ & =\  2(2n+5)U(X,Y)\ =\
-(2n+5)h^{-1}[\nabla dh-2h^{-1}dh\otimes dh]_{[3][0]}(X,Y).
\end{align}
\noindent Furthermore, when  the scalar curvature of $\eta$ is a constant
then Theorem \ref{t:horizontal system} gives
\begin{equation}\label{hordiv}
\nabla^*T^0=(n+2)\mathbb A, \qquad \nabla^* U = \frac{(1-n)}{2}\mathbb A.
\end{equation}
If $n>1$ and  either the vertical space of $\eta$ is an integrable
distribution or $\eta$ is qc-pseudo Einstein $U=0$, then  \eqref{hordiv}
 shows that $\mathbb A =0$ and the divergences of
$T^0$ and $U$ vanish $\nabla^*\, T^0\ =\ 0$ and $\nabla^*\, U\ =\ 0$. The
same conclusion can be reached in the case $n=1$ assuming the
integrability of the vertical space (recall that always $U=0$ when $n=1$).
We shall see that, in fact, $T^0$ and $U$ vanish, i.e., $\eta$ is also
qc-Einstein. Consider first the $[-1]$ component. Taking norms,
multiplying by $h$ and integrating, the divergence formula gives
\begin{multline*}
\int_{M} h\  \abs {\, [Ric_o]_{[-1]}\ }^2\, \eta\wedge\omega^{2n} \ =\ (2n+2)\ \int
\langle [ Ric_o ]_{[-1]}, \nabla dh] \rangle  \,
\eta\wedge\omega^{2n}\\
 =(2n+2) \int_{M} \langle \nabla^*\ [ Ric_o ]_{[-1]},
\nabla h] \rangle  \, \eta\wedge\omega^{2n}\ =\ 0.\hskip1truein
\end{multline*}
Thus, the $[-1]$ component of the qc-Einstein tensor vanishes $\abs {\,  [Ric_o]_{[-1]}\
}\ =\ 0$.  Define $h=\frac {1}{2u}$, inserting \eqref{transf} into \eqref{ric-3} one gets
$$[Ric_0]_{[3]}\ =\  2(2n+5)U\ =\ -(2n+5)[\nabla
du]_{[3][0]},$$ from where, arguing as before we get $[Ric_0]_{[3]}\ =\ 0$.
Theorem~\ref{Ein2MO} completes the proof.
\end{proof}

\begin{cor}
Let $\bar\eta=\frac{1}{2h}\eta$ be a conformal deformation of a compact qc-Einstein
manifold of dimension $(4n+3)$ and suppose $\bar\eta$ has  constant qc-scalar curvature.
\begin{enumerate}
\item[i)] { If $n>1$
 and either the gradient $\nabla h$ or the gradient
$\nabla(\frac1h)$ is a QC  vector fields then $h$ is a constant.}
\item[ii)] { If $n=1$ and the gradient
$\nabla(\frac1h)$ is a QC  vector fields then $h$ is a constant.}
\end{enumerate}
\end{cor}
\begin{proof} Suppose $\nabla h$ is a QC-vector field. Corollary~\ref{c:gradient infinit vf}, b) yields
$[\nabla dh]_{[sym][-1]}=0$ since the torsion of Biquard connection vanishes due to
Proposition~\ref{tor-ein}. Then Proposition~\ref{qcYamab} and a) in
Corollary~\ref{c:gradient infinit vf} imply that on $H$ we have\\
\centerline{$dh\otimes
dh+d_1h\otimes d_1h+d_2h\otimes d_2h+d_3h\otimes d_3h=\frac{|dh|^2}{n}g.$} { If $n>1$
then $dh_{|H}=0$, which implies $dh=0$ using the bracket generating condition.}

Suppose $\nabla(\frac1h)$ is a QC vector field. Then Proposition~\ref{qcYamab},
\eqref{transf} combined with b) in Corollary~\ref{c:gradient infinit vf} show that on $H$
we have\\
\centerline{$3dh\otimes dh-d_1h\otimes d_1h-d_2h\otimes d_2h-d_3h\otimes d_3h=0.$} Define
$X=I_1X, Y=I_1Y$ etc. to get $dh\otimes dh=d_1h\otimes d_1h=d_2h\otimes d_2h=d_3h\otimes
d_3h$. Hence, $dh_{|H}=0$ since dim $Ker\, dh\ =\ 4n-1$ and $dh=0$ as above.
\end{proof}

\subsection{Proof of Theorem~\ref{t:Yamabe}}
\begin{proof}
We start the proof with the observation that from
Proposition~\ref{qcYamab} and Corollary \ref{333s} the new structure
$\eta$ is also qc-Einstein. Next we bring into consideration the
quaternionic Heisenberg group. Let us identify $\QH$ with the boundary
$\Sigma$ of a Siegel domain in $\Hn\times\mathbb{H}$,
\[
\Sigma\ =\ \{ (q',p')\in \Hn\times\mathbb{H}\ :\ \Re {\ p'}\ =\ \abs{q'}^2 \},
\]
by using the map $(q', \omega')\mapsto (q',\abs{q'}^2 - \omega')$. The standard contact
form, written as a purely imaginary quaternion valued form, is  given by (cf.
\eqref{e:Heisenbegr ctct forms}) \hspace{3mm}
$2\tilde{\Theta} =
(d\omega \ - \ q' \cdot d\bar q' \ + \ dq'\, \cdot\bar q'),$ \hspace{3mm}
where $\cdot$ denotes the quaternion multiplication. Since \hspace{3mm}
$dp'\ =\ q'\cdot d\bar q'\ +\ dq'\, \cdot\bar {q}'\ -\ d\omega',
$ 
\hspace{3mm} under the identification  of $\QH$ with $\Sigma$ we have also \hspace{3mm}
$2\tilde{\Theta}\ =\ - dp'\ +\ 2dq'\cdot\bar {q}'.
$ 
\hspace{3mm}
Taking into account that $\tilde{\Theta}$ is purely imaginary, the last equation can be
written also in the following form
\[
4\,\tilde{\Theta}\ =\ (d\bar p'\ -\ d p')\ +\ 2dq'\cdot\bar {q'}\ -\ 2q'\cdot
d\bar q'.
\]
Now, consider the Cayley transform as the map \hspace{3mm}
$\mathcal{C}:S\mapsto \Sigma$ \hspace{3mm}
from the sphere $S\ =\ \{\abs{q}^2+\abs{p}^2=1 \}\subset \Hn\times\mathbb{H}$ minus a
point to the Heisenberg group $ \Sigma$, with $\mathcal{C}$ defined by
\[
 (q', p')\ =\ \mathcal{C}\ \Big ((q, p)\Big), \qquad
q'\ =\ (1+p)^{-1} \ q, \qquad p'\ =\ (1+p)^{-1} \ (1-p)
\]
\noindent and with an inverse map $(q, p)\ =\ \mathcal{C}^{-1}\Big ((q', p')\Big)$ given
by
\[
q\ =\ \ 2(1+p')^{-1} \ q', \qquad  p\ =\ (1+p')^{-1} \ (1-p').
\]
\noindent The Cayley transform maps $S$ minus a point to $\Sigma$ since
$$
\Re {\ p'}\ =\ \Re { \frac {(1+\bar p) (1-p)} {\abs {1+p\,}^2}
 }
 \ =\ \Re { \frac {1- \abs  {p} } {\abs
{1+ p\,}^2} }\ =\ \frac {\abs{q}^2}{\abs {1+p\,}^2}\ =\ \abs {q'}^2.
$$
Writing the Cayley transform in the form \hspace{3mm}
$
(1+p)q'\ =\  \ q, \quad (1+p)p'\ =\  1-p,
$
\hspace{3mm} gives
\[
dp\cdot q'\ +\ (1+p)\cdot dq'\ =\ d q, \hskip.5truein dp\cdot p'\ +\ (1+p)\cdot dp'\ =\
-dp,
\]
from where we find
\begin{equation}\label{e:dp'}
\begin{aligned}
dp'\ & =\ -2(1+p)^{-1}\cdot dp \cdot (1+p)^{-1}\\   dq' \ & =\ (1+p)^{-1}\cdot [ dq\ -\
dp\cdot (1+p)^{-1}\cdot q ].
\end{aligned}
\end{equation}
\noindent The Cayley transform is a conformal\ quaternionic contact
diffeomorphism between the quaternionic Heisenberg group with its
standard quaternionic contact structure $\tilde\Theta$ and the
sphere minus a point with its standard structure $\tilde\eta$, a
fact which can be seen as follows. Equations \eqref{e:dp'} imply the
following identities
\begin{multline}\label{e:Cayley transf ctct form}
2\,\mathcal{C}^*\,\tilde{\Theta}\ =\ -(1+\bar p)^{-1}\cdot d\bar p
\cdot (1+\bar p)^{-1}\ +\ (1+ p)^{-1}\cdot dp \cdot  (1+ p)^{-1}\\
+\ (1+p)^{-1}\, [dq\ -\ dp\cdot (1+p)^{-1}\cdot q ]\cdot \bar q\cdot (1+\bar p)^{-1}\\ -\
(1+p)^{-1}\, q\cdot [d\bar q\
 -\ \bar q\cdot (1+\bar
p)^{-1}\cdot
d\bar p\, ]\cdot (1+\bar p)^{-1}\\
=\ (1+p)^{-1}\, \Big [dp\cdot (1+p)^{-1}\cdot (1+\bar p)\ -\
\abs{q}^2 \, dp\cdot(1+p)^{-1} \Big ] (1+\bar p)^{-1}\\
+\ (1+p)^{-1}\, \Big [- (1+p)\cdot (1+\bar p)^{-1}\cdot d\bar p\ +\
\abs{q}^2 \,(1+p)^{-1}d\bar p \, \Big ] (1+\bar p)^{-1}\\
+(1+p)^{-1}\, \Big [dq\cdot \bar q\ -\ q\cdot d\bar q\, \Big ] (1+\bar p)^{-1}\ =\ \frac
{1}{\abs {1+p\, }^2}\, \lambda\, \tilde\eta\, \bar \lambda,
\end{multline}
\noindent where $\lambda\ = {\abs {1+p\,}}\, {(1+p)^{-1}}$ is a unit quaternion and
$\tilde\eta$ is the standard contact form on the sphere,
\begin{equation}\label{e:stand cont form on S}
\tilde\eta\ =\ dq\cdot \bar q\ +\ dp\cdot \bar p\ -\ q\cdot d\bar q -\ p\cdot d\bar p.
\end{equation}
\noindent Since $\abs{1+p}\ =\ \frac {2}{\abs {1+p'}}$ we have $\lambda\ =\ \frac
{1+p'}{\abs{1+p'\,}}$ equation \eqref{e:Cayley transf ctct form} can be put in the form
\[
\lambda\ \cdot (\mathcal{C}^{-1})^*\, \tilde\eta\ \cdot \bar\lambda\ =\ \frac
{8}{\abs{1+p'\, }^2}\, \tilde\Theta.
\]
We see that up to a constant multiplicative factor and a quaternionic contact
automorphism the forms $(\mathcal{C}^{-1})^*\tilde\eta$ and $\tilde\Theta$ are conformal
to each other. It follows that the same is true for $(\mathcal{C}^{-1})^*\eta$ and
$\tilde\Theta$.

In addition, $\tilde\Theta$ is qc-Einstein by definition, while
$\eta$ and hence also $(\mathcal{C}^{-1})^* \eta$ are qc-Einstein as
we observed at the beginning of the proof. Now we can apply Theorem
\ref{t:einstein preserving} according to which  up to a
multiplicative constant factor the forms
$(\mathcal{C}^{-1})^*\tilde\eta$ and $(\mathcal{C}^{-1})^*\eta$ are
related by a translation or dilation on the Heisenebrg group. Hence,
we conclude that up to a multiplicative constant, $\eta$ is obtained
from $\tilde\eta$ by a conformal quaternionic contact automorphism,
see Definition \ref{d:3-ctct auto}.
\end{proof}
Let us note that the Cayley transform defined in the setting of groups of Heisenberg type
is also a conformal transformation on $H$, see cf. \cite[Lemma 2.5]{ACD}. One can write
the above transformation formula in this more general setting.

\end{document}